\documentclass[11pt]{article}
\usepackage{amsmath,amsthm,amssymb,color,
cite
}
\usepackage[dvipdfmx]{graphicx
}
\usepackage{mathrsfs}

\renewcommand{\baselinestretch}{1.2}

\def\baselinestretch{1.4}
\newlength{\minitwocolumn}
\newcommand{\A}{{\mathbb A}}

\newcommand{\Z}{{\Bbb Z}} 
\newcommand{\C}{{\Bbb C}} 
\newcommand{\N}{{\Bbb N}} 
\newcommand{\PP}{{\Bbb P}} 

\newcommand{\F}{{\cal F}}

\newcommand{\cC}{{\cal C}}

\newcommand{\cG}{{\cal G}}

\newcommand{\cU}{{\cal U}}

\newcommand{\cN}{{\cal N}}

\newcommand{\cP}{{\cal P}}
\newcommand{\cE}{{\cal E}}

\newcommand{\gC}{\mathfrak{C}}

\newcommand{\cK}{{\cal K}}

\newcommand{\la}{\lambda}

\newcommand{\al}{\alpha}

\newcommand{\ep}{\epsilon}
\newcommand{\vep}{\varepsilon}

\newcommand{\tA}{{A}}

\newcommand{\bo}{{\bf 1}}

\newcommand{\nn}{{\nonumber}}
\newcommand{\bea}{\begin{eqnarray}}
\newcommand{\ena}{\end{eqnarray}}
\newcommand{\beit}{\begin{itemize}}
\newcommand{\enit}{\end{itemize}}

\newcommand{\be}{\begin{eqnarray*}}
\newcommand{\en}{\end{eqnarray*}}
\newcommand{\lb}[1]{\label{#1}}

\newcommand{\ds}[1]{{\displaystyle #1 }}

\newcommand{\id}{{\rm id}}

\newcommand{\Stab}{\mathrm{Stab}}

\newcommand{\bra}[1]{\langle #1 |}        
\newcommand{\ket}[1]{{| #1 \rangle}}      

\def\infq4p#1{{(#1;q^4,p)_\infty}}

\newcommand{\tPsi}{\widetilde{\Psi}}
\newcommand{\tPhi}{\widetilde{\Phi}}

\newcommand{\Tphi}[3]{
\phi\left(#1; 
\mmatrix{#2\\[-2mm] #3\cr}
\right)}

\newcommand{\eTphi}[4]{
\phi^{(#1)}\left(#2; 
\mmatrix{#3\\[-2mm] #4\cr}
\right)}

\newcommand{\tot}{\widetilde{\otimes}}

\newcommand{\tr}{{\rm tr}}

\newcommand{\mmatrix}[1]{\begin{matrix} #1 \end{matrix}}

\font\teneufm=eufm10
\font\seveneufm=eufm7
\font\fiveeufm=eufm5
\newfam\eufmfam
\textfont\eufmfam=\teneufm
\scriptfont\eufmfam=\seveneufm
\scriptscriptfont\eufmfam=\fiveeufm
\def\frak#1{{\fam\eufmfam\relax#1}}
\let\goth\mathfrak
\newcommand{\slth}{\widehat{\goth{sl}}_2}

\newcommand{\slnh}{\widehat{\goth{sl}}_N}
\newcommand{\sln}{\goth{sl}_N}
\newcommand{\g}{\goth{g}}

\newcommand{\Aqp}{{\cal A}_{q,p}}

\newcommand{\Bqla}{{{\cal B}_{q,\lambda}}}

\newcommand{\gl}{{\goth{gl}}}

\newcommand{\glnh}{\widehat{\goth{gl}}_N}

\newcommand{\gh}{\widehat{\goth{g}}}

\makeatletter
\@addtoreset{equation}{section}
\makeatother

\newtheorem{thm}{Theorem}[section]
\newtheorem{prop}[thm]{Proposition}
\newtheorem{lem}[thm]{Lemma}
\newtheorem{cor}[thm]{Corollary}

\newtheorem{df}[thm]{Definition}

\begin{document}

\vspace{-1cm}
\begin{center}
{\bf\Large  Elliptic Quantum Toroidal Algebra $U_{q,t,p}(\gl_{1,tor})$ \\and Affine Quiver Gauge Theories
\\[7mm] }
{\large  Hitoshi Konno${}^{\dagger}$ and Kazuyuki Oshima${}^{\star}$ }\\[6mm]
${}^\dagger${\it  Department of Mathematics, Tokyo University of Marine Science and 
Technology, \\Etchujima, Koto, Tokyo 135-8533, Japan\\
       hkonno0@kaiyodai.ac.jp}\\
${}^\star${\it Center for General Education, Aichi Institute of Technology, \\ Yakusa-cho, Toyota
470-0392, Japan }\\
oshima@aitech.ac.jp
\\[7mm]
\end{center}

\begin{center}
Dedicated to Professor Michio Jimbo on the occasion of his 70th birthday. 
\end{center}

\begin{abstract}
\noindent 
We introduce a new  elliptic quantum toroidal algebra $U_{q,t,p}(\gl_{1,tor})$. Various representations  in the quantum toroidal algebra $U_{q,t}(\gl_{1,tor})$ are extended to the elliptic case including the level (0,0) representation realized by using the elliptic Ruijsenaars difference operator. Various intertwining operators of $U_{q,t,p}(\gl_{1,tor})$-modules w.r.t. the Drinfeld comultiplication  are also constructed. 
We show that $U_{q,t,p}(\gl_{1,tor})$ gives a realization of the affine quiver $W$-algebra 
$W_{q,t}(\Gamma(\widehat{A}_0))$ proposed by Kimura-Pestun. This realization turns out to be useful to derive the Nekrasov instanton partition functions,  i.e. the $\chi_y$- and elliptic genus,  of the 5d and 6d lifts of the 4d ${\cal N}=2^*$ $U(M)$ theories
 and  provide  a new Alday-Gaiotto-Tachikawa correspondence.  
\end{abstract}
\nopagebreak

\section{Introduction}
The aim of this paper is to introduce a new elliptic quantum toroidal algebra $U_{q,t,p}(\gl_{1,tor})$
associated with the  toroidal algebra of type $\gl_1$ and  show its connection to the affine quiver $W$-algebra $W_{q,t}(\Gamma(\widehat{A}_0))$ 
 and the affine quiver gauge theories in 5d and 6d. 

We formulate $U_{q,t,p}(\gl_{1,tor})$ by generators and relations in the same scheme as the elliptic quantum group $U_{q,p}(\gh)$ associated with the affine Lie algebra $\gh$\cite{K98,JKOS,FKO,Konno16,KonnoBook}. 
The latter is the Drinfeld realization of the face type elliptic quantum group\cite{JKOStg} in 
terms of the Drinfeld generators, which are a deformation of the loop generators of the affine Lie algebras\cite{Drinfeld}. 
One of the important properties of the elliptic algebra $U_{q,p}(\gh)$ is 
that it gives a realization of  the deformation of the $W$ algebras of the Goddard-Kent-Olive (GKO) coset type 
$\gh_{r-h^\vee-1}\oplus \gh_1\supset \gh_{r-h^\vee}$\cite{GKO}.  Note that the GKO coset type $W$ algebras are isomorphic to those obtained by the quantum Hamiltonian reduction for simply laced $\gh$ but not for non-simply laced ones. See for example \cite{BS}. 
For simplicity let us consider the deformed $W$ algebra in  the simply laced case  and use the same notation $W_{q,t}(\g)$ as  the one  for the Hamiltonian reduction type\cite{FrRe}.  
Then the level-1 representation of  $U_{q,p}(\gh)$ realizes $W_{q,t}(\g)$\cite{K98,JKOS,KK03,FKO,Konno14} in the following way\footnote{The same was confirmed  for the non-simply laced case, say for $U_{q,p}(B^{(1)}_l)$\cite{FKO,Konno14}, which realizes the deformation of Fateev-Lukyanov's $W\hspace{-1mm}B_l$ algebra\cite{FaLu}.   
}.  
\begin{itemize}
\item[a)] the  parameters  $(p,p^*)$ of $U_{q,p}(\gh)$ are  identified with  $(q,t)$ of $W_{q,t}(\g)$, where $p^*=pq^{-2}$.  Here  $q$ appearing in the both sides are  different.   
\item[b)] the elliptic Drinfeld currents $E_j(z), F_j(z)$ of $U_{q,p}(\gh)$ realize the screening currents $S^+_j(z), S^-_j(z)$ of $W_{q,t}(\g)$.
\item[c)]  the basic generating functions $\Lambda_i(z)$ of $W_{q,t}(\g)$ are realized by 
  certain compositions of the intertwiners of the $U_{q,p}(\gh)$-modules w.r.t. the Drinfeld comultiplication.  
\end{itemize}
Furthermore for $W_{q,t}(\sln)$ a deformation of the primary fields constructed for example in \cite{AY} are realized as certain fusions of  intertwiners of the $U_{q,p}(\slnh)$-modules.       

The elliptic quantum toroidal algebra $U_{q,t,p}(\gl_{1,tor})$ is an elliptic quantum group associated with the toroidal algebra $\gl_1$. It is an elliptic analogue of the quantum  toroidal algebra $U_{q,t}(\gl_{1,tor})$ introduced by Miki  as a deformation of the $W_{1+\infty}$ algebra\cite{Miki}. 
Various representations of $U_{q,t}(\gl_{1,tor})$ have been studied by many 
papers such as \cite{AFS,Miki, FT, FHHSY, FFJMM, FJMM12, FJMM15, BFM}.  
It is also remarkable  that  $U_{q,t}(\gl_{1,tor})$ is  isomorphic to the elliptic Hall algebra introduced by Schiffmann and Vasserot\cite{SV1,SV2,SV3,Sc1,FFJMM}.  See also \cite{Negut,FHHSY,FT} 
for the shuffle algebra formulation. 

Representations of  $U_{q,t}(\gl_{1,tor})$  have many  interesting applications 
to the 5d  lifts of the 4d ${\cN=2}$ SUSY gauge theories, which are the gauge theories  associated with the linear quivers, 
such as a calculation of the Nekrasov instanton partition functions and a study of  
the 5d analogue of the Alday-Gaiotto-Tachikawa (AGT) correspondence\cite{AGT}. 
See for example \cite{AFS,MMZ,Zenke,AKMMMMOZ,BFMZZ,BFHMZ}. 
The essential  properties of $U_{q,t}(\gl_{1,tor})$ for applications to show AGT correspondence are summarized into the three points:
\begin{itemize}
\item[1)] the deformed $W$-algebra $W_{q,t}(\sln)$ is realized on the $N$-fold tensor product  of  $U_{q,t}(\gl_{1,tor})$\cite{ Miki, FHHSY, BFM, AKMMMMOZ}. 
\item[2)] the intertwiners of  $U_{q,t}(\gl_{1,tor})$-modules w.r.t the Drinfeld comultiplication realize the refined topological vertex \cite{IKV, AK} as the vertex operators on the bosonic Fock space\cite{AFS}. 
\item[3)] a certain block of composition of the intertwiners of  $U_{q,t}(\gl_{1,tor})$-modules  realizes the intertwiner (screened vertex operator) of the elliptic quantum group $U_{q,p}(\slnh)$\cite{Zenke, FOS19, FOS20}.  
\end{itemize}

Hence by 2) taking  a composition of  the intertwiners of $U_{q,t}(\gl_{1,tor})$-modules  according to the web diagram\cite{IKV,AK,AFS}, the expectation value of it gives an instanton partition function of the 5d lift of the 4d $\cN=2$ theory and basically by 3) the same quantity is identified with  an expectation value of the intertwiners of $U_{q,p}(\slnh)$ i.e. 
a deformation of the conformal block.  This is the 5d analogue of the AGT correspondence. 
It is also remarkable that thus obtained deformed conformal blocks coincides with the vertex functions introduced by Okounkov\cite{Ok,AFO,KPSZ}. See for example (6.4) in \cite{Konno17}   
and take its trigonometric limit $q^{\kappa}\to 0$.   
There are some different approaches based on the elliptic Hall algebra \cite{SV3} and the shuffle algebra\cite{NegutAGT}. For the affine Yangian case see \cite{MO}. See also \cite{Na16}. 

In stead of taking an expectation value, one can take a trace of the composed intertwininig 
operators of $U_{q,t}(\gl_{1,tor})$\cite{Saito,Zhu,FZ}. Then it gives an instanton partition function of the 6d lift of the theory\cite{Nieri,IKY,FZ} and the same arguments as 1) - 3) yield that it is equivalent to the trace of the 
 intertwiners of $U_{q,p}(\slnh)$  i.e.  the trace of the vertex operator of $W_{q,t}(\sln)$\cite{FOS20}. Note that the integral expressions for the latter is given as a part of the integral solution to the elliptic $q$-KZ equation\cite{Konno17,KonnoBook}, which has exactly the same structure as discussed as quantum $q$-Langlands correspondence in \cite{AFO}. 
 This is the 6d analogue of the AGT correspondence. 
%
In this way, the instanton calculus and the AGT correspondence in 5d and 6d lifts of the 4d $\cN=2$ SUSY gauge theories  associated with the linear quivers can be treated by using the quantum toroidal algebra $U_{q,t}(\gl_{1,tor})$.

 In this paper, we show that the elliptic quantum toroidal algebra $U_{q,t,p}(\gl_{1,tor})$ provides  a relevant quantum group structure to treat 
 the 5d and 6d lifts of the 4d $\cN=2^*$ SUSY gauge theories, which are the 
  $\cN=2$ SUSY gauge theories coupled with the adjoint matter\cite{Nekrasov04,Nekrasov}. 
 These are gauge theories associated with the Jordan quiver.  
We also formulate the AGT correspondence between these  5d and 6d lifts 
 and the affine quiver $W$ algebra $W_{q,t}(\Gamma(\widehat{A}_0))$.  For this purpose we 
 construct some useful  representations  of  $U_{q,t,p}(\gl_{1,tor})$ and give a realization of  $W_{q,t}(\Gamma(\widehat{A}_0))$\cite{KimPes} in the same way as the above b) and c) for $U_{q,p}(\gh)$.  
Namely the  level  $(1,N)$ elliptic currents $x^\pm(z)$ of  $U_{q,t,p}(\gl_{1,tor})$ give the screening currents $S^\pm(z)$ 
of $W_{q,t}(\Gamma(\widehat{A}_0))$, and 
 a certain  composition of  the ``intertwiners'' of  $U_{q,t,p}(\gl_{1,tor})$\footnote{A composition of 
the type I intertwiner $\Phi(u)$ of the $U_{q,t,p}(\gl_{1,tor})$-module and its shifted inverse $\Phi^*(u)$. See $\S$\ref{WA0}.} gives a realization of  the generating function $T(u)$.  
  There  one of the key properties of  $U_{q,t,p}(\gl_{1,tor})$ is  that the level  $(1,N)$ representation   possesses  the  four parameters $q, t, p, p^*$ satisfying  $p/p^*=q/t$. 
They play the role of  the $SU(4)$ $\Omega$ deformation parameters introduced by Nekrasov\cite{Nekrasov}, i.e.  two of three independent parameters are  the 5d analogue the  $\Omega$ deformation parameters $\ep_1, \ep_2$ \cite{Nekrasov04} and the remaining one is a  deformation of the adjoint mass parameter in the 4d  $\cN=2^*$ SUSY gauge theories.   

Realizing the generating function $T(u)$ of $W_{q,t}(\Gamma(\widehat{A}_0))$, it becomes  trivial  that the vacuum expectation value of $T(u)$ gives the Nekrasov instanton partition function of the 5d lift of the 4d 
$\cN=2^*$ $U(1)$ theory, i.e. the generating function of the $\chi_y$ genus of the Hilbert scheme $\mathrm{ Hilb}_n(\C^2)$ of $n$ points on $\C^2$ \cite{KimPes,Kimura}, where we take $y=p$. 
However this fact and our realization of $T(u)$ as a composition of the ``intertwiners'' 
lead us to identify $T(u)$ with a basic refined topological vertex operator corresponding to the  diagram Fig.\ref{3-3vertex} studied in \cite{IKV,HIV,HIKLV}. Then such realization 
makes calculation of compositions of $T(u)$ easy and allows us to apply them to various instanton calculus in the 5d and 6d  lifts of the 4d $\cN=2^*$ theories. This includes $\cN=2^*$ $U(M)$ theories, 
whose instanton partition function is the generating function of the $\chi_y$ genus of the moduli space of the rank $M$ instantons with charge $n$.  Their 6d lifts give the generating functions of the elliptic genus of the same moduli space. Hence we establish a new AGT correspondence 
between the 5d and 6d lifts of the 4d $\cN=2^*$ SUSY gauge theories and the affine quiver $W$ algebra $W_{q,t}(\Gamma(\widehat{A}_0))$. 

This paper is organized as follows. 
In $\S$2,  after reviewing basic facts in  the quantum toroidal algebra $U_{q,t}(\gl_{1,tor})$, we give a definition of the elliptic quantum toroidal algebra $U_{q,t,p}(\gl_{1,tor})$. The $Z$-algebra structure,  an elliptic analogue of Miki's automorphism and Hopf algebroid structure are also given. In $\S$3 we construct the level $(1,N)$, $(0,1)$  and $(0,0)$ representations. 
In particular, the level $(0,0)$ representation is given by using the elliptic Ruijsenaars difference operator on $\C[[x_1,\cdots,x_N]]$. In $\S$4 we construct the type I and the type II dual vertex operator as intertwining operators of the $U_{q,t,p}(\gl_{1,tor})$-modules constructed in $\S$3. 
We also give a definition of their shifted inverse operators. 
In $\S$5 we give a realization of the affine quiver $W$ algebra $W_{q,t}(\Gamma(\widehat{A}_0))$ by using the level $(1,N)$ representation of $U_{q,t,p}(\gl_{1,tor})$. In $\S$6 several calculations of the instanton partition functions of the 5d and 6d lifts of the 4d $\cN=2^*$ SUSY gauge theories are given. The results indicate a new AGT correspondence between those gauge theories and $W_{q,t}(\Gamma(\widehat{A}_0))$. 
In $\S$7 we summarize calculations of basic correlation functions of our intertwining operators, which should provide the $\cN=2^*$ analogues of the instanton partition functions of the 5d and 6d lifts of the 4d $\cN=2$ pure $SU(N)$ theory and  $SU(N)$ theory coupled with $2N$ fundamental matters obtained, for example, in \cite{AK,AFS}. 
Appendix A is a list of formulas, which is used to discuss the $Z$-algebra structure in $\S$\ref{Sec:Zalg}. Appendix B is a direct check of the statement on the level $(0,1)$ representation.
In Appendix C and D, proofs of Theorem \ref{typeIVO} and \ref{typeIdual} are given.  
In Appendix \ref{NekF} some useful formulas of the Nekrasov function are collected. 

A partial results in $\S$2, 3, 4, 5 have been presented by H.K. at several workshops \cite{K19}. 
 
%
%

\section{Elliptic Quantum Toroidal Algebra $\mathcal{U}_{q,t,p}(\gl_{1,tor})$}
After reviewing the quantum toroidal algebra ${U}_{q,t}({\gl}_{1,tor})$, we introduce 
the elliptic quantum toriodal algebra  $\mathcal{U}_{q,t,p}(\gl_{1,tor})$. 
The $Z$-algebra structure,  an elliptic analogue of Miki's automorphism and Hopf algebroid structure are also given.

\subsection{The quantum toroidal algebra ${U}_{q,t}({\gl}_{1,tor})$ }
This section is a review  of Miki's results\cite{Miki} in our notation.
\subsubsection{Definition}
\begin{df}
Let $q, t$ be nonzero complex numbers such that $q, t, q/t$ are not roots of unity, and set
\be
&&\kappa_m=(1-q^{m})(1-t^{-m})(1-(t/q)^{m}),\\
&&G^\pm(z)=(1-q^{\pm1}z)(1-t^{\mp1}z)(1-(q/t)^{\mp1}z). 
\en
The quantum toroidal algebra $\mathcal{U}_{q,t}=U_{q,t}(\gl_{1,tor})$ is a  $\C$-algebra generated by
\be
a_m, \quad X_{n}^{\pm}, \quad \gamma^{\pm 1/2},\quad (\psi_0^{\pm})^{\pm1}
\quad
( m \in \mathbb{Z}\setminus \{ 0 \}, \, n \in \mathbb{Z}). 
\en
By using the generating functions
\be
&&X^\pm(z)=\sum_{n\in \Z}X_{n}^\pm z^{-n},\quad\\
&&\phi(z)= \psi_0^+\exp\left\{\sum_{m > 0} a_m\gamma^{m/2}z^{-m}\right\},\quad\\
&&\psi(z)= \psi_0^-\exp\left\{-\sum_{m > 0} a_{-m} \gamma^{m/2}z^{m}\right\},
\en
the defining  relations  are 
\bea
&& \gamma^{\pm 1/2} : \ \mbox{central}, \\
&& [a_m, a_n]=-\delta_{m+n,0} \frac{\kappa_m}{m}(\gamma^m-\gamma^{-m})\gamma^{-|m|}, \lb{amam}\\
&& [a_m,X^{+}(w)]=-\frac{\kappa_m }{m}\gamma^{-|m|-m/2 }w^m X^+(w), \\
&&[a_m,X^{-}(w)]=\frac{\kappa_m }{m} \gamma^{-m/2 }w^m X^-(w),\lb{amXm}\\
&& [X^+(z),X^-(w)]=\frac{(1-q)(1-1/t)}{(1-q/t)}
\left(\delta(\gamma^{-1}z/w) \phi(w) -\delta(\gamma z/w) \psi(\gamma^{-1}w)\right), \lb{relxpxm}\\
&&z^3G^\pm(w/z)X^\pm(z)X^\pm(w)=-w^3G^\pm(z/w)X^\pm(w)X^\pm(z), \\
&& {\rm Sym}_{z_1, z_2, z_3}z_2z_3^{-1}[X^{\pm}(z_1),[X^{\pm}(z_2),X^{\pm}(z_3)]]=0,
\ena
where $ \delta(z)=\sum_{n \in \mathbb{Z}}z^n $. 
\end{df}
\noindent
{\it Remark.}\ By rescaling the generators $X^\pm_m$ as
\be
&&X^{'+}(z)={(1-t/q)}X^+(z),\qquad X^{'-}(z)={(1-q/t)}X^-(z), 
\en
the relation \eqref{relxpxm} can be rewritten as
\be
&& [X^{'+}(z),X^{'-}(w)]={\kappa_1}
\left(\delta(\gamma^{-1}z/w) \phi(w) -\delta(\gamma z/w) \psi(\gamma^{-1}w)\right)
\en
so that the whole relations of $\mathcal{U}_{q,t}$ are symmetric 
under any permutations among $q_1=q, q_2=t^{-1}, q_3=t/q$.   
Note also that \eqref{amam}-\eqref{amXm} are equivalent to
\bea
&& \phi(z)\phi(w)=\phi(w)\phi(z), \qquad
 \psi(z)\psi(w)=\psi(w)\psi(z), \label{uqg1psps}\\
&& \frac{G^+(\gamma^{-1} w/z)}{G^+(\gamma w/z)}\phi(z)\psi(w)=\frac{G^+(\gamma z/w)}{G^+(\gamma^{-1} z/w)}\psi(w)\phi(z), \label{uqg2}\\
&& z^3{G^+(\gamma^{-1} w/z)}\phi(z)X^+(w)=-\gamma^{-3}w^3G^+(\gamma z/w)^{-1}X^+(w)\phi(z),\quad\\
&&z^3{G^+(w/z)}\psi(z)X^+(w)=-w^3 G^+(z/w)X^+(w)\psi(z), \label{uqg3}\\
&&
z^{-3}{G^+(w/z)^{-1}}\phi(z)X^-(w)=-w^{-3}G^+(z/w)^{-1}X^-(w)\phi(z), \quad\\
&&z^{-3}{G^+(\gamma^{-1}w/z)^{-1}}\psi(z)X^-(w)=-\gamma^3 w^{-3}G^+(\gamma z/w)^{-1}X^-(w)\psi(z). \label{uqg4}
\ena

\subsubsection{Automorphism of $\cU_{q,t}$}\lb{MikiAuto}
Let us set $\omega=(1-q)(1-t^{-1})$ and define $Y^\pm_l\in \cU_{q,t}\ (l\in \Z)$ by
\be
&&Y^+_l=\left\{\mmatrix{(\psi^+_0)^l({\rm ad}\, X^+_0)^{l-1}X^+_{-1}/(-\omega)^{l-1}&(l>0)\cr 
\mbox{$\ds{ -\frac{\gamma^{1/2}}{1-t/q}a_{-1}}$}&(l=0)\cr
-{(q/t)^{|l|}}\gamma ({\rm ad}\, X^-_0)^{|l|-1}X^-_{-1}/(-\omega)^{|l|-1}&(l<0)\cr}
\right.\\
&&Y^-_l=\left\{\mmatrix{
-(t/q)\gamma^{-1} ({\rm ad}\, X^+_0)^{l-1}X^+_{1}/\omega^{l-1}&(l>0)\cr
\mbox{$\ds{\frac{\gamma^{1/2}}{1-q/t}a_{1}}$}&(l=0)\cr
{(q/t)^{|l|-1}}(\psi^+_0)^l({\rm ad}\, X^-_0)^{|l|-1}X^-_{1}/\omega^{|l|-1}&(l<0)\cr 
}
\right.,
\en
where $(ad\, x)y=[x,y]$ for $x,y\in \cU_{q,t}$. In addition, for $k\in \Z\backslash\{0\}$ we set 
\be
&&h_k=\left\{\mmatrix{[X^+_{-1},\overbrace{X^+_0,\cdots, X^+_0}^{k-2},X^+_1]&(k\geq 2)\cr
X^+_0&(k=1)\cr
-(q/t)X^-_0&(k=-1)\cr
-(q/t)^{|k|}[X^-_{1},\underbrace{X^-_0,\cdots, X^-_0}_{|k|-2},X^-_{-1}]&(k\leq -2)\cr
} \right.,
\en
where for $x_1, x_2,\cdots,x_n\in \cU_{q,t}$ 
\be
&&[x_n,\cdots,x_1]=[x_n,[x_{n-1},\cdots,[x_2,x_1]\cdots]].
\en
In terms of $h_k$, we deine $b_k\in \cU_{q,t}\ (k\in \Z\backslash\{0\})$ by 
the relations
\be
&&\exp\left(\sum_{k>0}b_k(\psi^+_0)^{k/2}z^{-k}\right)=1+\kappa_1\sum_{k\geq 1}\frac{h_k}{\omega^k}z^{-k},\\
&&\exp\left(-\sum_{k>0}b_{-k}(\psi^+_0)^{k/2}z^{k}\right)=1+\kappa_1\sum_{k\geq 1}\frac{h_{-k}}{\omega^k}z^{k},
\en
where we set $\omega=(1-q)(1-t^{-1})$. Then we have
\begin{thm}[Miki]\lb{Thm:Miki}
The map 
\be
\Psi&:&a_1\mapsto (1-t/q)\gamma^{-1/2}X^+_0,\quad 
a_{-1}\mapsto -(1-q/t)\gamma^{-1/2}X^-_0,\quad 
X^+_0\mapsto -\frac{\gamma^{1/2}}{1-t/q}a_{-1},\quad 
X^-_0\mapsto \frac{\gamma^{1/2}}{1-q/t}a_{1},\quad \\
&&\gamma\mapsto \psi^+_0,\quad \psi^+_0\mapsto \gamma^{-1}
\en
gives an automorphism of $\cU_{q,t}$ satisfying $\Psi^4=1$. 
Moreover $\Psi$ maps
\be
&&a_k\mapsto b_k,\quad X^\pm_l\mapsto Y^\pm_l, \quad 
b_k\mapsto a_{-k},\quad Y^\pm_l\mapsto X^\mp_{-l}
\en
\end{thm}

For readers convenience, we listed below the correspondence between  our notations and those used in \cite{Miki}. The symbols with $M$ denote the ones in \cite{Miki}.  
\be
&&\gamma_M=q^{1/2},\qquad q_M=(t/q)^{1/2},\qquad C^{\pm 1}_M=\gamma^{\pm 1},\qquad \cC^{\pm1}_M=\psi_0^\pm,
\\
&&(q_M-q_M^{-1})(\gamma_M^k-\gamma_M^{-k})a^M_k=\gamma^{|k|/2}a_k,\qquad
(q_M-q_M^{-1})(\gamma_M^k-\gamma_M^{-k})b^M_k=(\psi^+_0)^{|k|/2}b_k,\\
&&h^M_k=\left(\frac{t^{1/2}}{1-q}\right)^{|k|}h_k,\\
&&X^{M+}_l=\frac{t^{1/2}}{1-q}X^+_l,\quad X^{M-}_l=-\frac{t^{1/2}q/t}{1-q}X^-_l,\qquad
Y^{M+}_l=\frac{t^{1/2}}{1-q}Y^+_l,\quad Y^{M-}_l=-\frac{t^{1/2}q/t}{1-q}Y^-_l.
\en

\subsection{The elliptic quantum toroidal algebra $U_{q,t,p}({\gl}_{1,tor})$ }
Let $p$ be an indeterminate and set 
\be
&&(z;p)_n=(1-z)(1-zp)\cdots (1-zp^{n-1}), \quad (z;p)_0=1, 
\en
for $n\in \N$, $z\in \C$. 
We define $U_{q,t,p}({\gl}_{1,tor})$ by generators and relations as follows.

\begin{df}\lb{def:EQTAgl1}
The elliptic quantum toroidal algebra  $\mathcal{U}_{q,t,p}=U_{q,t,p}({\gl}_{1,tor})$ is a
$\C[[p]]$-algebra generated by $\al_m, x^\pm_n, 
$ $(m\in \Z\backslash\{0\}, n\in \Z)$
and invertible elements $\psi^{\pm}_0, \gamma^{1/2}$.  
The defining relations can be conveniently expressed in terms of the  generating functions,  
which we  call the elliptic currents, 
\be
&&x^{\pm}(z)=\sum_{n\in \Z}x^\pm_nz^{-n}, \\
&&\psi^{+}(z)=\psi^+_0\exp\left(-\sum_{m>0}\frac{p^m}{1-p^m}\al_{-m}(\gamma^{-1/2}z)^m\right)
\exp\left(\sum_{m>0}\frac{1}{1-p^m}\al_{m}(\gamma^{-1/2}z)^{-m}\right),\\
&& \psi^{-}(z)=\psi^{-}_0\exp\left(-\sum_{m>0}\frac{1}{1-p^m}\al_{-m}(\gamma^{1/2}z)^m\right)
\exp\left(\sum_{m>0}\frac{p^m}{1-p^m}\al_{m}(\gamma^{1/2}z)^{-m}\right).
\en
The defining relations are
\begin{align}
& \psi^\pm_0, \gamma^{1/2}\ :\ \mbox{central}, 
\\
&[\al_m,\al_n]=-\frac{\kappa_m}{m}(\gamma^m-\gamma^{-m})\gamma^{-m}\frac{1-p^m}{1-p^{*m}}\delta_{m+n,0}\lb{alal},\\
& [\al_m, x^+(z)] =-\frac{\kappa_m}{m}\frac{1-p^m}{1-p^{*m}}\gamma^{-3m/2} z^m x^+(z) \quad (m\not=0). \lb{euqg6}\\
& [\al_m, x^-(z)] =\frac{\kappa_m}{m}\gamma^{-m/2} z^m x^-(z) \quad (m\not=0), \lb{euqg7-2}\\
& [x^+(z),x^-(w)]=\frac{(1-q)(1-1/t)}{(1-q/t)}
\left(\delta(\gamma^{-1}z/w) \psi^+(w) -\delta(\gamma z/w) \psi^-(\gamma^{-1}w)\right), \lb{ellrelxpxm}\\
&z^3G^+(w/z)g(w/z;p^*)
 x^+(z)x^+(w)=-w^3G^+(z/w)g(z/w;p^*)
x^+(w)x^+(z), \lb{euqg9}\\
&z^{3}G^-( w/z)g(w/z;p)^{-1}
 x^-(z)x^-(w)=-w^{3}G^-(z/w)g(z/w;p)^{-1}
 x^-(w)x^-(z), \lb{euqg10}\\
 & g(w/z;p^*)g(u/w;p^*)g(u/z;p^*)
\left(\frac{w}{u}+\frac{w}{z}-\frac{z}{w}-\frac{u}{w}\right)
x^+(z)x^+(w)x^+(u) \nn \\
&+ \mbox{permutations in $z,w,u$}=0, \lb{euqg11}\\
& g(w/z;p)^{-1}
g(u/w;p)^{-1}
g(u/z;p)^{-1}
\left(\frac{w}{u}+\frac{w}{z}-\frac{z}{w}-\frac{u}{w}\right)
x^-(z)x^-(w)x^-(u) \nn \\
&+ \mbox{permutations in $z,w,u$}=0, \lb{euqg12}
\end{align}
where we set $p^*=p \gamma^{-2}$ and 
\bea
&&g(z;s)=\exp\left(\sum_{m>0}\frac{\kappa_m}{m}\frac{s^m}{1-s^m}z^m\right)\in \C[[z]]
\lb{gpm}
\ena
for $s=p, p^*$.
\end{df}
We treat these relations 
 as formal Laurent series in $z, w$ and $u$. 
 All the coefficients in $z, w, u$ are well defined in the $p$-adic topology\cite{FIJKMY,Konno16}.  

It is sometimes convenient to set 
\be
&&\al'_m=\frac{1-p^{*m}}{1-p^{m}}\gamma^{m}\al_m \qquad (m\in \Z_{\not=0})
\en
which satisfy 
\bea
&&[\al'_m,\al'_n]=-\frac{\kappa_m}{m}(\gamma^m-\gamma^{-m})\gamma^{m}\frac{1-p^{*m}}{1-p^{m}}\delta_{m+n,0}.\lb{alpalp}
\ena

\noindent
{\it Remark.}\ 
In the same way as $\cU_{q,t}$, by rescaling the generators $x^\pm_m$ as
\be
&&x^{'+}(z)={(1-t/q)}x^+(z),\qquad x^{'-}(z)={(1-q/t)}x^-(z), 
\en
the relation \eqref{ellrelxpxm} can be rewritten as
\be
&& [x^{'+}(z),x^{'-}(w)]={\kappa_1}
\left(\delta(\gamma^{-1}z/w) \psi^+(w) -\delta(\gamma z/w) \psi^-(\gamma^{-1}w)\right)
\en
so that the whole relations of $\mathcal{U}_{q,t,p}$ are symmetric 
under any permutations among $q_1=q, q_2=t^{-1}, q_3=t/q$.

\noindent
{\it Remark.}\ 
On $\cU_{q,t,p}$-modules, where the central element $\gamma^{1/2}$ takes a complex value, we 
 regard $p, p^*=p\gamma^{-2}$ as a generic complex number with $|p|<1, |p^*|<1$, 
and have 
\be
&&g(z;p)=\frac{(pq z;p)_\infty}{(pq^{-1} z;p)_\infty}
\frac{(pt^{-1} z;p)_\infty}{(pt z;p)_\infty}\frac{(p(t/q) z;p)_\infty}{(p(t/q)^{-1} z;p)_\infty}.
\en
Here we set
\be
&&(z;p)_\infty=\prod_{n=0}^\infty(1-zp^n)\qquad |z|<1. 
\en
From this one can rewrite \eqref{alal}-\eqref{euqg7-2}, \eqref{euqg9}-\eqref{euqg10} 
in the sense of analytic continuation as  follows.
\be
& \psi^\pm(z)\psi^\pm(w)=\frac{g_\theta(w/z;p^*)}{g_\theta(w/z;p)}\psi^\pm(w)\psi^\pm(z),\\
& \psi^+(z)\psi^-(w)=\frac{g_\theta(\gamma^{-1} {w}/{z};p^*)}{g_\theta(\gamma{w}/{z};p)} \psi^-(w)\psi^+(z), \lb{euqg1b}\\
& \psi^+(z)x^+(w)=g_\theta(\gamma^{-1}{w}/{z};p^*) x^+(w)\psi^+(z), \lb{2euqg2}\\
& \psi^-(z)x^+(w)=g_\theta({w}/z;p^*)x^+(w)\psi^-(z), \lb{2euqg3}\\
& \psi^+(z)x^-(w)=g_\theta({z}/w;p) x^-(w) \psi^+(z), \lb{2euqg4}\\
& \psi^-(z)x^-(w)=g_\theta(\gamma{z}/{w};p)x^-(w)\psi^-(z), \lb{2euqg5}\\
& x^+(z)x^+(w)=g_\theta(w/{z};p^*)x^+(w)x^+(z), \lb{2euqg9}\\
& x^-(z)x^-(w)=g_\theta({z}/{w};p) x^-(w)x^-(z), \lb{2euqg10}
\en
where 
\bea
&&g_\theta(z;p)=\frac{\theta_p(q^{-1}z)\theta_p((q/t)z)\theta_p(t  z)}
{\theta_p(q z)\theta_p((q/t)^{-1} z)\theta_p(t^{-1} z)},\lb{fncg}\\
&&\theta_p(z)=(z;p)_\infty(p/z;p)_\infty,\\
&&g_\theta(z;p^*)=g_\theta(z;p)\bigl|_{p\mapsto p^*}.
\ena
Note that 
\be
&&g_\theta(z^{-1};p)=g_\theta(z;p)^{-1},\qquad g_\theta(pz;p)=g_\theta(z;p). 
\en

\subsection{The $Z$ algebra structure}\lb{Sec:Zalg}
We next show that the $Z$ algebra structure of $\cU_{q,t,p}$ remains the same as the one of 
$\cU_{q,t}$, i.e. it is not elliptically deformed. This is a common feature in  the elliptic algebras $U_{q,p}(\gh)$\cite{FKO}. 

In this subsection we assume $\gamma\not=1$. 
Set
\begin{align}
& E^{\pm}(\alpha,z)=\exp\left\{\pm \sum_{n>0}\frac{1}{\gamma^n-\gamma^{-n}}\alpha_{\pm n} (\gamma^{-1/2}z)^{\mp n} \right\},\lb{Epm}\\
& E^{\pm}(\alpha',z)=\exp\left\{\mp \sum_{n>0}\frac{1}{\gamma^n-\gamma^{-n}}\alpha'_{\pm n} (\gamma^{-1/2} z)^{\mp n} \right\}, \lb{Epmp}
\end{align}
and define $\mathcal{Z}^{\pm}(z)$ by
\begin{align}
& \mathcal{Z}^+(z) = E^-(\alpha,z)x^+(z)E^+(\alpha,z), \\
& \mathcal{Z}^-(z) = E^-(\alpha',z)x^-(z)E^+(\alpha',z).
\end{align}
Then from \eqref{alEp}-\eqref{alEm} in Appendix \ref{Zalg}, we have
\be
&& [\alpha_m, \mathcal{Z}^{\pm}(z)]=0, \quad (m \in \mathbb{Z}_{\ne 0}). \lb{zalg1}
\en
Furthermore  the relations in Definition \ref{def:EQTAgl1} and Lemma \ref{lem:sec3-2} leads to the following theorem.
\begin{thm}\lb{prop:sec3-2}
\begin{align}
& z^3G^-(w/z)g(w/z;\gamma^2)^{-1}\mathcal{Z}^\pm(z) \mathcal{Z}^\pm(w)=
-w^3G^-(z/w)g(z/w;\gamma^2)^{-1}\mathcal{Z}^\pm(w) \mathcal{Z}^\pm(z), \lb{zalg2}\\
& g(\gamma w/z)g(\gamma w/z;\gamma^2) \mathcal{Z}^+(z) \mathcal{Z}^-(w) 
-g(\gamma z/w)g(\gamma z/w;\gamma^2) \mathcal{Z}^-(w)\mathcal{Z}^+(z) \nn \\
& \qquad \qquad =\frac{(1-q)(1-1/t)}{1-q/t}
\left\{ \delta(\gamma^{-1} z/w) \psi_0^+ -\delta(\gamma z/w) \psi_0^- \right\}, \lb{zalg4}\\
& g(w/z)^{-1}g(u/w)^{-1}g(u/z)^{-1}g(w/z;\gamma^2)^{-1}\, g({u}/{w};\gamma^2)^{-1}\, 
g({u}/{z};\gamma^2)^{-1}\nn\\
&\qquad \times\left(\frac{w}{u}+\frac{w}{z}-\frac{z}{w}-\frac{u}{w} \right) 
\mathcal{Z}^+(z)\mathcal{Z}^+(w)\mathcal{Z}^+(u) 
+ \mbox{permutations in $z,w,u$}=0, \lb{zalg5}\\
& g(w/z;\gamma^2)^{-1}\, g({u}/{w};\gamma^2)^{-1}\, g({u}/{z};\gamma^2)^{-1}\left(\frac{w}{u}+\frac{w}{z}-\frac{z}{w}-\frac{u}{w} \right) 
\mathcal{Z}^-(z)\mathcal{Z}^-(w)\mathcal{Z}^-(u) \nn \\
&+ \mbox{permutations in $z,w,u$}=0. \lb{zalg6}
\end{align}
Here $g(z;\gamma^2)$ is given in \eqref{gpm} with $s=\gamma^2$ and we also set 
\bea
&&g(z)=\exp\left(\sum_{m>0}\frac{\kappa_m}{m}z^m\right)\quad \in \C[[z]].\lb{def:g}
\ena
\end{thm}
We call the algebra generated by the coefficients of $\mathcal{Z}^\pm(z)=\sum_{\in n\in \Z}
\mathcal{Z}^\pm_n z^{-n}$ in $z$ satisfying the relations in Theorem \ref{prop:sec3-2} the $Z$ algebra associated with $\cU_{q,t,p}$.

\noindent
{\it Remark.}\ Noting that the whole relations in this theorem are independent from the elliptic nome $p$, 
one finds that these relations coincides with those in Proposition 4.2 in \cite{Miki} under the identification 
\be
&&\cK_M=\gamma,\quad K^{-1}_M=\psi_0,\\
&&T^\pm(z)_M=\pm Z^\pm(z)/((1-q^{\pm 1})(1-t^{\mp1}))
\en 
and the correspondence listed below Theorem \ref{Thm:Miki}. 

Hence   
the elliptic quantum toroidal algebra $\cU_{q,t,p}$ possesses  the same deformed Virasoro and $W_3$ algebra structures through the same $Z$ algebra structure as $\cU_{q,t}$ 
 discussed in Theorem 4.1 and Remark 4.7 of \cite{Miki}.

\subsection{Algebra homomorphism from  $U_{q,t,p}(\gl_{1,tor})$ to $U_{q,t}(\gl_{1,tor})$}\lb{JKOShom}
In \cite{JKOS}, an homomorphism from the  elliptic quantum algebra $U_{q,p}(\gh)$  to the quantum affine algebra $U_q(\gh)$ 
with $\gh$ being an untwisted affine Lie algebra 
was given. This allows us to extend any representations of $U_q(\gh)$ with generic $q$ to those 
of $U_{q,p}(\gh)$ with generic $p, q$. 
 Such homomorphism  can  be  generalized to the cases with $\gh$ being any twisted 
 affine Lie algebras, for example see \cite{KK04}, as well as any toroidal Lie algebras $\g_{tor}$ \cite{FHHSY,K19}.  They include the  cases with $\gh$ being (affine or toroidal) Lie 
 super algebras. See for example \cite{Kojima}.  
 \begin{prop}\lb{evhom} 
 Let us define $u^\pm(z,p)$ by
 \be
 &&u^+(\gamma^{-1/2}z,p)=\exp\left(-\sum_{m>0} \frac{p^{*m}}{1-p^{*m}}a_{-m}z^m\right),\\
 &&u^-(\gamma^{1/2}z,p)=\exp\left(\sum_{m>0} \frac{p^{m}}{1-p^m}a_{m}z^{-m}\right),
 \en
 where $p^*=p\gamma^{-2}$ as above. 
 Then the following gives a homomorphism from $\cU_{q,t,p}$ to $\cU_{q,t}\otimes \C[[p]]$.
 \bea
 &&x^+(z)\mapsto u^+(z,p)X^+(z),\lb{upxp}\\
 &&x^-(z)\mapsto X^-(z)u^-(z,p),\lb{xmum}\\
 &&\psi^+(z)\mapsto u^+(\gamma z,p)\phi(z)u^-(z,p),\lb{uppum}\\
&&\psi^-(z)\mapsto u^+( z,p)\psi(z)u^-(\gamma z,p). \lb{uppsium}
\ena
In particular, \eqref{uppum} (or \eqref{uppsium}) is equivalent to 
\be
&&\al_m\mapsto a_m,\qquad \al_{-m}\mapsto \frac{1-p^m}{1-p^{*m}}a_{-m}\qquad (m>0).
\en
 \end{prop}
The statement is proved by using the following Lemma.  
\begin{lem}\lb{relupm}
\begin{align}
& [a_n,u^+(z,p)]=\frac{\kappa_n}{n}\frac{\gamma^n-\gamma^{-n}}{1-p^{*n}}(p^* \gamma^{-1/2} z)^n u^+(z,p) \quad (n >0), \\
&[a_{-n},u^-(z,p)]=\frac{\kappa_n}{n}\frac{\gamma^n-\gamma^{-n}}{1-p^n }(p^{-1} \gamma^{1/2} z)^{-n} u^-(z,p) \quad (n >0), \\
& u^+(z,p)X^+(w)=g( {z}/{w}; p^*)X^+(w)u^+(z,p),\\
& u^+(z,p)X^-(w)=g(\gamma {z}/{w};p^*)^{-1}X^-(w)u^+(z,p),\\
& u^-(z,p)X^+(w)=g(\gamma^{-1} {w}/{z}; p)^{-1}X^+(w)u^-(z,p),\\
& u^-(z,p)X^-(w)=g( {w}/{z};p)X^-(w)u^-(z,p),
\\
& \phi(z)u^+(w,p)={g(\gamma {w}/{z};p^*)}{g(\gamma^{-1} {w}/{z};p^*)^{-1}}u^+(w,p)\phi(z),\\
& u^-(z,p)\psi(w)={g(\gamma {w}/{z};p)}{g(\gamma^{-1} {w}/{z};p)^{-1}}\psi(w)u^-(z,p),
\\
& \phi(z)u^-(w,p)=u^-(w,p)\phi(z),\\
& u^+(z,p)\psi(w)=\psi(w)u^+(z,p),\\
& u^+(z,p)u^-(w,p)={g(p^*\gamma {z}/{w};p^*)}{g(p\gamma^{-1}z/w;p)^{-1}}u^-(w,p)u^+(z,p).
\end{align}
Here $g(z;s), \ s=p,\ p^*$ is given in \eqref{gpm}. 

\end{lem}

It is obvious that the homomorphism in Proposition \ref{evhom} is invertible. 
Hence one has the following isomorphism. 

\begin{thm}\lb{iso}
For generic $q, t, p$, 
\be
&&\cU_{q,t,p} \cong \cU_{q,t}\otimes \C[[p]]
\en
as a topological algebra. 
\end{thm}

\subsection{Elliptic analogue of Miki's automorphism}
By using a set of new generators $b_m, Y^\pm_l$ of $U_{q,t}(\gl_{1,tor})$ defined in $\S$ \ref{MikiAuto}, let us consider their generating functions.  
\be
&&Y^\pm(z)=\sum_{l\in \Z}Y^\pm_l z^{-l},\\
&&\phi_b(z)= \gamma^{-1}\exp\left\{\sum_{m > 0} b_m(\psi^+_0)^{m/2}z^{-m}\right\},\quad\\
&&\psi_b(z)= \gamma\exp\left\{-\sum_{m > 0} b_{-m} (\psi^+_0)^{m/2}z^{m}\right\}.
\en
Applying the homomorphism in $\S$\ref{JKOShom}, one obtains 
a new elliptic quantum toroidal algebra $U^b_{q,t,p}(\gl_{1,tor})$ generated by  
$\beta_m, y^\pm_l, \psi_0^+, \gamma^{-1}$, defined by
\be
&& y^\pm(z)=\sum_{l\in \Z}y^\pm_lz^{-l},\\
&&\psi_b^+((\psi_0^+)^{1/2}z)=\gamma^{-1}\exp\left(-\sum_{m>0}\frac{p^m}{1-p^m}\beta_{-m}z^m\right)\exp\left(\sum_{m>0}\frac{1}{1-p^m}\beta_{m}z^{-m}\right),
\en
where
\bea
&&y^+(z):=u_b^+(z,p)Y^+(z),\lb{upyp}\\
&&y^-(z):=Y^-(z)u_b^-(z,p),\lb{ymum}\\
&&\psi_b^+(z):=u^+_b(\psi^+_0 z,p)\phi_b(z)u^-_b(z,p),\lb{uppbum}\\
&&\psi_b^-(z):=u^+_b( z,p)\psi_b(z)u^-_b(\psi^+_0 z,p)\lb{uppsibum}
\ena
with 
\be
&&u^+_b((\psi^+_0)^{-1/2}z,p):=\exp\left(-\sum_{m>0} \frac{p_b^{*m}}{1-p_b^{*m}}b_{-m}z^m\right),\\
 &&u^-_b((\psi^+_0)^{1/2}z,p):=\exp\left(\sum_{m>0} \frac{p^{m}}{1-p^m}b_{m}z^{-m}\right),
\en
and $p^*_b=p(\psi^+_0)^{-2}$.

Then we obtain the following isomorphism. 
\begin{cor}
\be
&&U_{q,t,p}(\gl_{1,tor})\cong U^b_{q,t,p}(\gl_{1,tor}),
\en
by
\be
&&\al_m\mapsto \beta_m, \quad x^\pm_l \mapsto y^\pm_l,\quad  
\gamma \mapsto \psi^+_0,\quad \psi^+_0 \mapsto \gamma^{-1}.
\en
\end{cor}
\noindent
{\it Proof.}\ The statement follows from Theorem \ref{Thm:Miki} and \ref{iso}.
\qed

\subsection{Hopf algebroid structure}\lb{coalgstr}
In this section, we introduce a Hopf algebroid structure into $\cU_{q,t,p}$. 
This structure was introduced by Etingof and Varchenko \cite{EV} and developed in \cite{KR,Konno09}. 

For $F(z,p)\in \C[[z,z^{-1}]][[p]]$, 
 let  $\tot$ denote the usual tensor product with the following extra condition
\bea
&&F(z,p^*)a\tot b=a\tot F(z,p)b \lb{reltot}
\ena 
Here $F(z,p^*)$ denotes the same $F$ with replacing $p$ by $p^*=p\gamma^{-2}$. 

Define two moment maps $\mu_l, \mu_r : \C[[z,z^{-1}]][[p]] \to \cU_{q,t,p}[[z,z^{-1}]]$   by 
\be
&\mu_l(F(z,p))=F(z,p),\qquad \mu_r(F(z,p))=F(z,p^*).
\en
Let $\gamma_{(1)}=\gamma\tot 1, \gamma_{(2)}=1\tot \gamma$ and  $p^{*}_{(i)}=p\gamma^{-2}_{(i)}$ $(i=1,2)$. Let us define two algebra homomorphisms $\Delta: \cU_{q,t,p}\to \cU_{q,t,p}\tot \cU_{q,t,p}$  and
$\vep: \cU_{q,t,p}\to \C$ by
\bea
&&\Delta(\gamma^{\pm 1/2})=\gamma^{\pm 1/2}\tot\gamma^{\pm 1/2},\lb{co1}\\
&&\Delta(\psi^\pm(z))=\psi^\pm(\gamma^{\mp1/2}_{(2)}z)\tot \psi^\pm(\gamma^{\pm 1/2}_{(1)}z),\\
&&\Delta(x^+(z))=1\tot x^{+}(\gamma^{-1/2}_{(1)}z)+x^{+}(\gamma^{1/2}_{(2)}z)\tot \psi^-(\gamma^{- 1/2}_{(1)}z),\\
&&\Delta(x^-(z))= x^{-}(\gamma^{-1/2}_{(2)}z)\tot 1+\psi^+(\gamma^{- 1/2}_{(2)}z)\tot x^{-}(\gamma^{1/2}_{(1)}z),\lb{co4}\\
&&\Delta(\mu_l(F(z,p)))=\mu_l(F(z,p))\tot 1,\qquad \Delta(\mu_r(F(z,p)))=1\tot \mu_r(F(z,p)),\\
&&\vep(\gamma^{1/2})=\vep(\psi^+_0)=1,\qquad \vep(\psi^\pm(z))=1,\qquad \vep(x^\pm(z))=0,\\
&&\vep(\mu_l(F(z,p)))=\vep(\mu_r(F(z,p)))=F(z,p).
\ena
Note that $\Delta$ is so-called the Drinfeld comultiplication.  Then we have
\begin{prop}
The maps $\vep$ and $\Delta$ satisfy
\bea
&&(\Delta\tot \id)\circ \Delta=(\id \tot \Delta)\circ \Delta,\lb{coaso}\\
&&(\vep \tot \id)\circ\Delta =\id =(\id \tot \vep)\circ \Delta.\lb{vepDelta}
\ena
\end{prop}
We also define an algebra anti-homomorphism $S: \cU_{q,t,p}\to \cU_{q,t,p}$ by
\be
&&S(\gamma^{1/2})=\gamma^{-1/2},\\
&&S(\psi^\pm(z))=\psi^\pm(z)^{-1},\\
&&S(x^+(z))=-x^+(z)\psi^-(z)^{-1},\\
&&S(x^-(z))=-\psi^+(z)x^-(z),\\
&&S(\mu_l(F(z,p)))=\mu_r(F(z,p)),\qquad S(\mu_r(F(z,p)))=\mu_l(F(z,p)).
\en
Then one can check the following.
\begin{prop}
\be
&&m\circ(\id\tot S)\circ \Delta(a)=\mu_l(\vep(a)1)\qquad \forall a\in \cU_{q,t,p},\\
&&m\circ( S\tot\id)\circ \Delta(a)=\mu_r(\vep(a)1).
\en
\end{prop}
These Propositions indicate that $(\cU_{q,t,p},\Delta,\vep,\mu_l,\mu_r,S)$ is a Hopf algebroid\cite{EV,KR,Konno08,Konno09}. 

\noindent
{\it Remark.} For representations on which $\gamma^{1/2}=1$ hence $p=p^*$, the following gives an opposite Drinfeld comultiplication.
\bea
&&\Delta^{op}(\psi^\pm(z))=\psi^\pm(z)\tot \psi^\pm(z),\lb{opco1}\\
&&\Delta^{op}(x^+(z))=x^{+}(z)\tot  1+\psi^-(z) \tot x^{+}(z),\\
&&\Delta^{op}(x^-(z))= 1\tot x^{-}(z)+ x^{-}(z)\tot \psi^+(z).\lb{opco3}
\ena 

\section{Representations of $U_{q,t,p}(\gl_{1,tor})$}
Let ${\cal V}$ be a $\cU_{q,t,p}$ module. For $(k,l)\in \C^2$, we say that ${\cal V}$ has level 
$(k,l)$, if the central elements $\gamma^{1/2}$ and $\psi^+_0$ act as
\be
&&\gamma^{1/2} v=(t/q)^{k/4}v, \qquad \psi^+_0 v=(t/q)^{-l/2} v\qquad \forall v\in {\cal V}.
\en

In the rest of this paper, we regard $p, p^*=p\gamma^{-2}$ as a generic complex number
 with $|p|<1, |p^*|<1$. 

More generally, for  $p_1,\cdots p_r
\in \C$ with $|p_i|<1$ $(1\leq i\leq r,\  r=1,2,\cdots )$, we set
\be
&&(z;p_1,\cdots,p_r)_\infty=\prod_{n_1,\cdots,n_r=0}^\infty(1-zp_1^{n_1}\cdots p_r^{n_r}). 
\en
We use the following multiple elliptic Gamma functions defined by \cite{Nishizawa,Narukawa}
\be
&&\Gamma_r(z;p_1,\cdots,p_{r})=(z;p_1,\cdots, p_{r})_\infty^{(-1)^{r-1}}{(p_1\cdots p_{r}/z;p_1,\cdots,p_{r})_\infty}.
\en
Note 
\be
&&\Gamma_1(z;p_1)=\theta_{p_1}(z),\\
&&\Gamma_2(z;p_1,p_{2})\equiv\Gamma(z;p_1,p_2)=\frac{(p_1p_2/z;p_1,p_2)_\infty}{(z;p_1,p_2)_\infty}.
\en
are Jacobi's odd theta function and  the elliptic Gamma function\cite{Ruijs}, respectively. 
We have
\begin{prop}
\be
&&\Gamma_r(p_jz;p_1,\cdots,p_{r})=\Gamma_{r-1}(z;p_1,\cdots,\check{p}_j,\cdots,p_{r})\Gamma_r(x;p_1,\cdots,p_{r})
\en
for $r\geq 2$, where $\check{p}_j$ means the excluding of $p_j$. 
\end{prop}

\subsection{Level $(1,N)$ representation of  ${U}_{q,t,p}({\gl}_{1,tor})$}\lb{sec:level1N}
For $u\in \C^*$, let $\F^{(1,N)}_u=\C[\al_{-m}\ (m>0)]1^{(N)}_u$ be a Fock space on which the Heisenberg algebra $\{\al_m\ (m\in \Z_{\not=0}) \}$ and the central elements $\gamma^{1/2}, \psi^+_0$ act as
\be
&&\gamma^{1/2}\cdot 1^{(N)}_u=(t/q)^{1/4}1^{(N)}_u,\quad \psi^+_0\cdot 1^{(N)}_u=(t/q)^{-N/2}1^{(N)}_u,\quad \al_{-m}\cdot 1^{(N)}_u=0,\\
&&\al_{-m}\cdot \xi =  \al_{-m}\xi,\\
&&\al_m\cdot \xi =-\frac{\kappa_m}{m}(1-(q/t)^m)\frac{1-p^m}{1-p^{*m}}\frac{\partial}{\partial \al_{-m}}\xi
\en
for $m>0$, $\xi\in \F^{(1,N)}_u$. Note that $p^*=pq/t$ on $\F^{(1,N)}_u$. 
\begin{thm}\lb{level1N}
The following assignment gives a level $(1,N)$ representation of $\cU_{q,t,p}$ on $\F^{(1,N)}_u$. 
\begin{align}
& x^+(z)=u z^{-N}(t/q)^{3N/4}\exp \left\{ -\sum_{n>0}\frac{(t/q)^{n/4}}{1-(t/q)^n}\al_{-n}z^n  \right\}
\exp \left\{ \sum_{n>0}\frac{(t/q)^{3n/4}}{1-(t/q)^n}\al_{n}z^{-n} \right\}, \\
& x^-(z)=u^{-1}z^N(t/q)^{-3N/4}\exp \left\{ \sum_{n>0}\frac{(t/q)^{n/4}}{1-(t/q)^n}\al'_{-n}z^n \right\}
\exp \left\{ -\sum_{n>0}\frac{(t/q)^{3n/4}}{1-(t/q)^n}\al'_{n}z^{-n} \right\}, \\
& \psi^+(z)=(t/q)^{-N/2}\exp \left\{-\sum_{n>0}\frac{p^n(t/ q)^{-n/4}}{1-p^{n}}\al_{-n} z^n \right\}
\exp \left\{ \sum_{n >0} \frac{(t/q)^{n/4}}{1-p^n} \al_n z^{-n} \right\}, \\
& \psi^-(z)=(t/q)^{N/2}\exp \left\{-\sum_{n>0} \frac{(t/q)^{n/4}}{1-p^{n}} \al_{-n} z^n \right\}
\exp \left\{ \sum_{n>0} \frac{p^{n}(t/q)^{-n/4}}{1-p^n} \al_n z^{-n} \right\}.
\end{align}
\end{thm}

\noindent
{\it Proof.} \ 
Let us set $x^\pm(z)=u^\pm z^{\mp N}(t/q)^{\pm 3N/4}\widetilde{x}^\pm(z)$, $\psi^\pm(z)=(t/q)^{\mp N/2}\widetilde{\psi}^\pm(z)$. 
The statement follows from the operator product expansion (OPE) formulas listed below. 
\be
&&\widetilde{\psi}^+(z)x^+((t/q)^{1/2}w)=f(w/z;p^*):\widetilde{\psi}^+(z)x^+((t/q)^{1/2}w):\qquad |w/z|<1,\nn\\
&&x^+((t/q)^{1/2}w)\widetilde{\psi}^+(z)=f(p^*z/w;p^*):\widetilde{\psi}^+(z)x^+((t/q)^{1/2}w):\qquad |z/w|<1,\nn\\
&&\widetilde{\psi}^-(z)x^+((t/q)^{1/2}w)=f(p^*(t/q)^{1/2}w/z;p^*):\widetilde{\psi}^-(z)x^+((t/q)^{1/2}w):\qquad |w/z|<1,\nn\\
&&x^+((t/q)^{1/2}w)\widetilde{\psi}^-(z)=f((t/q)^{-1/2}z/w;p^*):\widetilde{\psi}^-(z)x^+((t/q)^{1/2}w):\qquad |z/w|<1,\nn\\
&&\widetilde{\psi}^-(z)x^-(w)=f(p(t/q)^{-1/2}w/z;p)^{-1}:\widetilde{\psi}^-(z)x^-(w):\qquad |w/z|<1,\nn\\
&&x^-(w)\widetilde{\psi}^-(z)=f((t/q)^{1/2}z/w;p)^{-1}:\widetilde{\psi}^-(z)x^-((t/q)^{1/2}w):\qquad |z/w|<1,\nn\\
&&x^+(z)x^-(w)=(w/z)^N\frac{(1-q\gamma w/z)(1- \gamma w/tz)}{(1-\gamma w/z)(1-q\gamma w/tz)}:\widetilde{x}^+(z)\widetilde{x}^-(w):\qquad |w/z|<1,\nn\\
&&x^-(w)x^+(z)=(w/z)^N\frac{(1-\gamma^{-1} z/qw)(1-t\gamma^{-1} z/w)}{(1-\gamma^{-1} z/w)(1-t\gamma^{-1} z/qw)}:\widetilde{x}^+(z)\widetilde{x}^-(w):\qquad |z/w|<1,\nn
\en
and the formulas
\bea
&&:\widetilde{x}^+(\gamma z)\widetilde{x}^-(w):=\widetilde{\psi}^+(z),\\
&&:\widetilde{x}^+(\gamma^{-1} z)\widetilde{x}^-(w):=\widetilde{\psi}^-(\gamma^{-1}z).
\ena
Here $f(z;s)$ is given by
\be
&&f(z;s)=\frac{(q^{-1}z;s)_\infty(tz;s)_\infty(qt^{-1}z;s)_\infty}{(qz;s)_\infty(t^{-1}z;s)_\infty(tq^{-1}z;s)_\infty}
\en
for $s=p, p^*$. 
\qed

A similar representation is  given in Lemma A.16  in \cite{FHHSY}.

\subsection{Vector representation and the $q$-Fock space representation
 }
 We next consider the elliptic analogue of a representation given in \cite{FT,FFJMM} , i.e.  the $q$-Fock space representation  of $\cU_{q,t,p}$. 
 
For $u\in \C^*$, let $V(u)$ be a vector space spanned by $[u]_j\ (j\in \Z)$. 
\begin{prop}\lb{vecrep}
By the following action, $V(u)$ is a level-$(0,0)$ $\cU_{q,t,p}$-module. We call $V(u)$ the vector representation.
\be
&&x^+(z)[u]_j=a^+(p)\delta(q^ju/z)[u]_{j+1},\\
&&x^-(z)[u]_j=a^-(p)\delta(q^{j-1}u/z)[u]_{j-1},\\   
&&\psi^\pm(z)[u]_j=
\frac{\theta_p(q^jt^{-1}u/z)\theta_p(q^{j-1}tu/z)}{\theta_p(q^ju/z)\theta_p(q^{j-1}u/z)}\Biggl|_\pm[u]_j,\\
&&\psi^\pm_0[u]_j=[u]_j,
\en
where we set
\bea
a^+(p)&=&(1-t)\frac{(pt/q;p)_\infty(p/t;p)_\infty}{(p;p)_\infty(p/q;p)_\infty},\lb{def:ap}\\
a^-(p)&=&(1-t^{-1})\frac{(pq/t;p)_\infty(pt;p)_\infty}{(p;p)_\infty(pq;p)_\infty}.\lb{def:am}
\ena
\end{prop}
\noindent
{\it Proof.}\ 
Use a formula in Lemma \ref{limDiffTheta}. 

\qed

Let ${\cal F}_u$  be a vector space spanned by $\ket{\la}_u$ $(\la\in {\cal P}^+)$, where 
\bea
&&{\cal P}^+=\{ \la=(\la_1,\la_2,\cdots)\ |\ \la_i\geq \la_{i+1},\ \la_i\in \Z,\  \la_l=0\ \mbox{for sufficiently large}\ l\ \} .\lb{partPp}
\ena
We denote by  $\ell(\la)$  the length of $\la\in {\cal P}^+$ i.e. $ \la_{\ell(\la)}>0$ and $\la_{\ell(\la)+1}=0$. We also set $|\la|=\sum_{i\geq 1}\la_i$ and denote by $\la'$ the conjugate of $\la$.

\begin{thm}\lb{actionUqtpgl1}
 The following action  gives a level (0,1) representation of $\cU_{q,t,p}$ on $\F_u$. 
 We denote this by $\F^{(0,1)}_u$. 
\bea
&&\gamma^{1/2}\ket{\la}_u=\ket{\la}_u,\\
&&x^+(z)\ket{\la}_u=a^+(p)\sum_{i=1}^{\ell(\la)+1}A^+_{\la,i}(p)\delta(
u_i/z)\ket{\la+{\bf 1}_i}_u,\lb{actxp}\\
&&x^-(z)\ket{\la}_u=(q/t)^{1/2}a^-(p)\sum_{i=1}^{\ell(\la)}A^-_{\la,i}(p)\delta(
q^{-1}u_i/z)\ket{\la-{\bf 1}_i}_u,\lb{actxm}\\
&&\psi^+(z)\ket{\la}_u=(q/t)^{1/2}B^+_{\la}(u/z;p)\ket{\la}_u,\\
&&\psi^-(z)\ket{\la}_u=(q/t)^{-1/2}B^-_{\la}(z/u;p)\ket{\la}_u,\quad 
\ena
where  we set 
\bea
A^+_{\la,i}(p)&=&
\prod_{j=1}^{i-1}\frac{\theta_p(tu_i/u_j)\theta_p(qt^{-1}u_i/u_j)}{\theta_p(qu_i/u_j)\theta_p(u_i/u_j)},\lb{tildeAp=NAp}\\
A^-_{\la,i}(p)&=&
 \prod_{j=i+1}^{\ell(\la)}\frac{\theta_p(qt^{-1}u_j/u_i)}{\theta_p(u_j/u_i)}
  \prod_{j=i+1}^{\ell(\la)+1}\frac{\theta_p(tu_j/u_i)}{\theta_p(qu_j/u_i)},
  \lb{tildeAm=NAm}\\
B^+_\la(u/z;p)&=&
\prod_{i=1}^{\ell(\la)}
\frac{\theta_p(u_i/tz)}{\theta_p(u_i/qz)}
\prod_{i=1}^{\ell(\la)+1}\frac{\theta_p(tu_i/qz)}{\theta_p(u_i/z)},\\
B^-_\la(z/u;p)&=&
\prod_{i=1}^{\ell(\la)}
\frac{\theta_p(tz/u_i)}{\theta_p(qz/u_i)}
\prod_{i=1}^{\ell(\la)+1}\frac{\theta_p(qz/tu_i)}{\theta_p(z/u_i)}.
 \ena
\end{thm}
A direct proof showing these actions satisfy the defining relation of $\cU_{q,t,p}$ is given in Appendix \ref{proofVerticalRep}. 
In Appendix \ref{IndDer} we also give an inductive derivation of Theorem \ref{actionUqtpgl1}.

\bigskip
 \noindent
 {\it Remark 1.}\ 
In the trigonometric case obtained by taking $p\to0$,  the level (0,1) representation in Theorem \ref{actionUqtpgl1} 
is identified with the geometric representation of $\cU_{q,t}$ on $\bigoplus_N\mathrm{K}_T(\mathrm{Hilb}_{N}(\C^2))$, the direct sum of the  $T=\C^\times\times \C^\times$ equivariant $K$-theory of the Hilbert scheme of $N$ points on $\C^2$\cite{FT, FFJMM}. There the basis $\{\ket{\la}_u \}$ in $\F_u$ is identified with the fixed point basis $\{[\la]\}$  in $\mathrm{K}_T(\mathrm{Hilb}_{N}(\C^2))$. 
We conjecture \cite{K19} that the same is true in the elliptic case. Namely, if one could properly formulate  a geometric action of $\cU_{q,t,p}$ on the direct sum of the equivariant elliptic cohomology of the Hilbert scheme $\bigoplus_N\mathrm{E}_T(\mathrm{Hilb}_{N}(\C^2))$, it should be identified with the level (0,1) representation in Theorem \ref{actionUqtpgl1} by identifying $\ket{\la}_u$ with  the fixed point class $[\la]$ in $\bigoplus_N\mathrm{E}_T(\mathrm{Hilb}_{N}(\C^2))$, where the latter bases can be realized in terms of the elliptic stable envelopes\cite{AO} on  $\mathrm{E}_T(\mathrm{Hilb}_{N}(\C^2))$ constructed in \cite{Smirnov18}  by 
\be
&&[\la]=\sum_{\mu\in\cP^+}\Stab^{-1}_{\gC}(\mu)\vert_{\la} \Stab_{\gC}(\mu)
\en
in the similar way as the case of the equivariant elliptic cohomology of the partial flag variety \cite{Konno18}.
 
 \bigskip
 
  \noindent
 {\it Remark 2.}\ 
Again from the trigonometric limit of \eqref{actpsipn}, one finds the eigenvalue of $\al_1(=a_1)$ on $\ket{\la}_u^{(N)}$ is given by
\be
&&-(1-t/q)(1-t)t^{-N}\sum_{j=1}^Nq^{\la_j}t^{N-j}.
\en
Under the automorphism in Theorem \ref{Thm:Miki}, $a_1$ is identified with $(1-t/q)X^+_0$.
Hence we obtain 
the eigenvalue of $X^+_0$ as $-(1-t)t^{-N}\sum_{j=1}^Nq^{\la_j}t^{N-j}$. Then under the identification of $q^{N|\la|}\ket{\la}_u^{(N)}$ with the specialization of the Macdonald symmetric polynomial $P_\la(t^\rho q^\la)$, where the variables $x=(x_1,\cdots, x_N)$ are specialized by $t^\rho q^{\la}=(t^{N-1}q^{\la_1},t^{N-2}q^{\la_2},\cdots,q^{\la_N})$. This is consistent to the identification that the trigonometric limit of \eqref{actxpn} is the Pieri formula \cite{MacBook}:
\be
&&e_1(x)P_\la(x)=\sum_{i=1}^N\psi_{\la+\bo_i/\la}P_{\la+\bo_i}(x)
\en
where $e_1(x)=\sum_{j=1}^Nx_j$ and 
\be
&&\psi_{\la+\bo_i/\la}=\prod_{j=1}^{i-1}\frac{(1-tu_i/u_j)(1-qu_i/tu_j)}{(1-u_i/u_j)(1-qu_i/u_j)}.
\en
See also \cite{FFJMM} for such identification and its extension to the case $N\to \infty$ in \cite{FT}.  Hence it is natural to expect that the action of the elliptic generator $x^+_0=\oint_{0}\frac{dz}{2\pi i z} x^+(z)$ on $\ket{\la}_u^{(N)}$ in \eqref{actxpn} gives an elliptic analogue of the Pieri formula \cite{K19} and the eigenvalue of $x^+_0$ should be given in terms of an elliptic analogue of $e_1(t^\rho q^\la)$, which is unknown yet, times $t^{-N}a^+(p)$. 

These two remarks suggest that  $x^+_0$ behaves as the elliptic Ruijsenaars difference operator 
on $W^{(N)}(u)$ and the bases vector $\ket{\la}_u^{(N)}$ as an elliptic analogue of the Macdonald symmetric polynomials with the eigenvalue given by an elliptic analogue of $e_1(t^\rho q^\la)$, which should satisfy the elliptic analogue of the Pieri formula \eqref{actxpn}. 
See also \cite{LQW,Na14}, where the fixed point classes in the equivariant homology groups of Hilbert schemes of points on $\C^2$ are identified with the Jack symmetric functions.  
 In 
 $\S$\ref{ellRuijsenaars}, 
 we give more explicit relation of $x^+_0$ to the elliptic Ruijsenaars difference operator.

 \subsubsection{Some elliptic formulas}\lb{RelqFock}
 Let us 
consider the combinatorial factors $c_\la, c'_\la$ appearing in the 
inner product of the Macdonald symmetric functions as
\bea
&&\langle P_\la,P_\la \rangle_{q,t}=\frac{c'_\la}{c_\la},\\
&&c_\la=\prod_{\Box\in\la}(1-q^{a(\Box)}t^{\ell(\Box)+1})
=\prod_{1\leq i\leq j\leq \ell(\la)}\frac{(q^{\la_i-\la_{j}}t^{j-i+1};q)_\infty}{(q^{\la_i-\la_{j+1}}t^{j-i+1};q)_\infty},\nn\\
&&c'_\la=\prod_{\Box\in\la}(1-q^{a(\Box)+1}t^{\ell(\Box)})
=\prod_{1\leq i\leq j\leq \ell(\la)}\frac{(q^{\la_i-\la_{j}+1}t^{j-i};q)_\infty}{(q^{\la_i-\la_{j+1}+1}t^{j-i};q)_\infty}. \lb{clacpla}
\ena 
Here $a(\Box)\equiv a_\la(\Box)=\la_i-j,\ \ell(\Box)\equiv \ell_{\la}(\Box)=\la'_j-i$ for $\Box=(i,j)\in \la$. 
The second expressions for $c_\la, c'_{\la}$ follow from the formula 
due to Macdonald\cite{MacBook}
\bea
&&(1-q)\sum_{(i,j)\in \la}q^{\la_i-j}t^{\la'_j-i+1}=t\sum_{1\leq i\leq \ell(\la)}q^{\la_i-\la_j}t^{j-i}
-\sum_{1\leq i< \ell(\la)+1}q^{\la_i-\la_j}t^{j-i}.\lb{Macdoformula}
\ena

We introduce elliptic analogues of $c_\la, c'_\la$ as follows. 
\bea
&&c_\la(p)=\prod_{\Box\in\la}{\theta_p(q^{a(\Box)}t^{\ell(\Box)+1})}
=\prod_{1\leq i\leq j\leq \ell(\la)}\frac{\Gamma(q^{\la_i-\la_{j+1}}t^{j-i+1};q,p)}{\Gamma(q^{\la_i-\la_{j}}t^{j-i+1};q,p)},\lb{clap}\\
&&c'_\la(p)=\prod_{\Box\in\la}{\theta_p(q^{a(\Box)+1}t^{\ell(\Box)})}
=\prod_{1\leq i\leq j\leq \ell(\la)}\frac{\Gamma(q^{\la_i-\la_{j+1}+1}t^{j-i};q,p)}{\Gamma(q^{\la_i-\la_{j}+1}t^{j-i};q,p)}. \lb{cplap}
\ena
One can verify the following.
\begin{prop}
\bea
&&\frac{c_{\la+\bo_k}(p)}{c_{\la}(p)}=\prod_{i=1}^{k-1}\frac{\theta_p(tq^{-1}u_i/u_k)}{\theta_p(q^{-1}u_i/u_k)}\frac{\prod_{i=k+1}^{\ell(\la)+1}\theta_p(u_k/u_i)}{\prod_{i=k+1}^{\ell(\la)}\theta_p(tu_k/u_i)},\lb{def:cla}\\
&&\frac{c'_{\la+\bo_k}(p)}{c'_{\la}(p)}=\prod_{i=1}^{k-1}\frac{\theta_p(u_i/u_k)}{\theta_p(t^{-1}u_i/u_k)}\frac{\prod_{i=k+1}^{\ell(\la)+1}\theta_p(qt^{-1}u_k/u_i)}{\prod_{i=k+1}^{\ell(\la)}\theta_p(qu_k/u_i)}.
\ena
\end{prop}
The following formulas are useful in $\S$\ref{Sec:VOs}. 
\begin{prop}
Let us consider
\bea
A^{+'}_{\la,i}(p)
&=&  \prod_{j=i+1}^{\ell(\la)}\frac{\theta_p(tu_i/u_j)}{\theta_p(qu_i/u_j)}\prod_{j=i+1}^{\ell(\la)+1}\frac{\theta_p(qu_i/tu_j)}{\theta_p(u_i/u_j)},
\lb{tildepApp}\\
A^{-'}_{\la,i}(p)
&=&\prod_{j=1}^{i-1}\frac{\theta_p(tu_j/u_i)\theta_p(qu_j/tu_i)}{\theta_p(qu_j/u_i)\theta_p(u_j/u_i)}   \lb{tildepAmp}.
 \ena
We have  
\bea
A^{+'}_{\la,i}(p)
&=&\frac{c_\la(p) c'_{\la+\bo_i}(p)}{c'_\la(p) c_{\la+\bo_i}(p)}\tA^{+}_{\la,i}(p)
\lb{tApptAp}\\
A^{-'}_{\la,i}(p)
&=&\frac{c_\la(p) c'_{\la-\bo_i}(p)}{c'_\la(p) c_{\la-\bo_i}(p)}\tA^{-}_{\la,i}(p).
\lb{tAmptAm}
 \ena
\end{prop}

By a direct calculation one can also verify the following formulas. 
\begin{prop}  
\bea
&&\tA^{+'}_{\la,i}(p)=(q/t)\tA^{-}_{\la+\bo_i,i}(p),\lb{AppAm}\\
&&\tA^{-'}_{\la,i}(p)=(t/q)\tA^{+}_{\la-\bo_i,i}(p).\lb{AmpAp}
\ena
\end{prop}
Hence we have
\begin{prop}\lb{Ala2Alap1}
\bea
&&\frac{c_\la(p)}{c_{\la+{\bf1}_i}(p)}\frac{c'_{\la+{\bf1}_i}(p)}{c'_\la(p)}A^+_{\la,i}(p)
=(q/t)A^-_{\la+{\bf1}_i,i}(p),\\
&&\frac{c_\la(p)}{c_{\la-{\bf1}_i}(p)}\frac{c'_{\la-{\bf1}_i}(p)}{c'_\la(p)}A^-_{\la,i}(p)=(t/q)A^+_{\la-{\bf1}_i,i}(p).
\ena
\end{prop}

We also need to introduce $N_\la(p)$ and  $N'_\la(p)$ by 
 \bea
&&\frac{N_\la(p)}{N_{\la+{\bf1}_i}(p)}=\prod_{j=1}^{i-1}
\frac{(pu_j/tu_i;p)_\infty(ptu_j/qu_i;p)_\infty}{(pu_j/qu_i;p)_\infty(pu_j/u_i;p)_\infty}
\prod_{j=i+1}^{\ell(\la)}
\frac{(pqu_i/u_j;p)_\infty}{(ptu_i/u_j;p)_\infty}
\prod_{j=i+1}^{\ell(\la)+1}\frac{(pu_i/u_j;p)_\infty}{(pqu_i/tu_j;p)_\infty},\nn\\
&&\lb{defrel:N}\\
&&\frac{N'_\la(p)}{N'_{\la+{\bf 1}_i}(p)}=\prod_{j=1}^{i-1}
\frac{(pu_i/u_j;p)_\infty(pqu_i/u_j;p)_\infty}{(pqu_i/tu_j;p)_\infty(ptu_i/u_j;p)_\infty} 
\prod_{j=i+1}^{\ell(\la)}
\frac{(pu_j/tu_i;p)_\infty}{(pu_j/qu_i;p)_\infty}
\prod_{j=i+1}^{\ell(\la)+1}
\frac{(ptu_j/qu_i;p)_\infty}{(pu_j/u_i;p)_\infty}.\nn\\
&&\lb{defrel:Np}
\ena
One can derive the following expressions.
 \bea
 N_\la(p)&=&\prod_{1\leq i\leq j\leq \ell(\la)}\frac{(pqu_i/u_j;q,p)_\infty}{(ptu_i/u_j;q,p)_\infty}
 \prod_{1\leq i<j\leq \ell(\la)+1}\frac{(pu_i/u_j;q,p)_\infty}{(pqu_i/tu_j;q,p)_\infty}\nn\\
 &=&\prod_{\Box\in \la}\frac{(pq^{a(\Box)+1}t^{\ell(\Box)};p)_\infty}
 {(pq^{a(\Box)}t^{\ell(\Box)+1};p)_\infty},\\
 N'_\la(p)&=&\prod_{1\leq i\leq j\leq \ell(\la)}\frac{(pu_j/u_i;q,p)_\infty}{(pqu_j/tu_i;q,p)_\infty}
 \prod_{1\leq i<j\leq \ell(\la)+1}\frac{(pqu_j/u_i;q,p)_\infty}{(ptu_j/u_i;q,p)_\infty}\nn\\
 &=&\prod_{\Box\in \la}\frac{(pq^{-a(\Box)}t^{-\ell(\Box)-1};p)_\infty}
 {(pq^{-a(\Box)-1}t^{-\ell(\Box)};p)_\infty}.
\ena
Then one finds the following property.
\begin{prop}\lb{clapcplap}
\bea
&&\frac{c'_\la}{c_\la}\frac{N_\la(p)}{N'_\la(p)}=\frac{c'_\la(p)}{c_\la(p)}.
\ena
\end{prop}

\subsection{Level (0,0) representation and elliptic Ruijsenaars operators}\lb{ellRuijsenaars}
The two remarks after Proposition \ref{clapcplap}
suggest a connection of the level (0,1) or (0,0) representations  of $\cU_{q,t,p}$ to a possible  elliptic analogue of the Macdonald symmetric functions, which are expected to be eigen functions of the elliptic Ruijsenaars difference operator. 
In this subsection, we show a direct relation between the level (0,0) representation and 
the elliptic Ruijsenaars difference operator 
\be
&&D=\sum_{i=1}^N\prod_{j\not=i}\frac{\theta_p(tx_i/x_j)}{\theta_p(x_i/x_j)}T_{q,x_i}, 
\en
which acts on  
$\C[[x_1^{\pm1},\cdots,x_N^{\pm1}]]$. Here $T_{q,x_i}$ denotes the shift operator
\be
&&T_{q,x_i}f(x_1,\cdots,x_i,\cdots x_n)=f(x_1,\cdots,q x_i,\cdots x_n).
\en

The following theorem is an elliptic analogue of Proposition 3.3 in \cite{Miki}.
\begin{thm}
The following assignment gives a level $(0,0)$ representation of $\cU_{q,t,p}$ on 
$\C[[x_1^{\pm1},\cdots,x_N^{\pm1}]]$.
\bea
&&x^+(z)=a^+(p)\sum_{i=1}^N\prod_{j\not=i}\frac{\theta_p(tx_i/x_j)}{\theta_p(x_i/x_j)}\delta(x_i/z)T_{q,x_i},\\
&&x^-(z)=-a^-(p)\sum_{i=1}^N\prod_{j\not=i}\frac{\theta_p(t^{-1}x_i/x_j)}{\theta_p(x_i/x_j)}\delta(q^{-1}x_i/z)T^{-1}_{q,x_i},\\
&&\alpha_m=\frac{(1-t^{-m})(1-(q/t)^{-m})}{m}\sum_{j=1}^N x_j^m\qquad (m\in \Z\backslash \{0\}),
\ena
or 
\bea
&&\psi^+(z)=\prod_{j=1}^{N}\frac{\theta_p(t^{-1}x_j/z)\theta_p(q^{-1}tx_j/z)}{\theta_p(x_j/z)\theta_p(q^{-1}x_j/z)},\\
&&\psi^-(z)=\prod_{j=1}^{N}\frac{\theta_p(tz/x_j)\theta_p(qt^{-1}z/x_j)}{\theta_p(z/x_j)\theta_p(qz/x_j)}.
\ena
In particular,  the zero-mode $x^+_0=\oint_{|z|=0}\frac{dz}{2\pi i z} x^+(z)$ acts as the elliptic Ruijsenaars difference operator
\be
&&x^+_0=a^+(p)\sum_{i=1}^N\prod_{j\not=i}\frac{\theta_p(tx_i/x_j)}{\theta_p(x_i/x_j)}
T_{q,x_i}
\en

\end{thm}
\noindent
{\it Proof.}\ Let us check the relation \eqref{ellrelxpxm}. 
\be
[x^+(z),x^-(w)]
&=&-a^+(p)a^-(p)\delta(z/w)\sum_{i=1}^N\left(
\prod_{j\not=i}
\frac{\theta_p(tx_i/x_j)}{\theta_p(x_i/x_j)}
\prod_{k\not=i}\frac{\theta_p(qx_i/tx_k)}{\theta_p(qx_i/x_k)}\delta(x_i/z)\right.\\
&&\left.\qquad\qquad\qquad\qquad\qquad-\prod_{j\not=i}\frac{\theta_p(tx_i/qx_j)}{\theta_p(x_i/qx_j)}
\prod_{k\not=i}\frac{\theta_p(x_i/tx_k)}{\theta_p(x_i/x_k)}\delta(q^{-1}x_i/z)
\right).
\en
Then the relation \eqref{ellrelxpxm} holds by the following Lemma.
\begin{lem}\lb{limDiffTheta}
\be
&&
\prod_{j=1}^N\frac{\theta_p(qz/tx_j)\theta_p(tz/x_j)}{\theta_p(z/x_j)\theta_p(qz/x_j)}\Biggl|_+
-\prod_{j=1}^N\frac{\theta_p(qz/tx_j)\theta_p(tz/x_j)}{\theta_p(z/x_j)\theta_p(qz/x_j)}\Biggl|_-
\\
&&=-\frac{\theta_p(q/t)\theta_p(t)}{(p;p)_\infty^2\theta_p(q)}
\sum_{i=1}^N\left(
\prod_{j\not=i}
\frac{\theta_p(tx_i/x_j)}{\theta_p(qx_i/x_j)}
\prod_{k\not=i}\frac{\theta_p(qx_i/tx_k)}{\theta_p(x_i/x_k)}\delta(x_i/z)\right.\\
&&\left.\qquad\qquad\qquad\qquad\qquad-\prod_{j\not=i}\frac{\theta_p(tx_i/qx_j)}{\theta_p(x_i/qx_j)}
\prod_{k\not=i}\frac{\theta_p(x_i/tx_k)}{\theta_p(x_i/qx_k)}\delta(q^{-1}x_i/z)
\right).
\en
\end{lem}
\noindent
{\it Proof.}\ Note the partial fraction expansion formula, see for example \cite{Rosengren},  
\bea
&&\theta_p(s/b_{m+1})\prod_{j=1}^m\frac{\theta_p(s/b_j)}{\theta_p(s/a_j)}=-
\sum_{k=1}^m\frac{1}{\theta_p(a_k/s)}\frac{\prod_{j=1}^{m+1}\theta_p(a_k/b_j)}{\prod_{j=1\atop \not=k}^m\theta_p(a_k/a_j)}\lb{PFE}
\ena
with the balancing condition $b_1\cdots b_{m+1}=a_1\cdots a_m s$.
Let us consider the case $m=2N$ and  take 
\be
&&\hspace{-1cm}s=z,\quad a_j=x_j,\quad  a_{N+j}=q^{-1}x_j,\quad b_j=(t/q)x_j,\quad b_{N+j}=\beta^{-1} x_j\quad (j=1,\cdots,N), \quad b_{2N+1}=(\beta/t)^N z, 
\en
for some constant $\beta$. 
Then one has
\be
&&\prod_{j=1}^N\frac{\theta_p(qz/tx_j)}{\theta_p(z/x_j)} 
\prod_{j=1}^N\frac{\theta_p(\beta z/x_j)}{\theta_p(qz/x_j)} \\
&&=-\sum_{i=1}^N\frac{\theta_p((t/\beta)^{N} x_i/z)\theta_p(q/t)}{\theta_p((t/\beta)^{N})\theta_p(x_i/z)}
\prod_{j=1\atop \not=i}^N\frac{\theta_p(qx_i/tx_j)}{\theta_p(x_i/x_j)} 
\prod_{k=1\atop \not=i}^N\frac{\theta_p(\beta x_i/x_k)}{\theta_p(qx_i/x_k)}\times\frac{\theta_p(\beta)}{\theta_p(q)} \\
&&-\sum_{i=1}^N\frac{\theta_p((t/\beta)^{N} x_i/qz)\theta_p(\beta/q)}{\theta_p((t/\beta)^{N})\theta_p(x_i/qz)}
\frac{\theta_p(t^{-1})}{\theta_p(q^{-1})}
\prod_{j=1\atop \not=i}^N\frac{\theta_p(x_i/tx_j)}{\theta_p(x_i/qx_j)} 
\prod_{k=1\atop \not=i}^N\frac{\theta_p(\beta x_i/qx_k)}{\theta_p(x_i/x_k)}.
\en
Hence we obtain
\be
&&\prod_{j=1}^N\frac{\theta_p(qz/tx_j)}{\theta_p(z/x_j)} 
\prod_{j=1}^N\frac{\theta_p(\beta z/x_j)}{\theta_p(qz/x_j)} \Biggl|_+- 
\prod_{j=1}^N\frac{\theta_p(qz/tx_j)}{\theta_p(z/x_j)} 
\prod_{j=1}^N\frac{\theta_p(\beta z/x_j)}{\theta_p(qz/x_j)} \Biggl|_- 
\\
&&=-\frac{\theta_p(\beta)\theta_p(q/t)}{(p;p)_\infty^2\theta_p(q)}\sum_{i=1}^N
\delta(x_i/z)\prod_{j=1\atop \not=i}^N\frac{\theta_p(qx_i/tx_j)}{\theta_p(x_i/x_j)} 
\prod_{k=1\atop \not=i}^N\frac{\theta_p(\beta x_i/x_k)}{\theta_p(qx_i/x_k)}\\
&&\quad-\frac{\theta_p(t^{-1})\theta_p(\beta/q)}{(p;p)_\infty^2\theta_p(q^{-1})}
\sum_{i=1}^N\delta(q^{-1}x_i/z)\prod_{j=1\atop \not=i}^N\frac{\theta_p(x_i/tx_j)}{\theta_p(x_i/qx_j)} 
\prod_{k=1\atop \not=i}^N\frac{\theta_p(\beta x_i/qx_k)}{\theta_p(x_i/x_k)}
\en
Then taking the limit $\beta \to t$, one obtains the desired formula. \qed

\section{The Vertex Operators}\lb{Sec:VOs}
We construct the two vertex operators $\Phi(u)$ and $\Psi^*(u)$  of $\cU_{q,t,p}$ called the type I  and the type II dual vertex operators\cite{JM} as intertwining operators of $\cU_{q,t,p}$-modules. 
These two vertex operators are the elliptic analogues of those constructed in \cite{AFS}, 
whose matrix elements reproduce the refined  topological vertex in \cite{Taki, IKV, AK}.     
We also construct a shifted inverse 
of them denoted by $\Phi^*(u)$ and $\Psi(u)$, respectively. 
These vertex operators turn out to be useful to realize the affine quiver $W$ algebra and 
instanton calculus in the affine quiver gauge theories. See $\S$\ref{WA0} and $\S$\ref{InstPF}.

\subsection{The type I vertex operator}
The type I vertex operator is the intertwining operator 
\be
\Phi(u) : \F^{(1,N+1)}_{-uv}\to \F^{(0,1)}_u\tot \F^{(1,N)}_v
\en
 w.r.t. the comultiplication $\Delta$ satisfying 
\bea
&&\Delta(x)\Phi (u)=\Phi (u)x\qquad \qquad (\forall x\in \cU_{q,t,p}).\lb{intrelPhis}
\ena
We define the components of $\Phi(u)$ by
\bea
&&\Phi (u)\ket{\xi}=\sum_{\la\in \cP^+}\ket{\la}'_u\tot \Phi _\la(u)\ket{\xi} \qquad 
\forall \ket{\xi}\in \F^{(1,N+1)}_{-uv},\lb{compPhis}
\ena
where we set
\bea
&&\ket{\la}'_u=\frac{c_\la(p)}{c'_\la(p)}\ket{\la}_u. \lb{ketlapketla}
\ena

\begin{lem}\lb{lemtypeI}
The intertwining relation \eqref{intrelPhis} reads 
\bea
&&\Phi _\la(u)\psi^+((t/q)^{1/4}z)=(q/t)^{1/2}B^+_\la(u/z;p)\psi^+((t/q)^{1/4}z)\Phi _\la(u),\lb{Phipsip}\\
&&\Phi _\la(u)\psi^-((t/q)^{-1/4}z)=(q/t)^{-1/2}B^-_\la(z/u;p)\psi^-((t/q)^{-1/4}z)\Phi _\la(u),\lb{Phipsim}\\
&&\Phi _\la(u)x^+((t/q)^{-1/4}z)=x^+((t/q)^{-1/4}z)\Phi _\la(u)\nn\\
&&\qquad\qquad+qf(1;p)^{-1}\psi^-((t/q)^{-1/4}z)\sum_{i=1}^{\ell(\la)+1}
a^-(p)\tA^-_{\la,i}(p)\delta(q^{-1}u_i/z)\Phi _{\la-{\bf1}_i}(u),\lb{Phixp}\\
&&\Phi _\la(u)x^-((t/q)^{1/4}z)=(q/t)^{1/2}B^+_{\la}(u/z;p)x^-((t/q)^{1/4}z)\Phi _\la(u)\nn\\
&&\qquad\qquad\qquad+q^{-1}f(1;p)(q/t)^{1/2}\sum_{i=1}^{\ell(\la)}
a^+(p) \tA^+_{\la,i}(p)\delta(u_i/z)\Phi _{\la+{\bf1}_i}(u).\lb{Phixm}
\ena
\end{lem}
\noindent
{\it Proof.}
For example, let us consider 
\be
&&\Delta(x^+(z))\Phi(u)\ket{\xi}=\Phi(u)x^+(z)\ket{\xi}\qquad \forall\ket{\xi}\in \F^{(1,N+1)}_{-uv}.
\en
\be
RHS&=&\sum_{\la}\ket{\la}'_u\tot \Phi_\la(u)x^+(z)\ket{\xi},\\
LHS&=&\sum_{\la}\ket{\la}'_u\tot x^+(z)\Phi_\la(u)\ket{\xi}+
\sum_{\la}x^+((t/q)^{1/4}z)\ket{\la}'_u\tot \psi^-(z)\tPhi_\la(u)\ket{\xi}\nn\\
&=&\sum_{\la} \ket{\la}'_u\tot x^+(z)\Phi_\la(u)\ket{\xi}\nn\\
&&\qquad+
\sum_{\la}\sum_{i=1}^{\ell(\la)+1}a^+(p)A^{+'}_{\la,i}(p)
\delta((t/q)^{-1/4}u_i/z)
\ket{{\la+{\bf1}_i}}'_u\tot \psi^-(z)\Phi_\la(u)\ket{\xi}\nn\\
&=&\sum_{\la} \ket{\la}'_u\tot x^+(z)\Phi_\la(u)\ket{\xi}\nn\\
&&\qquad+
qf(1,p)^{-1}\sum_{\la}\sum_{i=1}^{\ell(\la)+1}a^-(p)A^-_{\la+{\bf 1}_i,i}(p)\delta((t/q)^{-1/4}u_i/z)
\ket{\la+{\bf1}_i}'_u\tot \psi^-(z)\Phi_\la(u)\ket{\xi}\nn\\
&=&\sum_{\la} \ket{\la}'_u\tot \left\{x^+(z)\Phi_\la(u)
+qf(1,p)^{-1}\sum_{i=1}^{\ell(\la)}a^-(p)A^-_{\la,i}(p)\delta((t/q)^{-1/4}q^{-1}u_i/z)
\psi^-(z)\Phi_{\la-{\bf 1}_i}(u)\right\}\ket{\xi}
\en
The third equality follows from \eqref{AppAm}.
\qed

By using the representations in Theorem \ref{level1N} and Theorem \ref{actionUqtpgl1}, one can solve these intertwining relations, and obtain the following result. 
\begin{thm}\lb{typeIVO}
\be
&&{\Phi}_\la(u)=\frac{q^{n(\la')}
N_\la(p)
t^*(\la,u,v,N)}{c_\la }{\tPhi}_\la(u),\nn\\
&&{\tPhi}_\la(u)=:{\Phi}_\emptyset(u)\prod_{i=1}^{\ell(\la)}\prod_{j=1}^{\la_i} \widetilde{x}^-((t/q)^{1/4}q^{j-1}t^{-i+1}u):,\\
&&\Phi _\emptyset(u)
=\exp\left\{-\sum_{m>0}\frac{1}{\kappa_m}\al'_{-m} ((t/q)^{1/2}u)^m
\right\}\exp\left\{\sum_{m>0}\frac{1}{\kappa_m}\al'_{m}((t/q)^{1/2} u)^{-m}
\right\}
\en
where 
$x^-(z)=u^{-1}z^{N}(t/q)^{-3N/4}\widetilde{x}^-(z)$ on $\F^{(1,N)}_u$ and
\bea
&&n(\la)=\sum_{i\geq 1}(i-1)\la_i,\qquad n(\la')=\sum_{i\geq 1}(i-1)\la'_i=\sum_{i\geq 1}\frac{\la_i(\la_i-1)}{2},\\
&&t^*(\la,u,v,N)=(q^{-1}v)^{-|\la|}(-u)^{N|\la|}f_\la(q,t)^N,\lb{tsfac}\\
&&f_\la(q,t)=(-1)^{|\la|}q^{n(\la')+|\la|/2}t^{-n(\la)-|\la|/2}. \lb{framingfac}
\ena
\end{thm}
The factor $t^*(\la,u,v,N)$ was introduced in \cite{AFS} and in particular $f_\la(q,t)$ is called 
the framing factor\cite{Taki,AK}. Our vertex operator $\Phi(u)$ is the elliptic analogue of $\Phi^*(u)$ in \cite{AFS}. One should note that our comultiplication is opposite from the one in \cite{AFS}. A proof of the statement is given in Appendix \ref{prooftypeIVO}.

In later sections the following formula is useful. 
\begin{prop}\lb{tildePhi}
\bea
&&\widetilde{\Phi}_\la(u)=:\exp\left( \sum_{m\not=0} \frac{1-t^{m}}{\kappa_m}\cE_{\la,m}\alpha'_m((t/q)^{1/2}u)^{-m}\right) :,
\ena
where we set
\be
&&\cE_{\la,m}=\frac{1}{1-t^m}+\sum_{j=1}^{\ell(\la)}(q^{-m\la_j}-1)t^{m(j-1)}\qquad (m\in \Z_{\not=0}).
\en
\end{prop}
\subsection{The type II dual vertex operator}
The type II dual vertex operator is  the intertwining operator 
\bea
&&\Psi^*(v)\ :\ \F^{(1,N)}_u\tot\F^{(0,1)}_v \to \F^{(1,N+1)}_{-vu}
\ena
satisfying 
\bea
&&x\Psi^*(v)=\Psi^*(v)\Delta(x)\qquad\qquad \forall x \in \cU_{q,t,p}.\lb{IntrelII}
\ena
We call $\Psi^*(v)$ the type II dual vertex operator. 
We define its components by
\bea
&&\Psi^*_\la(v)\ket{\xi}=\Psi^*(v)\left(\ket{\xi}\tot \ket{\la}_v\right)\qquad \forall \ket{\xi}\in \F^{(1,N)}_u.
\lb{compPsis}
\ena 
\begin{lem}\lb{lemtypeIIdual}
The intertwining relation \eqref{IntrelII} is equivalent to
\bea
&&\psi^+((q/t)^{1/4}z)\Psi^*_\la(v)=(q/t)^{1/2}B^+_\la(v/z;p^*)\Psi^*_\la(v)\psi^+((q/t)^{1/4})),\lb{Psipsip}\\
&&\psi^-((q/t)^{-1/4}z)\Psi^*_\la(v)=(q/t)^{-1/2}B^-_\la(z/v;p^*)\Psi^*_\la(v)\psi^-((q/t)^{-1/4})),\\
&&x^-((q/t)^{1/4}z)\Psi^*_\la(v)=\Psi^*_\la(v)x^-((q/t)^{1/4}z)\nn\\
&&\hspace{3cm} +(q/t)^{1/2}a^-(p^*)\sum_{i=1}^{\ell(\la)}A^-_{\la,i}(p^*)\delta(q^{-1}v_i/z)\tPsi^*_{\la-{\bf 1}_i}(v)
\psi^+((q/t)^{1/4}z),\lb{Psixm}
\\
&&x^+((q/t)^{-1/4}z)\Psi^*_\la(v)=(q/t)^{-1/2}B^-_\la(z/v;p^*)\Psi^*_\la(v)x^+((q/t)^{-1/4}z)\nn\\
&&\hspace{3cm} +a^+(p^*)\sum_{i=1}^{\ell(\la)+1}A^+_{\la,i}(p^*)\delta(v_i/z)\tPsi^*_{\la+{\bf1}_i}(v).
\lb{Psixp}
\ena
Here $p^*=pq/t$ associated with $\F^{(1,N)}_u$. 
\end{lem}
One can prove this in the similar way to Lemma \ref{lemtypeI}.

By using the representations in Theorem \ref{level1N}, \ref{actionUqtpgl1} and \eqref{ketlapketla}, 
one can solve these intertwining relations, and obtain the following result. 
\begin{thm}\lb{typeIIdual}
\bea
&&{\Psi}^*_\la(v)=\frac{q^{n(\la')}t(\la,u,v,N)}{c_\la N_\la'(p^*)}{\tPsi}^*_\la(v),\nn\\
&&{\tPsi}^*_\la(v)=
\Psi^*_{\emptyset}(v)\prod_{i=1}^{\ell(\la)}\prod_{j=1}^{\la_i}\widetilde{x}^+((t/q)^{1/4}q^{j-1}t^{-i+1}v),\\
&&\Psi^*_{\emptyset}(v)=\exp\left\{\sum_{m>0}\frac{1}{\kappa_m}\al_{-m} ((t/q)^{1/2}v)^m
\right\}\exp\left\{-\sum_{m>0}\frac{1}{\kappa_m}\al_{m} ((t/q)^{1/2}v)^{-m}
\right\}.
\ena
where $x^+(z)=uz^{-N}(t/q)^{3N/4}\widetilde{x}^+(z)$  on $\F^{(1,N)}_u$, and
\be
&&t(\la,u,v,N)=(-uv)^{|\la|}(-v)^{-(N+1)|\la|}f_\la(q,t)^{-N-1}\lb{tfac}
\en
with $f_\la(q,t)$ given in \eqref{framingfac}.
\end{thm}
Our vertex operator $\Psi^*(u)$ is the elliptic analogue of $\Phi(u)$ in \cite{AFS}. 
A proof of the statement is similar to the one 
in Appendix \ref{prooftypeIVO}.

We have a similar formula to Proposition \ref{tildePhi}. 
\begin{prop}\lb{tildePsi}
\bea
&&\widetilde{\Psi}^*_\la(v)=:\exp\left( -\sum_{m\not=0} \frac{1-t^{m}}{\kappa_m}\cE_{\la,m}\alpha_m((t/q)^{1/2}v)^{-m}\right) :.
\ena
\end{prop}

\subsection{The shifted inverse of  $\Phi(u)$ and $\Psi^*(v)$}
We next introduce the shifted inverse of $\Phi(u)$ and $\Psi^*(v)$ denoted by $\Phi^*(u)$ and $\Psi(v)$, respectively.  
In $\S$\ref{WA0} and $\S$\ref{InstPF}, we show that the vertex operators here and in the previous subsections 
are useful to derive physical and mathematical quantities associated with the Jordan quivers.

\subsubsection{The  vertex operator $\Phi^*(u)$}\lb{Sec:typeIdual} 
Let us consider the linear map
\bea
&&\Phi^*(u) :  \F^{(0,1)}_u\tot \F^{(1,N)}_v \to \F^{(1,N+1)}_{-uv}, \lb{conjtypeIdual}
\ena
whose components are defined by
\bea
&&\Phi^*(u)\left(\ket{\la}'_u\tot \ket{\xi}\right)=\Phi^*_\la(u)\ket{\xi},\qquad \forall \ket{\xi}\in \F^{(1,N)}_v ,\lb{compPhidual}\\
&&{\Phi}^*_\la(u)=\frac{q^{n(\la')}
N'_\la(p)
t(\la,v,up^{-1},N)}{c'_\la }:{\tPhi}_\la(p^{-1}u)^{-1}:.
\ena
Note that  from Proposition \ref{tildePhi} we have
\be
&&:{\tPhi}_\la(p^{-1}u)^{-1}:=: 
\exp\left( -\sum_{m\not=0} \frac{1-t^{m}}{\kappa_m}p^m\cE_{\la,m}\alpha'_m((t/q)^{1/2}u)^{-m}\right)
:
\en

\begin{prop}\lb{typeIdual}
The  vertex operator ${\Phi}^*_\la(u)$ satisfies the following relations. 
\bea
&&\psi^+((t/q)^{1/4}z)\Phi^* _\la(u)=(t/q)^{-1/2}B^+_\la(p^{-1}u/z;p)\Phi^*_\la(u)\psi^+((t/q)^{1/4}z),\lb{Phispsip}\\
&&\psi^-((t/q)^{-1/4}z)\Phi^*_\la(u)=(t/q)^{1/2}B^-_\la(pz/u;p)\Phi^*_\la(u)\psi^-((t/q)^{-1/4}z),\lb{Phispsim}\\
&&x^+((t/q)^{-1/4}z)\Phi^*_\la(u)=\Phi^*_\la(u)x^+((t/q)^{-1/4}z)\nn\\
&&\qquad\qquad+
(t/q)^{-1/2}a^+(p)\sum_{i=1}^{\ell(\la)+1}
\tA^{+'}_{\la,i}(p)\delta(p^{-1}tu_i/qz)\Phi^*_{\la+{\bf1}_i}(u)\psi^+((t/q)^{1/4}qz/t),\lb{Phisxp}\\
&&x^-((t/q)^{1/4}z)\Phi^*_\la(u)=(t/q)^{-1/2}B^+_{\la}(p^{-1}u/z;p)\Phi^*_\la(u)x^-((t/q)^{1/4}z)\nn\\
&&\qquad\qquad\qquad+
(t/q)a^-(p)\sum_{i=1}^{\ell(\la)}
 \tA^{-'}_{\la,i}(p)\delta(p^{-1}q^{-1}u_i/z)\Phi^*_{\la-{\bf1}_i}(u).\lb{Phisxm}
\ena

\end{prop}
A proof of this statement is given in Appendix \ref{prooftypeIdualVO}.

These relations allow us to expect that  $\Phi^*(u)$ is the intertwinner 
satisfying 
\bea
&&\Phi^*(u)\Delta(x)=x\Phi^*(u)\qquad \qquad \forall x\in \cU_{q,t,p}.\lb{intrelPhidual}
\ena
However it turns out this is not the case. In fact one can derive the following relations from 
\eqref{intrelPhidual}, which are slightly different from those in Proposition \ref{typeIdual}.  
\begin{prop}
The intertwining relation \eqref{intrelPhidual} reads 
\bea
&&\psi^+((t/q)^{1/4}z)\Phi^* _\la(u)=(t/q)^{-1/2}B^+_\la(u/z;p)\Phi^*_\la(u)\psi^+((t/q)^{1/4}z),\lb{Phidualpsip}\\
&&\psi^-((t/q)^{-1/4}z)\Phi^*_\la(u)=(t/q)^{1/2}B^-_\la(z/u;p)\Phi^*_\la(u)\psi^-((t/q)^{-1/4}z),\lb{Phidualpsim}\\
&&x^+((t/q)^{-1/4}z)\Phi^*_\la(u)=\Phi^*_\la(u)x^+((t/q)^{-1/4}z)\nn\\
&&\qquad\qquad+a^+(p) 
\sum_{i=1}^{\ell(\la)+1}
\tA^{+'}_{\la,i}(p)\delta(u_i/z)\Phi^*_{\la+{\bf1}_i}(u)\psi^-((t/q)^{-1/4}z),\lb{Phidualxp}\\
&&x^-((t/q)^{1/4}z)\Phi^*_\la(u)=(t/q)^{-1/2}B^+_{\la}(u/z;p)\Phi^*_\la(u)x^-((t/q)^{1/4}z)\nn\\
&&\qquad\qquad\qquad+
(t/q)^{-1/2}a^-(p) \sum_{i=1}^{\ell(\la)}
\tA^{-'}_{\la,i}(p)\delta(q^{-1}u_i/z)\Phi^*_{\la-{\bf1}_i}(u).\lb{Phidualxm}
\ena
\end{prop}
This discrepancy is probably due to a lack of understanding the dual representation to 
$\F^{(0,1)}_u$.  It is hence  an open problem to find a representation theoretical meaning of $\Phi^*(u)$. 

\subsubsection{The vertex operator $\Psi(v)$}\lb{Sec:typeII}
Similar to $\Phi^*(u)$, 
we  
 consider the linear map
 \be
 &&\Psi(v) :  \F^{(1,N+1)}_{-uv} \to  \F^{(1,N)}_u\tot \F^{(0,1)}_v 
 \en
 and define its components  by
\bea
&&\Psi(v)\ket{\xi}=\sum_{\la\in \cP^+}\Psi_\la(v)\ket{\xi}\tot \ket{\la}_v \qquad \forall \ket{\xi}\in \F^{(1,N+1)}_{-uv},
\lb{compPsi}\\
&&{\Psi}_\la(v)=\frac{q^{n(\la')}t^*(\la,p^*v,u,N)}{c'_\la N_\la(p^*)}:{\tPsi}^*_\la(p^*v)^{-1}:.
\ena
Note that  from Proposition \ref{tildePsi} we have
\be
&&:{\tPsi}^*_\la(p^{*}v)^{-1}:=: 
\exp\left( -\sum_{m\not=0} \frac{1-t^{m}}{\kappa_m}p^{*-m}\cE_{\la,m}\alpha_m((t/q)^{1/2}v)^{-m}\right)
:.
\en
One finds the following relations satisfied. 
\begin{prop}\lb{conjtypeII}
\bea
&&\Psi_\la(v)\psi^+((t/q)^{-1/4}z)=(t/q)^{-1/2}B^+_\la(p^{*}v/z;p^*)\psi^+((t/q)^{-1/4}z)\Psi_\la(v),\lb{Psi2psip}\\
&&\Psi_\la(v)\psi^-((t/q)^{1/4}z)=(t/q)^{1/2}B^-_\la(p^{*-1}z/v;p^*)\psi^-((t/q)^{1/4}z)\Psi_\la(v),\lb{Psi2psim}\\
&&\Psi_\la(v)x^-((t/q)^{-1/4}z)=x^-((t/q)^{-1/4}z)\Psi_\la(v)\nn\\
&&\qquad
+(t/q)a^-(p^*)\sum_{i=1}^{\ell(\la)+1}
\tA^{+'}_{\la,i}(p^*)\delta(p^{*}tv_i/qz)\psi^-((t/q)^{-1/4}z)\Psi_{\la+{\bf1}_i}(v),\lb{Psi2xm}\\
&&\Psi_\la(v)x^+((t/q)^{1/4}z)=(t/q)^{1/2}B^-_{\la}(p^{*-1}z/v;p^*)x^+((t/q)^{1/4}z)\Psi_\la(u)\nn\\
&&\qquad\qquad\qquad+
(t/q)^{-1/2}a^+(p^*)\sum_{i=1}^{\ell(\la)}
 \tA^{-'}_{\la,i}(p^*)\delta(p^*v_i/qz)\Psi_{\la-{\bf1}_i}(v).\lb{Psi2xp}
\ena
\end{prop}
A proof of this is similar to the one 
 in Appendix \ref{prooftypeIdualVO}. 

Again these relations are not coincides with  the intertwining relation 
\bea
&&\Psi (v)x=\Delta(x)\Psi (v)\qquad \qquad \forall x\in \cU_{q,t,p},\lb{intrelPsi}
\ena
which is equivalent to the following relations. 
\begin{prop}
\bea
&&\Psi _\la(v)\psi^+((t/q)^{-1/4}z)=(t/q)^{-1/2}B^+_\la(v/z;p^*)\psi^+((t/q)^{-1/4}z)\Psi _\la(v),\lb{Psipsip}\\
&&\Psi _\la(v)\psi^-((t/q)^{1/4}z)=(t/q)^{1/2}B^-_\la(z/v;p^*)\psi^-((t/q)^{1/4}z)\Psi _\la(v),\lb{Psipsim}\\
&&\Psi _\la(v)x^-((t/q)^{1/4}z)=x^-((t/q)^{1/4}z)\nn\\
&&\qquad\qquad+
(t/q)^{1/2}a^-(p^*)\sum_{i=1}^{\ell(\la+\bo_i)}
\tA^{+'}_{\la,i}(p^*)\delta(q^{\la_i}t^{-i+1}v/z)\psi^+((t/q)^{-1/4}z)\Psi_{\la+{\bf1}_i}(v),\lb{Psixm}\\
&&\Psi _\la(v)x^+((t/q)^{1/4}z)=(t/q)^{1/2}B^-_{\la}(z/v;p^*)x^+((t/q)^{1/4}z)\Psi _\la(v)\nn\\
&&\qquad\qquad\qquad+
a^+(p^*)\sum_{i=1}^{\ell(\la-\bo_i)+1}
\tA^{-'}_{\la,i}(p^*)\delta(q^{-1}q^{\la_i}t^{-i+1}v/z)\Psi_{\la-{\bf1}_i}(v).\lb{Psixp}
\ena
\end{prop}
Hence again it is an open problem to find a representation theoretical meaning of $\Psi(u)$.

\section{Affine quiver $W$-algebra $W_{p,p^*}(\Gamma({\widehat{A}_0}))$}\lb{WA0}
One of the importance to consider the elliptic quantum group is that it gives a realization of the  deformed $W$ algebras such as $W_{p,p^*}(\g)$\cite{FrRe} and provides an algebraic structure  i.e. a co-algebra structure,  which enables us to define the intertwining operators  ( the vertex operators ) as  deformation of the primary fields in CFT. In this section, we realize  the deformed 
$W$ algebra $W_{p,p^*}(\Gamma({\widehat{A}_0}))$ associated with the Jordan quiver  
$\widehat{A}_0$\cite{KimPes} by using the level $(1,N)$ representation of $\cU_{q,t,p}$ given in $\S$\ref{sec:level1N}, where in particular $\gamma=(t/q)^{1/2}$, in the same scheme as it was done for $W_{p,p^*}(\g)$ 
in terms of  the level 1 representation of the elliptic quantum group $U_{q,p}(\gh)$\cite{K98, JKOS,KK03,FKO,Konno14}.

\subsection{Screening currents} 
Let us set 
\be
&&s^+_m=\frac{(t/q)^{m/2}}{1-(t/q)^{m}}\al_m, \qquad s^-_m=\frac{(t/q)^{m/2}}{1-(t/q)^{m}}\al'_m.
\en
Then from \eqref{alal} and \eqref{alpalp} with $\gamma=(t/q)^{1/2}$, hence $p^*=pq/t$, one can show the following commutation relations
\be
&&[s^+_m,s^+_n]=-\frac{1}{m}\frac{1-p^m}{1-p^{*m}}(1-q^{m})(1-t^{-m})\delta_{m+n,0},\\
&&[s^-_m,s^-_n]=-\frac{1}{m}\frac{1-p^{*-m}}{1-p^{-m}}(1-q^{m})(1-t^{-m})\delta_{m+n,0}.
\en
Moreover one can rewrite the elliptic currents $x^\pm(z)$ in Theorem \ref{level1N} as
\be
&&x^\pm((t/q)^{1/4}z)=
\left(({t}/{q}\right)^{N/2}{u}/{z^N})^{\pm1}:\exp\left\{\pm\sum_{m\not=0}s^\pm_m \ 
z^{-m} \right\}:.
\en
Hence one of $x^\pm((t/q)^{1/4}z)$ coincides with the screening currents of  $W_{p,p^*}(\Gamma({\widehat{A}_0}))$  \cite{KimPes,Kimura} with 
the {$SU(4)$ $\Omega$-deformation parameters $p, p^*, q, t$
\cite{Nekrasov} satisfying 
\be
p/p^{*}=t/q.
\en
In our knowledge the elliptic quantum toroidal algebra $\cU_{q,t,p}$ is the first 
quantum group structure 
which possesses the {$SU(4)$ $\Omega$-deformation parameters. 

One should also note that in \cite{KimPes,Kimura} Kimura and Pestun constructed only one type of screening current. However it is natural for the (deformed) $W$ algebras that there are two types of screening currents\cite{FrRe,AKOS,FeiFr,BS}. In this sense our realization completes their construction. 

\subsection{Generating function}\lb{genfW}
To obtain the generating function of $W_{p,p^*}(\Gamma({\widehat{A}_0}))$, we apply the 
same scheme as used in \cite{Konno14}. Namely consider the composition of $\Phi(u)$ and $\Phi^*(u)$, which are given in Theorem \ref{typeIVO} and $\S$ \ref{Sec:typeIdual}, respectively. 
 \be
T(u)&=&\Phi^*(u)\Phi(u)=
\sum_{\la\in \cP^+}\Phi^*_\la(u)\Phi_\la(u)\quad : \ \F^{(1,N+1)}_{-uv}\to \F^{(1,N+1)}_{-uv},
\en
Note that one can chose $v\in \C^*$, $N\in \Z$ arbitrarily. 
 Taking the normal ordering one obtains 
\be
\Phi^*_\la(u)\Phi_\la(u)&=&{\cal C}_\la(q,t,p):\tPhi_\la(u)\tPhi^*_\la(u):. 
\en
Then one finds that  the operator part is given by
\bea
&&:\tPhi^*_\la(u)\tPhi_\la(u):=:\prod_{\Box\in A(\la)}Y( u/q^\Box)\prod_{\blacksquare\in R(\la)}
Y((q/t) u/q^\blacksquare)^{-1}:\lb{YYinv}
\ena
with 
\bea
&&Y(u)=:\exp\left\{\sum_{m\not=0} y_m
u^{-m}\right\}:. \lb{Yop}
\ena
Here we set $q^\Box\equiv t^{i-1}q^{-j+1}$ for $\Box=(i,j)\in \la$ etc. and  {$\ds{y_m= \frac{1-p^m}{\kappa_m}(t/q)^{-m/2}\al'_m}$}. 
The symbols $R(\la)$ and $A(\la)$ denote the set of removable  and  addable boxes in the Young diagram $\la$, respectively. The main structure of \eqref{YYinv} is due to  Proposition \ref{tildePhi}, which yields
\bea
&&:\tPhi^*_\la(u)\tPhi_\la(u):=:\exp\left( \sum_{m\not=0} \frac{(1-t^{m})(1-p^m)}{\kappa_m}\cE_{\la,m}\alpha'_m((t/q)^{1/2}u)^{-m}\right) :,\lb{PhilasPhila}
\ena
and the following combinatorial formula.
\begin{prop}\lb{cE2AR}
\bea
&&\cE_{\la,m}=\frac{1}{1-t^{m}}\left(
\sum_{\Box\in A(\la)}q^{m\Box}
-(t/q)^{m}\sum_{\blacksquare\in R(\la)}q^{m\blacksquare}
\right)
\ena
\end{prop} 
\noindent
{\it Proof.}\ 
The statement follows from
\be
&&\cE_{\la,m}=\frac{1}{1-t^m}\left(1-(1-q^{-m})(1-t^m)\sum_{\Box\in \la}q^{m\Box}\right)
\en
 and  
\begin{itemize} 
\item  $(i,\la_i)\in R(\la)$ $\Longleftrightarrow$ $\la_i>\la_{i+1}$
\item if $(i,\la_i)\in R(\la)$ $\Longleftrightarrow$ $(i+1,\la_{i+1}+1)\in A(\la)$
\end{itemize}
for $\la\in \cP^+$. 
\qed

Moreover from \eqref{alpalp} with $\gamma=(t/q)^{1/2}$, one finds the following commutation relation.
\be
&&[y_m, y_n]=-\frac{1}{m}\frac{(1-p^{*m})(1-p^{-m})}{(1-q^m)(1-t^{-m})}\delta_{m+n,0}.
\en
This agrees with the one  in \cite{KimPes}.  

The coefficient part in $\Phi^*_\la(u)\Phi_\la(u)$ can be calculated by combining the normalization factors of the vertex operators and the OPE coefficient.  
The calculation of the latter coefficient is essentially due to the following formula\cite{AK}  obtained by considering a $q$-analogue of (3.6) in \cite{NO}. 
\begin{prop}\lb{cEcE}
\bea
&&-\frac{1-t^m}{1-q^m}\cE_{\la,-m}\cE_{\mu,m}\nn\\
&&=\frac{t^m}{(1-q^m)(1-t^m)}+\sum_{\Box\in \mu}q^{ma_\la(\Box)}t^{m(\ell_\mu(\Box)+1)}
+\sum_{\blacksquare\in \la}q^{-m(a_\mu(\blacksquare)+1)}t^{-m\ell_\la(\blacksquare)}.
\ena
\end{prop}
Then one finds  
\be
&&{\cal C}_\la(q,t,p)
={\cal C}{\frak q}^{|\la|}{\cal Z}^{\widehat{A}_0}_\la(t,q^{-1},p), 
\en
where 
\bea
&&{\frak q}=p^{*-1}p^{N-1}(t/q)^{1/2},\\
&&{\cal Z}^{\widehat{A}_0}_\la(t,q^{-1},p)=\prod_{\Box\in \la}\frac{(1-p\, q^{a(\Box)+1}t^{\ell(\Box)})(1-pq^{-a(\Box)}t^{-\ell(\Box)-1})}
{(1-q^{a(\Box)+1}t^{\ell(\Box)})(1-q^{-a(\Box)}t^{-\ell(\Box)-1})},\lb{ZA0}\\
&&{\cal C}=\frac{(p^{-1}t;q,t,p)_\infty}{(q;q,t,p)_\infty}.
\ena
Note that the sum $\sum_{\la, |\la|=n}{\cal Z}^{\widehat{A}_0}_\la(t,q^{-1},p)$ coincides with  {the equivariant 
$\chi_y$-genus}  of the Hilbert scheme of $n$ points on $\C^2$, 
$
{\rm Hilb}_n(\C^2)$,  \cite{LiLiuZhou, IKV} with $y=p$. 
Note also that one can rewrite \eqref{ZA0}  as
\bea
{\cal Z}^{\widehat{A}_0}_\la(t,q^{-1},p)=
\frac{N_{\la\la}(pq/t)}{N_{\la\la}(q/t)}
\ena
in terms of  the 5d analogue of the Nekrasov function  $N_{\la\mu}(x)$  given by
\bea
&&N_{\la\mu}(x)=
\prod_{\Box\in \la}(1-xq^{-a_\mu(\Box)-1}t^{-\ell_\la(\Box)})
\prod_{\blacksquare\in \mu}(1-xq^{a_\la(\blacksquare)}t^{\ell_\mu(\blacksquare)+1}).\lb{5dNekrasov}
\ena

Hence the whole operator 
\bea
T(u)&=&{\cal C}\, \sum_{\la}{\frak q}^{|\la|}{\cal Z}^{\widehat{A}_0}_\la(t,q^{-1},p)
:\prod_{\Box\in A(\la)}Y( u/q^\Box)\prod_{\blacksquare\in R(\la)}
Y((q/t) u/q^\blacksquare)^{-1}:
\lb{affWop}
\ena
 agrees with the generating function of $W_{p,p^*}(\Gamma({\widehat{A}_0}))$ 
  in  \cite{KimPes,Kimura} up to an over all constant factor.  Note that we have the symmetry
  \be
 &&{\cal Z}^{\widehat{A}_0}_\la(t,q^{-1},p)=(pp^*)^{|\la|}{\cal Z}^{\widehat{A}_0}_\la(t,q^{-1},p^{*-1})
 =p^{2|\la|}{\cal Z}^{\widehat{A}_0}_\la(t^{-1},q,p^{-1}) =(t/q)^{|\la|}{\cal Z}^{\widehat{A}_0}_\la(t^{-1},q,p^{*}).  
 \en

From \eqref{affWop} it is immediate to obtain the rank 1 instanton partition function of the 5d lift of the 4d $\cN=2^*$  theory\cite{AK,HIV,HI,HIKLV} 
by taking the  vacuum expectation value :  
\bea
\bra{0} T(u)\ket{0}={\cal C}\,\sum_\la{\frak q}^{|\la|}{\cal Z}^{\widehat{A}_0}_\la(t,q^{-1},p).
\lb{chiygenus}
\ena
It is then important to recognize that this result and our realization  $T(u)=\sum_\la\Phi^*_\la(u)\Phi_\la(u)$  lead to the identification of  
$T(u)$ with the basic refined topological vertex  depicted in Fig.\ref{3-3vertex}, which 
was introduced in \cite{IKV,HIV,HIKLV}.   
Once obtaining such basic operator one can apply it to various calculations 
 presented in the subsequent sections. 

\begin{figure}[htbp]
\begin{center}
\includegraphics[
scale=0.3, 
height=60mm, width=75mm]{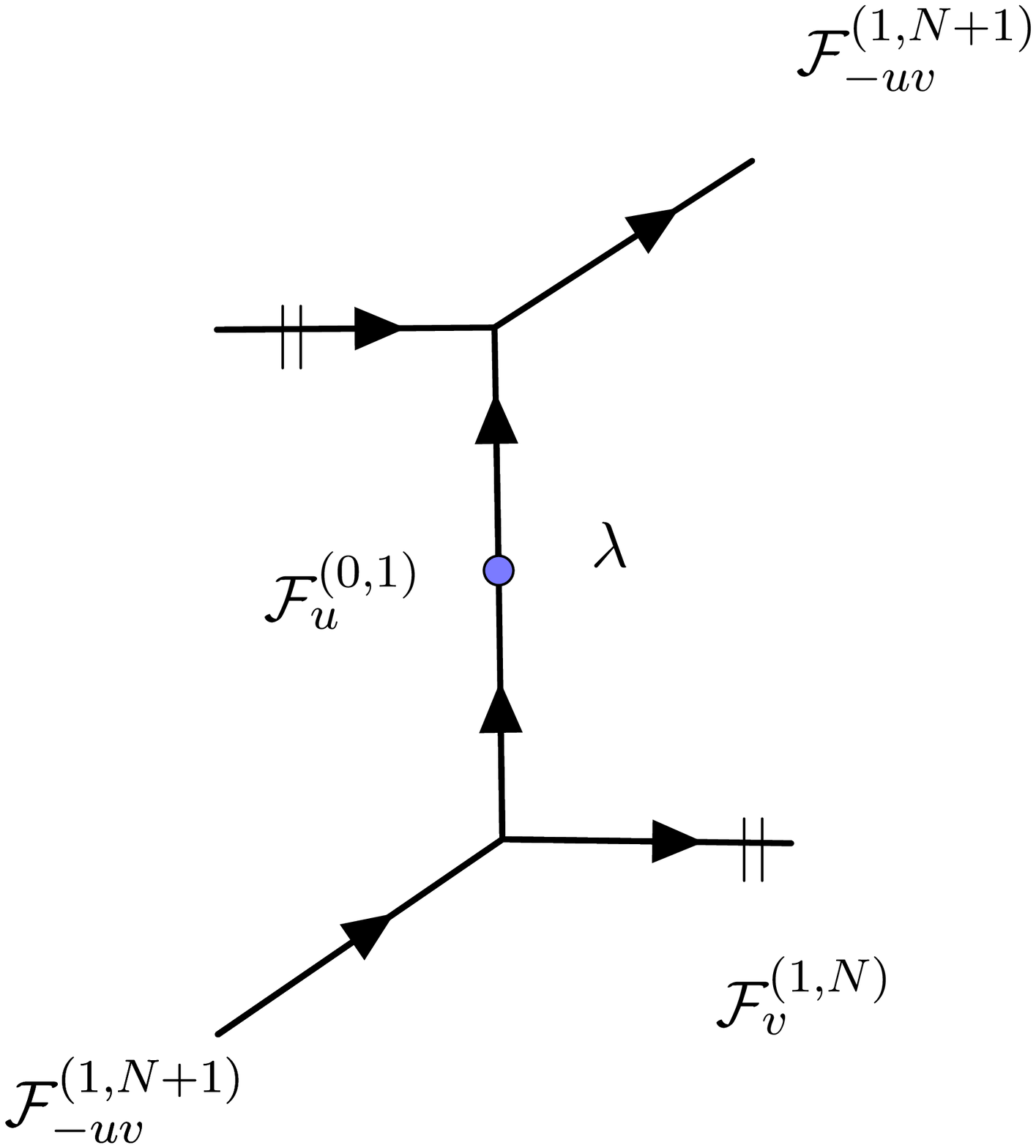}
\end{center}
\caption{Graphical expression of $\sum_\la\Phi^*_\la(u)\Phi_\la(u)$. The two horizontal lines with $||$ are glued together. }
 \lb{3-3vertex} 
\end{figure}

\subsection{The higher rank extension}
To extend $W_{p,p^*}(\Gamma(\widehat{A}_0))$ to the one associated with  the higher rank instantons,  
one need to take a composition of  $T(u)$'s. By using \eqref{PhilasPhila} and Propositions \ref{cE2AR} and \ref{cEcE}, one obtains the following expression. 
\bea
T(u_1)\cdots T(u_M)&=&{\cal C}_M\sum_{k=0}^\infty {\frak q}_M^{k}
\sum_{\la^{(1)},\cdots,\la^{(M)}\atop
\sum_j|\la^{(j)}|=k}\prod_{i,j=1}^M\frac{N_{\la^{(i)}\la^{(j)}}(pqu_{i,j}/t)}{N_{\la^{(i)}\la^{(j)}}(qu_{i,j}/t)}
\nn\\
&&\qquad\quad \times:\prod_{l=1}^M\prod_{\Box\in A(\la^{(l)})}Y(u_l/q^{\Box})
\prod_{\blacksquare\in R(\la^{(l)})}Y((q/t)u_l/q^{\blacksquare})^{-1}:,\lb{TT}
\ena
where we set $u_{j,i}=u_j/u_i$, and 
\bea
&&{\frak q}_M={\frak q}p^{-(M-1)}=p^{*-1}p^{M+N}(t/q)^{1/2},\\
&&{\cal C}_M=\left(\frac{(p^{-1}t;q,t,p)_\infty}{(q;q,t,p)_\infty}\right)^M
\prod_{1\leq i<j\leq M}
\frac{(p^{-1}tu_{j,i};q,t)_\infty(pqu_{j,i};q,t)_\infty}{(tu_{j,i};q,t)_\infty(qu_{j,i};q,t)_\infty}.
\ena
In fact, the sum
\bea
\chi_p({\frak M}_{k,M})&=&
\sum_{\la^{(1)},\cdots,\la^{(M)}\atop
\sum_j|\la^{(j)}|=k}\prod_{i,j=1}^M\frac{N_{\la^{(i)}\la^{(j)}}(pqu_{i,j}/t)}{N_{\la^{(i)}\la^{(j)}}(qu_{i,j}/t)}
\lb{chiyMkM}
\ena
coincides with the equivariant $\chi_y$ $(y=p)$ genus 
of the moduli space of rank $M$ instantons with charge $k$\cite{HIKLV}. 

Note that $T(u)$ gives a non-commutative  5d-analogue of 
Nekrasov's $qq$-character of the $\cN=2^*$ $U(1)$ theory\cite{KimPes}.   
We expect that $T(u_1)\cdots T(u_M)$ gives a non-commutative  5d-analogue of 
Nekrasov's $qq$-character of the $\cN=2^*$ $U(M)$ theory. 

\section{Instanton calculus in  the Jordan quiver gauge theories}\lb{InstPF}
By using the generating function $T(u)$ of $W_{p,p^*}(\Gamma(\widehat{A}_0))$ constructed in the previous section, we show some demonstrations of deriving the 
instanton partition functions  of the 5d and 6d lifts of the 4d $\cN=2^*$ SUSY gauge theories.  
Hence the  elliptic quantum toroidal algebra $\cU_{q,t,p}$ is a  relevant quantum group structure of dealing with such theories. 

\subsection{Instanton partition functions for the 5d and 6d lifts of the $\cN=2^*$ $U(1)$ theory}
As mentioned in $\S$\ref{genfW}, the vacuum expectation value of  $T(u)$ gives 
 the  instanton partition functions of the 5d  lift of the 4d $\cN=2^*$ $U(1)$ gauge theory.
\bea
\bra{0} T(u)\ket{0}&=&{\cal C}\,\sum_\la{\frak q}^{|\la|}{\cal Z}^{\widehat{A}_0}_\la(t,q^{-1},p).
\lb{chiygenus}
\ena
We then identify $T(u)=\sum_\la\Phi^*_\la(u)\Phi_\la(u)$  with the basic topological vertex depicted  in Fig.\ref{3-3vertex}.

An immediate next application  is to take the trace of $T(u)$. 
Let us introduce the degree counting operator
\be
&&d=-\sum_{m>0}\frac{m^2}{\kappa_m(1-(t/q)^m)}\frac{1-p^m}{1-p^{*m}}\al'_{-m}\al'_m
\en
such that
\be
&&[d,\al'_m]=m\al'_m\qquad m\in \Z_{\not=0}.
\en
Then the following trace  yields the 6d version of the partition function of the rank 1 instantions\cite{AK,HIV,HI,HIKLV}. 
\bea
&&\tr_{\F^{(1,N+1)}_{-uv}} Q^{d}\,T(u)
={\cal C}_Q\sum_\la  {\frak q}^{|\la|}{\cal Z}^{\widehat{A}_0}_\la(t,q^{-1},p;Q),
\lb{ellgenus}
\ena
where 
\bea
&&{\cal C}_Q=\frac{1}{(Q;Q)_\infty}\frac{(p^{-1}t;q,t,p)_\infty}{(q;q,t,p)_\infty}
\frac{(p^{-1}tQ;q,t,Q)_\infty(pqQ;q,t,Q)_\infty}{(tQ;q,t,Q)_\infty(qQ;q,t,Q)_\infty},\\
&&{\cal Z}^{\widehat{A}_0}_\la(t,q^{-1},p;Q)=
\frac{N^\theta_{\la\la}(pq/t;Q)}{N^\theta_{\la\la}(q/t;Q)}.
\ena
Here $N^\theta_{\la\mu}(x;Q)$ denotes the theta function analogue of the Nekrasov function given by
\bea
&&N^\theta_{\la\mu}(x;Q)=
\prod_{\Box\in \la}\theta_Q(xq^{-a_\mu(\Box)-1}t^{-\ell_\la(\Box)})
\prod_{\blacksquare\in \mu}\theta_Q(xq^{a_\la(\blacksquare)}t^{\ell_\mu(\blacksquare)+1}).
\lb{thetaNekrasov}
\ena
In fact  the sum  $\sum_{\la, |\la|=n}{\cal Z}^{\widehat{A}_0}_\la(t,q^{-1},p;Q)$ gives 
 the equivariant elliptic genus of ${\rm Hilb}_n(\C^2)$ \cite{LiLiuZhou}. 

\subsection{The 5d and 6d lifts of  the $\cN=2^*$ $U(M)$ theory}
The higher rank instanton partition functions can be obtained from the composition 
$T(u_1)\cdots T(u_M)$ in \eqref{TT}. 
The  vacuum expectation value gives the instanton partition function of the 5d lift 
of the 4d $\cN=2^*$ $U(M)$ theory
\cite{HIKLV}. 
\bea
\bra{0}T(u_1)\cdots T(u_M)\ket{0}&=&{\cal C}_M\sum_{k=0}^\infty {\frak q}_M^{k}\,
\chi_p({\frak M}_{k,M}),
\ena
where $\chi_p({\frak M}_{k,M})$ is given by \eqref{chiyMkM}. 

Furthermore taking the trace  of \eqref{TT}, one obtains
\bea
\tr_{\F^{(1,N+1)}_{-u_1v_1}} Q^{d}\,T(u_1)\cdots T(u_M)
&=&{\cal C}_{Q,M}\sum_{k=0}^\infty {\frak q}_M^{k}\, 
{\cal E}_{p,Q}({\frak M}_{k,M}),\lb{trTT}
\ena
where  $u_1v_1=u_2v_2=\cdots=u_Mv_M$ with arbitrary $v_1, \cdots, v_M\in \C^*$.  
We here also set
\bea
&&{\cal E}_{p,Q}({\frak M}_{k,M})=\sum_{\la^{(1)},\cdots,\la^{(M)}\atop
\sum_j|\la^{(j)}|=k}\prod_{i,j=1}^M
\frac{N^\theta_{\la^{(i)}\la^{(j)}}(pqu_{i,j}/t;Q)}
{N^\theta_{\la^{(i)}\la^{(j)}}(qu_{i,j}/t;Q)}
,\lb{ellgenus}\\
&&{\cal C}_{Q,M}=\frac{1}{(Q;Q)_\infty}
\left(\frac{(t;q,t)_\infty \Gamma_3(p^{-1}t;q,t,Q)}{(p^{-1}t;q,t)_\infty \Gamma_3(t;q,t,Q)}\right)^M
\prod_{1\leq i< j\leq M}\frac{\Gamma_3(p^{-1}tu_{j,i};q,t,Q)\Gamma_3(pqu_{j,i};q,t,Q)}{ \Gamma_3(tu_{j,i};q,t,Q)\Gamma_3(qu_{j,i};q,t,Q)},\nn\\&&
\ena
The  sum $\cE_{p,Q}({\frak M}_{k,M})$ gives the equivariant elliptic  genus of the moduli space of 
rank $M$ instantons with charge $k$\cite{HI,HKLV}. Hence  \eqref{trTT} gives the  instanton partition function of the 6d lift of the $\cN=2^*$ 
$U(M)$ theory.

\section{Correlation Functions of $\Phi(u)$ and $\Psi^*(v)$}
We here give some $\cU_{q,t,p}$-analogues of the formulas of those obtained in Propositions 5.1 and 5.2 in \cite{AFS}. We expect that  they  give the $\cN=2^*$ theory (Jordan quiver gauge theory) analogue of the partition function of  the  5d  lift of the pure $SU(N)$ theory 
and the 5d and 6d lifts of the $SU(N)$ theory with $2N$ fundamental matters, respectively.  The latter theories are the $A_1$ linear quiver gauge theories.

\subsection{OPE formulas}
From Propositions \ref{tildePhi}, \ref{tildePsi}  and  \ref{cEcE}, we have the following OPE formulas for the operator parts of the type I $\Phi(u)$ and the type II dual $\Psi^*(v)$ intertwining operators. 
\begin{prop}\lb{OPEs}
\be
&&\tPhi_\la(u)\tPhi_\mu(v) 
=\frac{\cG(v/u,p)}{{N}_{\mu\la}(v/u;p)}:\tPhi_\la(u)\tPhi_\mu(v) :,\\
&&\tPsi^*_\la(u)\tPsi^*_\mu(v)
=\frac{\cG^*(qv/tu,p)}{{N}^*_{\mu\la}(qv/tu;p)}:\tPsi^*_\la(u)\tPsi^*_\mu(v):,\\
&&\tPhi_\la(u)\tPsi^*_\mu(v)
=\frac{N_{\mu\la}((t/q)^{-1/2}v/u)}{\cG(v/u)}:\tPhi_\la(u)\tPsi^*_\mu(v)
:,\\
&&\tPsi^*_\la(u)\tPhi_\mu(v)
=\frac{N_{\mu\la}((t/q)^{-1/2}v/u)}{\cG(v/u)}:\tPsi^*_\la(u)\tPhi_\mu(v):.
\en
Here 
${N}_{\la\mu}(u;p)$ and ${N}^*_{\la\mu}(u;p)$  denote  $p$-deformations of $N_{\la\mu}(u)$ 
in \eqref{5dNekrasov} given by
\bea
{N}_{\la\mu}(u;p)&=&\prod_{\Box\in \la}\frac{(uq^{-a_\mu(\Box)-1}t^{-\ell_\la(\Box)};p)_\infty}
{(p^*uq^{-a_\mu(\Box)-1}t^{-\ell_\la(\Box)};p)_\infty}
\prod_{\blacksquare\in \mu}\frac{(uq^{a_\la(\blacksquare)}t^{\ell_\mu(\blacksquare)+1};p)_\infty}{(p^*uq^{a_\la(\blacksquare)}t^{\ell_\mu(\blacksquare)+1};p)_\infty},
\\
{N}^*_{\la\mu}(u;p)&=&N_{\la\mu}(u;p)\Biggr|_{p\leftrightarrow p^*}. 
\ena
We also set
\be
\cG(u)&=&\exp\left\{\sum_{n>0}\frac{1}{n}\frac{t^{n}}{(1-q^n)(1-t^{n})}u^n\right\}=\frac{1}{(tu;q,t)_\infty},\\
\cG(u,p)&=&\exp\left\{\sum_{n>0}\frac{1}{n}\frac{(1-p^{*n})t^n}{(1-q^n)(1-t^{n})(1-p^n)}u^n\right\}=\frac{(p^*tu;q,t,p)_\infty}{(tu;q,t,p)_\infty},\\
\cG^*(u,p)&=&\cG(u,p)\Biggr|_{p\leftrightarrow p^*}.
\en
\end{prop}

\subsection{The four points operator}
Remember
\be
&&\Phi(-x)\ :\ \F^{(1,N+1)}_{u}\ \to\ \F^{(0,1)}_{-x}\tot \F^{(1,N)}_{u/x},\\
&&\Psi^*(-x)\ :\  \F^{(1,N)}_{u}\tot\F^{(0,1)}_{-x} \ \to\ \F^{(1,N+1)}_{ux}.
\en 
Let us consider the following composition.
\be
&&\phi(w):=(\Psi^*(-w)\tot\id)\circ (\id\tot \Phi(-w))\ :\  \F^{(1,M)}_{v}\tot\F^{(1,L)}_{u} \ \to\ \F^{(1,M+1)}_{vw}\tot \F^{(1,L-1)}_{u/w}.
\en
This is an  intertwining operator satisfying
\be
&&\Delta(x)\phi(w)=\phi(w)\Delta(x)\qquad \forall x\in \cU_{q,t,p}. 
\en 
We call $\phi(w)$ the basic four points operator. 
 
\begin{prop}
The action of $\phi(w)$ on $\ket{\xi_v^M}\tot \ket{\eta_u^L}\in \F^{(1,M)}_{v}\tot\F^{(1,L)}_{u} $ is given by
\be
&&\phi(w)\ \ket{\xi_v^M}\tot \ket{\eta_u^L}
=\sum_{\la}\frac{c_\la(p^*)}{c'_\la(p^*)}\ \Psi^*_\la(-w)\ket{\xi_v^M}\tot \Phi_\la(-w) \ket{\eta_u^L}.
\en
\end{prop}
\noindent
{\it Proof.}\ Use \eqref{compPhis}, \eqref{ketlapketla}, \eqref{compPsis}  and \eqref{reltot}. \qed

\begin{prop}
The composition of $\phi(w)$'s  gives the following generalized  four points intertwining operator. 
\be
&&\phi(w_N)\circ \cdots\circ  \phi(w_1)\ :\   \F^{(1,M)}_{v}\tot\F^{(1,L)}_{u} \ \to\ \F^{(1,M+N)}_{vw_1\cdots w_N}\tot \F^{(1,L-N)}_{u/w_1\cdots w_N}.
\en
The action on $\ket{\xi_v^M}\tot \ket{\eta_u^L}\in \F^{(1,M)}_{v}\tot\F^{(1,L)}_{u} $ is given by
\be
&&\phi(w_N)\circ \cdots\circ  \phi(w_1)\ \ket{\xi_v^M}\tot \ket{\eta_u^L}\\
&&=\sum_{\la^{(1)},\cdots,\la^{(N)}}\prod_{a=1}^N\frac{c_{\la^{(a)}}(p^*)}{c'_{\la^{(a)}}(p^*)}\ 
\Psi^*_{\la^{(N)}}(-w_N)\cdots \Psi^*_{\la^{(1)}}(-w_1)\ket{\xi_v^M}\tot \Phi_{\la^{(N)}}(-w_N) \cdots \Phi_{\la^{(1)}}(-w_1)\ket{\eta_u^L}.
\en
\end{prop}

\begin{figure}[htbp]
\begin{center}
\includegraphics[
scale=0.3, 
height=70mm, width=100mm]{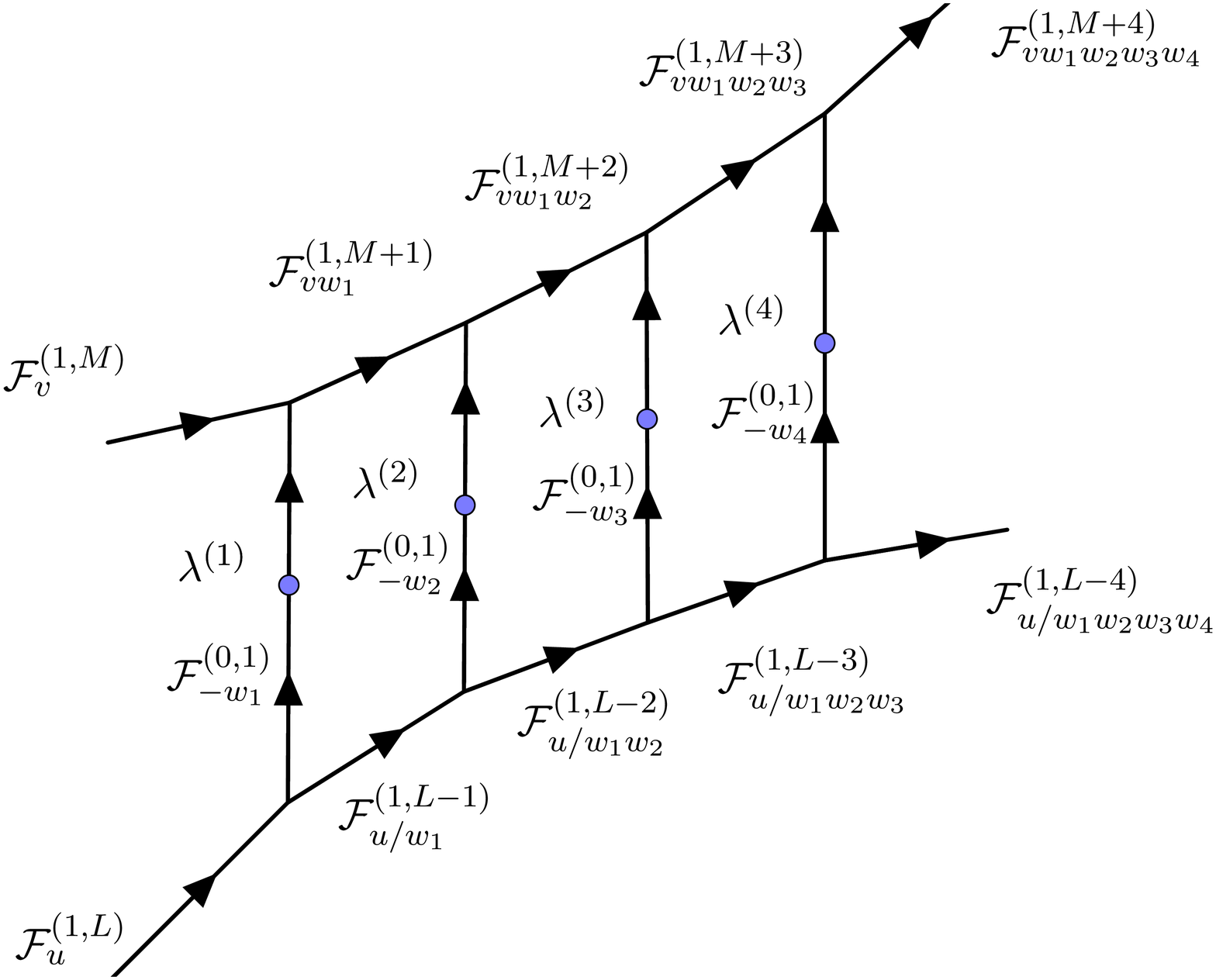}
\end{center}
\caption{Graphical expression of $\phi(w_4)\circ \cdots\circ  \phi(w_1)$}
 \lb{4-pt op} 
\end{figure}

\begin{cor}
For a vector $\ket{\omega^{M+N}_{vw_1\cdots w_N}}\tot\ket{\zeta^{L-N}_{u/w_1\cdots w_N}}
\in \F^{(1,M+N)}_{vw_1\cdots w_N}\tot \F^{(1,L-N)}_{u/w_1\cdots w_N}$, one has the 
following expectation value.
\be
&&\bra{\omega^{M+N}_{vw_1\cdots w_N}}\tot\bra{\zeta^{L-N}_{u/w_1\cdots w_N}}
\ \phi(w_N)\circ \cdots\circ  \phi(w_1)\ \ket{\xi_v^M}\tot \ket{\eta_u^L}\\
&&=\sum_{\la^{(1)},\cdots,\la^{(N)}}\prod_{a=1}^N\frac{c_{\la^{(a)}}(p^*)}{c'_{\la^{(a)}}(p^*)}\ 
\bra{\omega^{M+N}_{vw_1\cdots w_N}}\Psi^*_{\la^{(N)}}(-w_N)\cdots \Psi^*_{\la^{(1)}}(-w_1)\ket{\xi_v^M}\\
&&\qquad\qquad\qquad\qquad\qquad\qquad \tot 
\bra{\zeta^{L-N}_{u/w_1\cdots w_N}}\Phi_{\la^{(N)}}(-w_N) \cdots \Phi_{\la^{(1)}}(-w_1)\ket{\eta_u^L}.
\en
\end{cor}
Note that this  gives the $\cN=2^*$ theory analogue of the partition function of  the  5d pure $SU(N)$ theory obtained, for example, in \cite{AK,AFS}.

\subsection{The six points  and higher operators}

Next let us consider the operator 
\be
&&\phi\left(w; \mmatrix{y\\[-2mm] x\cr}
\right)=(\Phi(-y)\tot \id)\circ(\Psi^*(-w) \tot \Psi^*(-x))\circ (\id \tot \Phi(-w)\tot \id)\\
&&\qquad\qquad\qquad  :\  \F^{(1,0)}_{v}\tot\F^{(1,0)}_{u}\tot \F^{(0,1)}_{-x} \ \to\ \F^{(0,1)}_{-y}\tot\F^{(1,0)}_{vw/y}\tot \F^{(1,0)}_{ux/w}.
\en
This is an intertwining operator satisfying
\be
&&\phi\left(w; \mmatrix{y\\[-2mm] x\cr}\right)(\Delta\tot\id)\Delta(a)
=(\Delta\tot\id)\Delta(a)\ \phi\left(w; \mmatrix{y\\[-2mm] x\cr}\right) \qquad \forall a\in \cU_{q,t,p}.
\en
We call $\Tphi{w}{y}{x}$ the six points operator. 
 
\begin{prop}
For $\mu\in \cP^+$, the action of $\Tphi{w}{y}{x}$ on $\ket{\xi_v^0}\tot \ket{\eta_u^0}\tot \ket{\mu}_{-x}\in \F^{(1,0)}_{v}\tot\F^{(1,0)}_{u}\tot\F^{(0,1)}_{-x} $ is given by
\be
&&\Tphi{w}{y}{x}\ \ket{\xi_v^0}\tot \ket{\eta_u^0}\tot \ket{\mu}_{-x}
=\sum_{\la,\sigma\in \cP^+}\frac{c_\la(p^*)}{c'_\la(p^*)}\ 
\ket{\sigma}'_{-y}\tot \Phi_{\sigma}(-y)\Psi^*_\la(-w)\ket{\xi_v^0}\tot 
\Psi^*_\mu(-x)\Phi_{\la}(-w) \ket{\eta_u^0}.
\en
\end{prop}
Let us set
\be
&&\eTphi{j}{w}{y}{x}:=\underbrace{\id\tot\cdots\tot\id}_{j-1}\tot \Tphi{w}{y}{x}\tot\underbrace{\id\tot\cdots\tot\id}_{N-j}.
\en
\begin{prop}
By composing $\eTphi{j}{w_j}{y_j}{x_j}\ (j=1,\cdots,N)$  we obtain 
the following $2(N+2)$ points intertwining operator. 
\be
&&\eTphi{N}{w_N}{y_N}{x_N}\circ \cdots\circ  \eTphi{1}{w_1}{y_1}{x_1}\\
&&\qquad  :\ 
\F^{(1,0)}_{v}\tot\F^{(1,0)}_{u}\tot \F^{(0,1)}_{-x_1}\tot \cdots \tot\F^{(0,1)}_{-x_N} \ \to\ \F^{(0,1)}_{-y_1}\tot \cdots \tot \F^{(0,1)}_{-y_N}\tot\F^{(1,0)}_{v\frac{w_1\cdots w_N}{y_1\cdots y_N}
}\tot \F^{(1,0)}_{u\frac{x_1\cdots x_N}{w_1\cdots w_N}
}.
\en
For $\mu=(\mu^{(1)},\cdots,\mu^{(N)})\in (\cP^+)^N$, the action on $\ket{\xi^0_{v}}\tot \ket{\eta^0_u}\tot \ket{\mu^{(1)}}_{-x_1}\tot \cdots \tot \ket{\mu^{(N)}}_{-x_N}$\\ $\in \F^{(1,0)}_{v}\tot\F^{(1,0)}_{u}\tot\F^{(0,1)}_{-x_1}\tot \cdots \tot\F^{(0,1)}_{-x_N} $
is given by
\be
&&\eTphi{N}{w_N}{y_N}{x_N}\circ \cdots\circ  \eTphi{1}{w_1}{y_1}{x_1}
\ket{\xi^0_{v}}\tot \ket{\eta^0_u}\tot  \ket{\mu^{(1)}}_{-x_1}\tot \cdots \tot \ket{\mu^{(N)}}_{-x_N}\\
&&=\sum_{\la^{(1)},\cdots,\la^{(N)}\atop
\sigma^{(1)},\cdots,\sigma^{(N)}}\prod_{a=1}^N\frac{c_{\la^{(a)}}(p^*)}{c'_{\la^{(a)}}(p^*)}\ 
\ket{\sigma^{(1)}}'_{-y_1}\tot\cdots\tot \ket{\sigma^{(N)}}'_{-y_N}\\
&&\hspace{3cm}\tot \Phi_{\sigma^{(N)}}(-y_N)\Psi^*_{\la^{(N)}}(-w_N)\cdots \Phi_{\sigma^{(1)}}(-y_1)\Psi^*_{\la^{(1)}}(-w_1)\ket{\xi^0_v}\\
&&\hspace{3cm}\tot 
\Psi^*_{\mu^{(N)}}(-x_N)\Phi_{\la^{(N)}}(-w_N)\cdots \Psi^*_{\mu^{(1)}}(-x_1)\Phi_{\la^{(1)}}(-w_1) \ket{\eta^0_u},
\en
where $\la=(\la^{(1)},\cdots,\la^{(N)}), \sigma=(\sigma^{(1)},\cdots,\sigma^{(N)})\in (\cP^+)^N$. 
\end{prop}

\begin{figure}[htbp]
\begin{center}
\includegraphics[
scale=0.3, 
height=73mm, width=123mm]{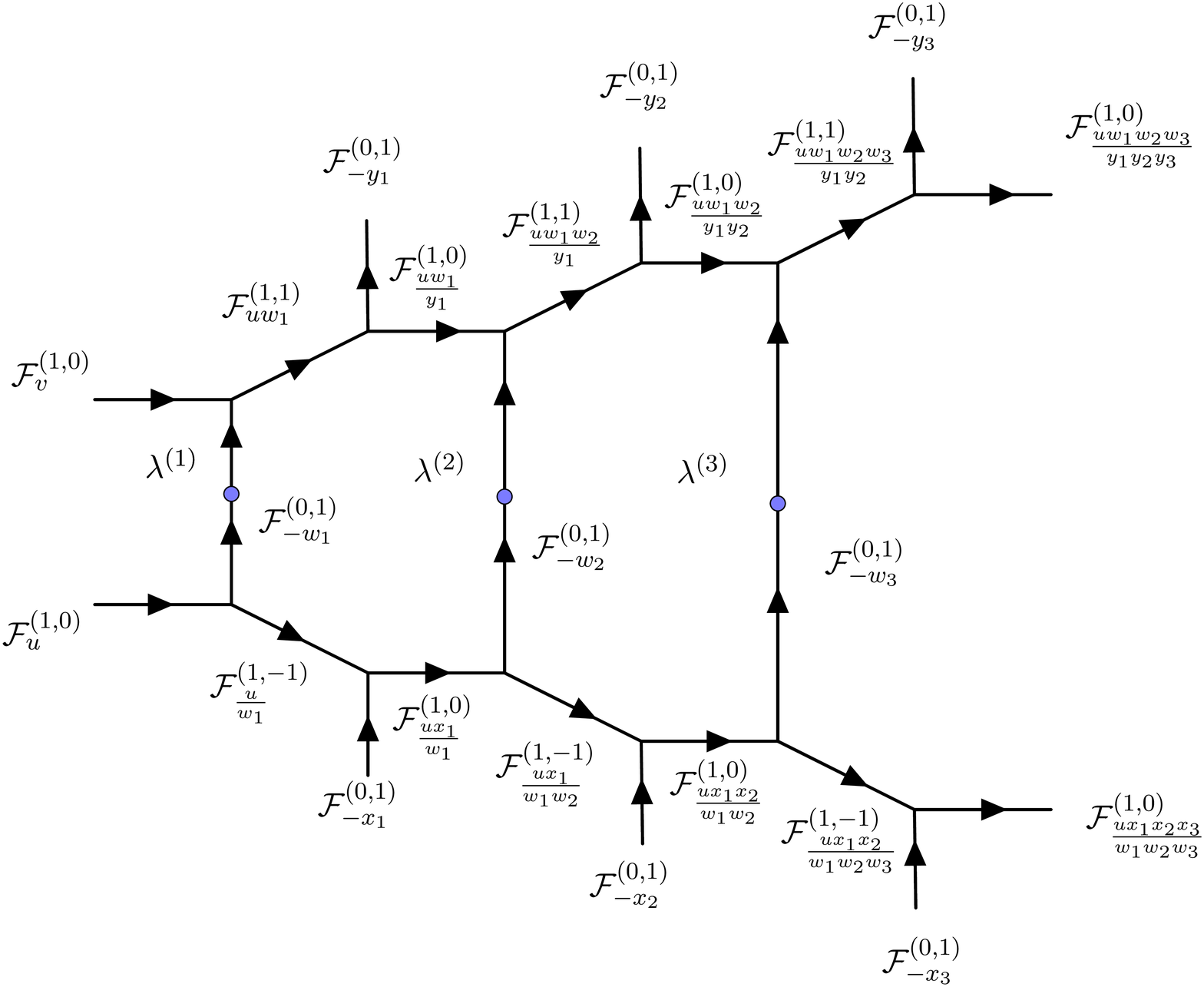}
\end{center}
\caption{Graphical expression of 
$\eTphi{3}{w_3}{y_3}{x_3}\circ \eTphi{2}{w_2}{y_2}{x_2}\circ  \eTphi{1}{w_1}{y_1}{x_1}$}
 \lb{10-pt op} 
\end{figure}

\begin{cor}\lb{expvalue}
 For a vector 
\be
&&\ket{\nu^{(1)}}'_{-y_1}\tot \cdots\tot \ket{\nu^{(N)}}'_{-y_N}\tot 
\ket{\omega^0_{v\frac{w_1\cdots w_N}{y_1\cdots y_N}}}\tot  
\ket{\zeta^0_{u\frac{x_1\cdots x_N}{w_1\cdots w_N}}}\in 
\F^{(0,1)}_{-y_1}\tot \cdots \tot \F^{(0,1)}_{-y_N}\tot\F^{(1,0)}_{v\frac{w_1\cdots w_N}{y_1\cdots y_N}}\tot \F^{(1,0)}_{u\frac{x_1\cdots x_N}{w_1\cdots w_N}}, 
\en
 we have the expectation value 
\be
&&
\hspace{-1cm}
{}_{-y_1}`\bra{\nu^{(1)}}\tot\cdots\tot {}_{-y_N}`\bra{\nu^{(N)}}\tot 
\bra{\omega^0_{v\frac{w_1\cdots w_N}{y_1\cdots y_N}}}\tot  
\bra{\zeta^0_{u\frac{x_1\cdots x_N}{w_1\cdots w_N}}}
\eTphi{N}{w_N}{y_N}{x_N}\circ \cdots\circ  \eTphi{1}{w_1}{y_1}{x_1}\\
&&\hspace{3cm}\qquad\qquad\qquad\qquad\qquad\qquad \times \ket{\xi^0_{v}}\tot \ket{\eta^0_u}\tot  
\ket{\mu^{(1)}}_{-x_1}\tot \cdots \tot \ket{\mu^{(N)}}_{-x_N}\\
&&\hspace{-1cm}=\sum_{\la^{(1)},\cdots,\la^{(N)}\atop
}\prod_{a=1}^N\frac{c_{\la^{(a)}}(p^*)}{c'_{\la^{(a)}}(p^*)}
\ \bra{\omega^0_{v\frac{w_1\cdots w_N}{y_1\cdots y_N}}}\Phi_{\nu^{(N)}}(-y_N)\Psi^*_{\la^{(N)}}(-w_N)\cdots \Phi_{\nu^{(1)}}(-y_1)\Psi^*_{\la^{(1)}}(-w_1)\ket{\xi^0_v}\\
&&\hspace{3cm}\tot \bra{\zeta^0_{u\frac{x_1\cdots x_N}{w_1\cdots w_N}}}
\Psi^*_{\mu^{(N)}}(-x_N)\Phi_{\la^{(N)}}(-w_N)\cdots \Psi^*_{\mu^{(1)}}(-x_1)\Phi_{\la^{(1)}}(-w_1) \ket{\eta^0_u}.
\en
Here ${}_y`\bra{\mu}$\  $(\mu\in \cP^+, y\in \C^*)$ is defined by
\be
&&{}_y`\bra{\mu}\ket{\nu}'_y=\delta_{\mu,\nu}\qquad  \nu\in \cP^+
\en
\end{cor}
\begin{cor}\lb{tracevalue}
When $x_1\cdots x_N=y_1\cdots y_N=w_1\cdots w_N$, one can take the following trace 
\be
&&\hspace{-2cm}\tr_{\F^{(1,0)}_v}Q^d\tot\ \tr_{\F^{(1,0)}_u}Q^d\left(
{}_{-y_1}`\bra{\nu^{(1)}}\tot\cdots\tot {}_{-y_N}`\bra{\nu^{(N)}}
\eTphi{N}{w_N}{y_N}{x_N}\circ \cdots\circ  \eTphi{1}{w_1}{y_1}{x_1}
\ket{\mu^{(1)}}_{-x_1}\tot \cdots \tot \ket{\mu^{(N)}}_{-x_N}
\right)\\
&&=\sum_{\la^{(1)},\cdots,\la^{(N)}
}\prod_{a=1}^N\frac{c_{\la^{(a)}}(p^*)}{c'_{\la^{(a)}}(p^*)}
\ \tr_{\F^{(1,0)}_v}\left(Q^d\Phi_{\nu^{(N)}}(-y_N)\Psi^*_{\la^{(N)}}(-w_N)\cdots \Phi_{\nu^{(1)}}(-y_1)\Psi^*_{\la^{(1)}}(-w_1)\right)\\
&&\hspace{4cm}\tot\ \tr_{\F^{(1,0)}_u}\left(Q^d
\Psi^*_{\mu^{(N)}}(-x_N)\Phi_{\la^{(N)}}(-w_N)\cdots \Psi^*_{\mu^{(1)}}(-x_1)\Phi_{\la^{(1)}}(-w_1) \right).
\en
\end{cor}
Note that Corollary \ref{expvalue} gives the $\cN=2^*$ theory analogue of the partition function of  the 5d $SU(N)$ theory with $2N$ fundamental matters  obtained, for example, in \cite{AK,AFS} for the case $\mu=\nu=(\emptyset,\cdots, \emptyset)$.  
The trace in Corollary \ref{tracevalue} gives its 6d lift.  It can be 
evaluated as follows. 
\begin{prop}
\be
&&\hspace{-2cm}\tr_{\F^{(1,0)}_v}Q^d\tot\ \tr_{\F^{(1,0)}_u}Q^d\left(
{}_{-y_1}\bra{\nu^{(1)}}\tot\cdots\tot {}_{-y_N}\bra{\nu^{(N)}}
\eTphi{N}{w_N}{y_N}{x_N}\circ \cdots\circ  \eTphi{1}{w_1}{y_1}{x_1}
\ket{\mu^{(1)}}'_{-x_1}\tot \cdots \tot \ket{\mu^{(N)}}'_{-x_N}
\right)\\
&&={\cal C}_N(q,t,p,Q){\cN}^N_{\mu}(u,x,w;q,t,p,Q){\cN}^N_{\nu}(v,y,w;q,t,p,Q)
\\
&&\times \prod_{1\leq a<b\leq N} \left(\frac{\Gamma_4(qw_a/w_b;q,t,p^{*},Q)}{\Gamma_4(p^{*}tw_a/w_b;q,t,p^{*},Q)}
\frac{\Gamma_4(tw_a/w_b;q,t,p^*,Q)}{\Gamma_4(p^*qw_a/w_b;q,t,p^*,Q)}\right)\\
&&\times\sum_{\la^{(1)},\cdots,\la^{(N)}}\prod_{a=1}^N
\left(-\frac{v}{u}\prod_{k=1}^{a-1}\frac{ w_k^2}{x_ky_k}
\right)^{|\la_a|}q^{-n(\la^{(a)'})}t^{n(\la^{(a)})}
\\
&&\times\prod_{1\leq a<b\leq N} {N^\theta_{\la^{(a)}\la^{(b)}}(w_a/w_b;p^*)}{N^\theta_{\la^{(a)}\la^{(b)}}(qw_a/tw_b;p^*)}\\
&&\times\prod_{a,b=1}^N \left(
N^\theta_{\la^{(b)}\nu^{(a)}}(\sqrt{q/t}w_b/y_a;Q)
N^\theta_{\mu^{(b)}\la^{(a)}}(\sqrt{q/t}x_b/w_a;Q)
N^{\Gamma}_{\la^{(a)}\la^{(b)}}(w_a/w_b;p^*,Q)N^{\Gamma}_{\la^{(a)}\la^{(b)}}(qw_a/tw_b;p^*,Q)
\right).
\en
Here we set $p^{**}=p^*(q/t)=p(q/t)^2$ and 
\be
&&{\cal C}_N(q,t,p,Q)\\
&&=\frac{1}{(Q;Q)^2_\infty}
\left(
\frac{(p^*tQ;q,t,p,Q)_\infty}{(tQ;q,t,p,Q)_\infty}\frac{(pqQ;q,t,p^*,Q)_\infty}{(qQ;q,t,p^*,Q)_\infty}\frac{(p^{**}tQ;q,t,p^*,Q)_\infty}{(tQ;q,t,p^*,Q)_\infty}
\frac{(p^{**}tQ;q,t,p^{**},Q)_\infty}{(qQ;q,t,p^{**},Q)_\infty}\right)^N,
\en
\be
&&{\cN}^N_{\mu}(u,x,w;q,t,p,Q)\\
&&=\prod_{a=1}^N
\left(-t^{-1}u\prod_{k=1}^a\frac{x_k}{w_k}\right)^{|\mu^{(a)}|}
t^{-n(\mu^{(a)})}c'_{\mu^{(a)}}(p^{**})
\\
&&\qquad\times 
\prod_{1\leq a<b\leq N} \Gamma(\sqrt{qt}x_a/w_b;q,t)_\infty\frac{\Gamma_4(qx_a/x_b;q,t,p^{**},Q)}{\Gamma_4(p^{**}x_a/x_b;q,t,p^{**},Q)}
{N^\theta_{\mu^{(a)}\mu^{(b)}}(x_a/x_b;p^{**})}
\\
&&\qquad\times 
\prod_{a,b=1}^N{N^\Gamma_{\mu^{(a)}\mu^{(b)}}(qx_a/tx_b;p^{**},Q)}{\Gamma_3(\sqrt{qt}x_a/w_b;q,t,Q)},
\\
&&{\cN}^N_{\nu}(v,y,w;q,t,p,Q)\\
&&=\prod_{a=1}^N
\left(-v^{-1}\prod_{k=1}^a\frac{y_k}{w_k}
\right)^{|\nu^{(a)}|}
t^{-n(\nu^{(a)})}
c'_{\nu^{(a)}}(p)
\\
&&\qquad\times 
\prod_{1\leq a<b\leq N} \Gamma(\sqrt{qt}w_b/y_a;q,t) \frac{\Gamma_4(ty_a/y_b;q,t,p,Q)}{\Gamma_4(pqy_a/y_b;q,t,p,Q)}
{N^\theta_{\nu^{(a)}\nu^{(b)}}(qy_a/ty_b;p)}
\\
&&\qquad\times 
\prod_{a,b=1}^N {N^\Gamma_{\nu^{(a)}\nu^{(b)}}(y_a/y_b;p,Q)}{\Gamma_3(\sqrt{qt}w_a/y_b;q,t,Q)}.
\en
Here 
$N^\Gamma_{\la\mu}(x;p,Q)$ denotes 
the elliptic Gamma- analogues of the Nekrasov function defined by 
\bea
&&N^\Gamma_{\la\mu}(x;p,Q)=
\prod_{\Box\in \la}\Gamma(xq^{-a_\mu(\Box)-1}t^{-\ell_\la(\Box)};p,Q)
\prod_{\blacksquare\in \mu}\Gamma(xq^{a_\la(\blacksquare)}t^{\ell_\mu(\blacksquare)+1};p,Q).
\ena
\end{prop}
\noindent
{\it Proof.}\ The statement follows from Proposition \ref{OPEs} and the trace formula given, for example, in Appendix E of \cite{KonnoBook}. 
\qed

As for the functions $N^\theta_{\la\mu}(x;p)$ and $N^\Gamma_{\la\mu}(x;p,Q)$, the following properties are useful.
\begin{prop}
\be
N^\theta_{\la\mu}(x;p)&=&\left(x\sqrt{\frac{t}{q}}\right)^{|\la|+|\mu|}\frac{f_\la(q,t)}{f_\mu(q,t)}
N^\theta_{\mu\la}(q/tx;p),\\
N^\Gamma_{\la\mu}(x;p,Q)&=&\frac{1}{N^\Gamma_{\mu\la}(pQq/tx;p,Q)}\\
&=&\left(x\sqrt{\frac{t}{q}}\right)^{|\la|+|\mu|}\frac{f_\la(q,t)}{f_\mu(q,t)}
\frac{1}{N^\theta_{\la\mu}(x;p)N^\theta_{\la\mu}(x;Q)N^\Gamma_{\mu\la}(q/tx;p,Q)}.
\en
\end{prop}

\section*{Acknowledgements}
The authors would like to thank Tatsuyuki Hikita, Michio Jimbo, Hiroaki Kanno, Taro Kimura, Evgeny Mukhin, Hiraku Nakajima, Yoshihisa Saito, Yasuhiko Yamada, Shintaro Yanagida, Yutaka Yoshida, Junichi Shiraishi, Andrey Smirnov for useful discussions. 
 HK is supported by the Grant-in-Aid for Scientific Research (C) 20K03507 JSPS, Japan.
KO is supported by the Grant-in-Aid for Scientific Research (C) 21K03191 JSPS, Japan. 


\vspace{1cm}
\noindent
{\bf\Large Declarations} \\[1mm]
{\bf Conflict of interest}\quad The authors declare that they have no conflict of interest.

\vspace{2cm}

\appendix
\setcounter{equation}{0}
\begin{appendix}



\section{Formulas for $E^\pm(\alpha,z)$ and $E^\pm(\alpha',z)$  }\lb{Zalg}

\begin{lem}\lb{lem:sec3-2}
The operators $E^\pm(\alpha,z)$ and $E^\pm(\alpha',z)$ defined in \eqref{Epm} and \eqref{Epmp} satisfy the following relations.
\begin{align}
& [\alpha_{-n}, E^+(\alpha,z)]=\frac{\kappa_n}{n}\frac{1-p^n}{1-p^{*n}}\gamma^{-n/2}z^{-n}E^+(\alpha,z) \quad (n>0), \lb{alEp}\\
& [\alpha_n, E^-(\alpha,z)]=\frac{\kappa_n}{n}\frac{1-p^n}{1-p^{*n}} \gamma^{-3n/2}z^n E^-(\alpha,z) \quad (n>0), \\
& [\alpha_{-n}, E^+(\alpha',z)]=-\frac{\kappa_n}{n} \gamma^{n/2} z^{-n}E^+(\alpha',z) \quad (n>0), \\
& [\alpha_n, E^-(\alpha',z)]=-\frac{\kappa_n}{n}\gamma^{-n/2}z^n E^-(\alpha',z) \quad (n>0), 
\lb{alEm}\\
& E^+(\alpha,z)E^-(\alpha,w)=g(w/z)^{-1}g(w/z;\gamma^2)^{-1}g(w/z;p^*)^{-1}E^-(\alpha,w)E^+(\alpha,z), \\
& E^+(\alpha',z)E^-(\alpha',w)=g(w/z;\gamma^2)^{-1}g(w/z;p)E^-(\alpha',w)E^+(\alpha',z), \\
& E^+(\alpha,z)E^-(\alpha',w)=g(\gamma w/z)g(\gamma w/z;\gamma^2) E^-(\alpha',w)E^+(\alpha,z), \\
& E^+(\alpha',z)E^-(\alpha,w)=g(\gamma w/z)g(\gamma w/z;\gamma^2) E^-(\alpha,w)E^+(\alpha',z), \\
& E^+(\alpha, z)x^+(w)=g( w/z)g(w/z;\gamma^2) g(w/z;p^*)x^+(w)E^+(\alpha,z), \\
& E^+(\alpha, z)x^-(w)=g(\gamma w/z)^{-1}g(\gamma w/z;\gamma^2)^{-1} x^-(w)E^+(\alpha,z), \\
& E^-(\alpha, z)x^+(w)=g( z/w)^{-1} g(z/w;\gamma^2)^{-1} g(z/w;p^*)^{-1} x^+(w)E^-(\alpha,z), \\
& E^-(\alpha, z)x^-(w)=g(\gamma z/w)g(\gamma z/w;\gamma^2)  x^-(w)E^-(\alpha,z), \\
& E^+(\alpha', z)x^+(w)=g(\gamma w/z)^{-1}g(\gamma w/z;\gamma^2)^{-1}  x^+(w)E^+(\alpha',z), \\
& E^+(\alpha', z)x^-(w)=g(w/z;\gamma^2)g(w/z;p)^{-1} f(p w/z;p) x^-(w)E^+(\alpha',z), \\
& E^-(\alpha', z)x^+(w)=g(\gamma z/w)g(\gamma z/w;\gamma^2)  x^+(w)E^-(\alpha',z), \\
& E^-(\alpha', z)x^-(w)=g( z/w;\gamma^2)^{-1} g( z/w;p) x^-(w)E^-(\alpha',z).
\end{align}
Here $g(z;s), \ s=p, p^*, \gamma^2$ and $g(z)$ are given in \eqref{gpm} and \eqref{def:g}, respectively.
\end{lem}

\noindent
{\it Proof.}\ The statement follows from  direct calculations using \eqref{alal}-\eqref{euqg7-2} 
and \eqref{alpalp}.
\qed

\section{Direct Check  of the Level (0,1) Representation} \lb{proofVerticalRep}
Let us check the relation \eqref{ellrelxpxm}. 
\bea
x^+(z)x^-(w)\ket{\la}_u&=&a^+(p)a^-(p)(t/q)^{-1/2}\sum_{i=1}^\ell
\sum_{j=1\atop j\not=i}^{\ell+1}\delta(q^{-1}u_i/w)\delta(u_j/z)\tA^{-}_{\la,i}(p)\tA^+_{\la-\bo_i,j}(p)\ket{\la-\bo_i+\bo_j}_u\nn\\
&&+a^+(p)a^-(p)(t/q)^{-1/2}\sum_{i=1}^\ell\delta(q^{-1}u_i/w)\delta(z/w)\tA^{-}_{\la,i}(p)\tA^+_{\la-\bo_i,i}(p)\ket{\la}_u,\lb{xpxm}\\
x^-(w)x^+(z)\ket{\la}_u
&=&a^+(p)a^-(p)(t/q)^{-1/2}\sum_{i=1}^\ell\sum_{j=1\atop j\not=i}^{\ell+1}
\delta(q^{-1}u_i/w)\delta(u_j/z)\tA^{+}_{\la,j}(p)\tA^-_{\la+\bo_j,i}(p)\ket{\la+\bo_j-\bo_i}_u\nn\\
&&+a^+(p)a^-(p)(t/q)^{-1/2}\sum_{i=1}^{\ell+1}\delta(u_i/w)\delta(z/w)\tA^{+}_{\la,i}(p)\tA^-_{\la+\bo_i,i}(p)\ket{\la}_u.\lb{xmxp}
\ena
One finds that the first terms in \eqref{xpxm} and \eqref{xmxp} coincide. Hence
\be
&&[x^+(z),x^-(w)]\ket{\la}_u\\
&&=a^+(p)a^-(p)(t/q)^{-1/2}\delta(z/w)\\
&&\qquad\times \left(\sum_{i=1}^\ell\delta(q^{-1}u_i/w)\delta(z/w)\tA^{-}_{\la,i}(p)\tA^+_{\la-\bo_i,i}(p)
-\sum_{i=1}^{\ell+1}\delta(u_i/w)\tA^{+}_{\la,i}(p)\tA^-_{\la+\bo_i,i}(p)\right)\ket{\la}_u.
\en
One also shows 
\be
a^+(p)a^-(p)&=&\frac{(1-q)(1-t^{-1})}{1-q/t}\frac{\theta_p(q/t)\theta_p(t)}{(p;p)_\infty^2\theta_p(q)},\\
\tA^{-}_{\la,i}(p)\tA^+_{\la-\bo_i,i}(p)&=&\prod_{j=1\atop \not=i}^\ell\frac{\theta_p(qu_j/tu_i)}{\theta_p(u_j/u_i)}\prod_{j=1\atop \not=i}^{\ell+1}\frac{\theta_p(tu_j/u_i)}{\theta_p(qu_j/u_i)}=(t/q)
\prod_{j=1\atop \not=i}^\ell\frac{\theta_p(tu_i/qu_j)}{\theta_p(u_i/u_j)}\prod_{j=1\atop \not=i}^{\ell+1}\frac{\theta_p(u_i/tu_j)}{\theta_p(u_i/qu_j)}\\
&=&-(t/q)\frac{(p;p)_\infty^2\theta_p(q)}{\theta_p(q/t)\theta_p(t)}\mathrm{Res}_{z=q^{-1}u_i}B^-(z/u;p)\frac{dz}{z},\\
\tA^{+}_{\la,i}(p)\tA^-_{\la+\bo_i,i}(p)&=&\prod_{j=1\atop \not=i}^\ell
\frac{\theta_p(u_j/tu_i)}{\theta_p(u_j/qu_i)}\prod_{j=1\atop \not=i}^{\ell+1}\frac{\theta_p(tu_j/qu_i)}{\theta_p(u_j/u_i)}=(t/q)
\prod_{j=1\atop \not=i}^\ell\frac{\theta_p(tu_i/u_j)}{\theta_p(qu_i/u_j)}\prod_{j=1\atop \not=i}^{\ell+1}\frac{\theta_p(qu_i/tu_j)}{\theta_p(u_i/u_j)}
\\
&=&(t/q)\frac{(p;p)_\infty^2\theta_p(q)}{\theta_p(q/t)\theta_p(t)}\mathrm{Res}_{z=u_i}B^-(z/u;p)\frac{dz}{z}.
\en
Hence
\be
&&[x^+(z),x^-(w)]\ket{\la}_u\\
&&=\frac{(1-q)(1-t^{-1})}{1-q/t}\delta(z/w)
\left(-\sum_{i=1}^\ell\delta(q^{-1}u_i/w)
\mathrm{Res}_{z=q^{-1}u_i}(t/q)^{1/2}B^-(z/u;p)\frac{dz}{z}\right.\\
&&\left.\quad \qquad\qquad\qquad\qquad\qquad\qquad -\sum_{i=1}^{\ell+1}\delta(u_i/w)
\mathrm{Res}_{z=u_i}(t/q)^{1/2}B^-(z/u;p)\frac{dz}{z}
\right)\ket{\la}_u.
\en
Using the partial fraction expansion formula \eqref{PFE},  
one can show
\be
&&B^-_\la(z/u;p)|_{|u/z|<1}-B^-_\la(z/u;p)|_{|u/z|>1}\\
&&=-\sum_{k=1}^\ell\delta(q^{-1}u_k/z)\mathrm{ Res}_{z=q^{-1}u_k}B^-_\la(z/u;p)\frac{dz}{z}
-\sum_{k=1}^{\ell+1}\delta(u_k/z)\mathrm{ Res}_{z=q^{-1}u_k}B^-_\la(z/u;p)\frac{dz}{z}
\en
Therefore we have
\be
&&[x^+(z),x^-(w)]\ket{\la}_u\\
&&=\frac{(1-q)(1-t^{-1})}{1-q/t}\delta(z/w)
\left((t/q)^{1/2}B^-_\la(z/u;p)|_{|u/z|<1}-(t/q)^{1/2}B^-_\la(z/u;p)|_{|u/z|>1}
\right)\ket{\la}_u\\
&&=\frac{(1-q)(1-t^{-1})}{1-q/t}\delta(z/w)
\left((t/q)^{1/2}(t/q)^{-1}B^+_\la(u/z;p)|_{|u/z|<1}-(t/q)^{1/2}B^-_\la(z/u;p)|_{|u/z|>1}
\right)\ket{\la}_u\\
&&=\frac{(1-q)(1-t^{-1})}{1-q/t}\delta(z/w)
\left(\psi^+(z)-\psi^-(z)\right)\ket{\la}_u.
\en
\qed

\section{Inductive Derivation of Theorem \ref{actionUqtpgl1}}\lb{IndDer}
We consider a tensor product of the vector representations in Proposition \ref{vecrep}.  
Define 
\be
&&V^{(N)}(u)=V(u)\tot V(u(t/q)^{-1})\tot V(u(t/q)^{-2})\tot\cdots\tot V(u(t/q)^{-N+1}).
\en
Set 
\be
&&{\cal P}^{(N)}=\{\la=(\la_1,\la_2,\cdots,\la_N)\in \Z^N\ |\  \la_1\geq \la_2\geq \cdots \geq \la_N\ \},\\
&&\ket{\la}^{(N)}_u=[u]_{\la_1}\tot [u(t/q)^{-1}]_{\la_2-1}\tot [u(t/q)^{-2}]_{\la_3-2}\tot \cdots \tot [u(t/q)^{-N+1}]_{\la_N-N+1} 
\en
and define $W^{(N)}(u)$ to be a subspace of $V^{(N)}(u)$ spanned by $\{\ket{\la}^{(N)}_u\ |\ \la\in {\cal P}^{(N)}\}$. 
An action of $\cU_{q,t,p}$ on the tensor product space can be constructed 
by using the  comultiplication 
$\Delta^{op}$ 
 in \eqref{opco1}-\eqref{opco3} repeatedly. 
 
One can verify the following two propositions.  
\begin{prop}
By applying  $\Delta^{op}$ repeatedly, the following gives a level-(0,0) $\cU_{q,t,p}$-module structure on $W^{(N)}(u)$.
 \bea
&&\gamma^{1/2}\ket{\la}^{(N)}_u=\ket{\la}^{(N)}_u,\\
&&x^+(z)\ket{\la}^{(N)}_u=a^+(p)\sum_{i=1}^{N}{A}^{(N)+}_{\la,i}(p)\delta(u_i/z)\ket{\la+{\bf 1}_i}^{(N)}_u,\lb{actxpn}\\
&&x^-(z)\ket{\la}^{(N)}_u=
a^-(p)\sum_{i=1}^{N}{A}^{(N)-}_{\la,i}(p)\delta(q^{-1}u_i/z)\ket{\la-{\bf 1}_i}^{(N)}_u,\lb{actxmn}\\
&&\psi^+(z)\ket{\la}^{(N)}_u=
B^{(N)+}_{\la}(u/z;p)\ket{\la}^{(N)}_u,\lb{actpsipn}\\
&&\psi^-(z)\ket{\la}^{(N)}_u=
B^{(N)-}_{\la}(z/u;p)\ket{\la}^{(N)}_u,
\ena
where $u_i=q^{\la_i}t^{-i+1}u$, $\la\pm{\bf 1}_i=(\la_1,\cdots,\la_i\pm1,\cdots,\la_N)$ and we set 
\be
{A}^{(N)+}_{\la,i}(p)&=&
\prod_{j=1}^{i-1}\frac{\theta_p(tu_i/u_j)\theta_p(qt^{-1}u_i/u_j)}{\theta_p(qu_i/u_j)\theta_p(u_i/u_j)}\nn\\
{A}^{(N)-}_{\la,i}(p)&=&
 \prod_{j=i+1}^{N}\frac{\theta_p(tu_j/u_i)\theta_p(qt^{-1}u_j/u_i)}{\theta_p(qu_j/u_i)\theta_p(u_j/u_i)}\\
 {B}^{(N)+}_{\la}(u/z;p)&=&
\prod_{j=1}^{N}\frac{\theta_p(t^{-1}u_j/z)\theta_p(q^{-1}tu_j/z)}{\theta_p(u_j/z)\theta_p(q^{-1}u_j/z)}\nn\\
 {B}^{(N)-}_{\la}(z/u;p)&=&
\prod_{j=1}^{N}\frac{\theta_p(tz/u_j)\theta_p(qt^{-1}z/u_j)}{\theta_p(z/u_j)\theta_p(qz/u_j)}\nn.
 \en
 In particular, we have
\bea
&&\alpha_m\ket{\la}^{(N)}_u=\frac{(1-t^{-m})(1-(q/t)^{-m})}{m}\sum_{j=1}^N u_j^m\ket{\la}^{(N)}_u
\qquad (m\in \Z\backslash \{0\}).
\ena
\end{prop}
\noindent

An inductive limit $N\to \infty$ can be taken  in the same way as in the trigonometric case\cite{FFJMM}. 
Let 
\be
{\cal P}^{(N),+}=\{ \la\in {\cal P}^{(N)}\ |\ \la_N\geq 0\ \}.
\en
and define $W^{(N),+}(u)$ to be the subspace of $W^{(N)}(u)$ spanned by $\{\ket{\la}^{(N)}_u,\ \la\in {\cal P}^{(N),+} \}$.
Let us define $\tau_N:{\cal P}^{(N),+}\to {\cal P}^{(N+1),+}$ by
\be
&&\tau_N(\la)=(\la_1,\la_2,\cdots,\la_N,0).
\en
This induces the embedding $W^{(N),+}\hookrightarrow W^{(N+1),+}$. We then define a semi-infinite tensor product space ${\cal F}_u$ by  the inductive limit 
\be
&&{\cal F}_u=\lim_{N\to \infty} W^{(N),+}(u).
\en
%
The action of $\cU_{q,t,p}$ on $\F_u$ via $\Delta^{op}$  is defined inductively as follows.
Define 
\be
&&x^{+[N]}(z)=x^+(z),\qquad x^{-[N]}(z)=
(q/t)^{1/2}\frac{\theta_p(q^{-1}t^{-N+1}u/z)}{\theta_p(t^{-N}u/z)}x^-(z),\\
&&\psi^{+[N]}(z)=(q/t)^{1/2}\frac{\theta_p(q^{-1}t^{-N+1}u/z)}{\theta_p(t^{N}u/z)}\psi^{+}(z),\qquad 
\psi^{-[N]}(z)=(q/t)^{-1/2}\frac{\theta_p(qt^{N-1}z/u)}{\theta_p(t^{-N}z/u)}\psi^{+}(z).
\en
Then we have
\begin{lem}
For $x=x^\pm, \psi^\pm$, we have 
\be
&&\tau_N(x^{[N]}(z)\ket{\la}_u)=x^{[N+1]}(z)\tau_N(\ket{\la}_u).
\en
\end{lem}
Thanks to this Lemma, one can define the action of $\cU_{q,t,p}$ on  ${\cal F}_u$ by
\be
&&x(z)\ket{\la}_u=\lim_{N\to \infty} x^{[N]}(z)\ket{(\la_1,\la_2,\cdots,\la_N)}. 
\en
This gives the level $(0,1)$ representation in Theorem \ref{actionUqtpgl1}.

\section{Proof of Theorem \ref{typeIVO}}\lb{prooftypeIVO}

\noindent
Proof of  \eqref{Phipsip}: \  
From Theorem \ref{level1N} one has
\be
&&\Phi_{\emptyset}(u)\psi^+((t/q)^{1/4}z)\Bigl|_{\F^{(1,N+1)}}=(t/q)^{-(N+1)/2}\frac{(pqz/tu;p)_\infty}{(pz/u;p)_\infty}:\Phi_{\emptyset}(u)\widetilde{\psi}^+((t/q)^{1/4}z):,\\
&&\psi^+((t/q)^{1/4}z)\Bigl|_{\F^{(1,N)}}\Phi_{\emptyset}(u)=(t/q)^{-N/2}\frac{(u/z;p)_\infty}{(tu/qz;p)_\infty}:\widetilde{\psi}^+((t/q)^{1/4}z)\Phi_{\emptyset}(u):.
\en
Hence we have
\be
&&\Phi_{\emptyset}(u)\psi^+((t/q)^{1/4}z)=(t/q)^{-1/2}\frac{\theta_p(tu/qz)}{\theta_p(u/z)}\psi^+((t/q)^{1/4}z)\Phi_{\emptyset}(u).
\en
Then \eqref{Phipsip} follows from 
\be
&&:\widetilde{\psi}^+((t/q)^{1/4}z)^{-1}\tPhi_\lambda(u)\widetilde{\psi}^+((t/q)^{1/4}z):\nn\\
&&=:\widetilde{\psi}^+((t/q)^{1/4}z)^{-1}\Phi_\emptyset(u)\widetilde{\psi}^+((t/q)^{1/4}z)\prod_{(i,j)\in \lambda}
\widetilde{\psi}^+((t/q)^{1/4}z)^{-1}\widetilde{x}^-((t/q)^{1/4}u_{i,j}
)\widetilde{\psi}^+((t/q)^{1/4}z):,
\en
where $u_{i,j}=q^{j-1}t^{-i+1}u$, and the following Proposition. 
\begin{prop}\lb{prodg=B}
\be
&&\prod_{i=1}^{\ell(\la)}\prod_{j=1}^{\la_i}g(q^{j-1}t^{-i+1}z;p)=\frac{\theta_p(z)}{\theta_p((t/q)z)}B^+_\la(z;p).
\en
\end{prop}
\qed

\noindent
Proof of  \eqref{Phixp}:\  From Theorem \ref{level1N} one has
\be
&&\Phi_\emptyset(u)\widetilde{x}^+((t/q)^{-1/4}z)=(1-qz/tu):\Phi_\emptyset(u)\widetilde{x}^+((t/q)^{-1/4}z):,\\
&&\widetilde{x}^+((t/q)^{-1/4}z)\Phi_\emptyset(u)=(1-tu/qz):\Phi_\emptyset(u)\widetilde{x}^+((t/q)^{-1/4}z):,\\
&&\widetilde{x}^+((t/q)^{-1/4}z)\widetilde{x}^-((t/q)^{1/4}w)=\frac{(1-tw/z)(1-w/qz)}{(1-tw/qz)(1-w/z)}:\widetilde{x}^+((t/q)^{-1/4}z)\widetilde{x}^-((t/q)^{1/4}w):
\en
Hence we have
\be
&&\tPhi_\la(u)\widetilde{x}^+((t/q)^{-1/4}z)=(1-t^{\ell(\la)}qz/tu)\prod_{j=1}^{\ell(\la)}\frac{1-qz/tu_j}{1-qz/u_j}:\tPhi_\la(u)\widetilde{x}^+((t/q)^{-1/4}z):\nn\\
&&\hspace{13cm}\mbox{for}\ |z/u_j|<1,\\
&&\widetilde{x}^+((t/q)^{-1/4}z)\tPhi_\la(u)=(1-t^{-{\ell(\la)}}tu/qz)\prod_{j=1}^{\ell(\la)}\frac{1-tu_j/qz}{1-u_j/qz}:\tPhi_\la(u)\widetilde{x}^+((t/q)^{-1/4}z):\nn\\
&&\hspace{13cm}\mbox{for}\ |u_j/z|<1
\en
and
\bea
&&:\tPhi_\la(u)\widetilde{x}^+((t/q)^{-1/4}q^{-1}u_k):\nn\\
&&=:\Phi_\emptyset(u)\prod_{(i,j)\in \la}\widetilde{x}^-((t/q)^{1/4}u_{i,j})\widetilde{x}^+((t/q)^{-1/4}q^{-1}u_k):\nn\\
&&=:\Phi_\emptyset(u)\prod_{(i,j)\in \la-{\bf1}_k}\widetilde{x}^-((t/q)^{1/4}u_{i,j})\cdot\widetilde{x}^-((t/q)^{1/4}u_{k,\la_k})\widetilde{x}^+((t/q)^{-1/4}q^{-1}u_k):\nn\\
&&=:\tPhi_{\la-{\bf1}_k}(u)\widetilde{\psi}^-((t/q)^{-1/4}q^{-1}u_k):.\lb{Philaxp}
\ena
The last equality follows from 
\be
&&:\widetilde{x}^-((t/q)^{1/4}z)\widetilde{x}^+((t/q)^{-1/4}z):
=\widetilde{\psi}^-((t/q)^{-1/4}z).
\en
Therefore noting
\be
&&x^+(z)\Bigl|_{\F^{(1,N)}_v}=vz^{-N}(t/q)^{3N/4}\widetilde{x}^+(z),\\
&&x^+(z)\Bigl|_{\F^{(1,N+1)}_{-uv}}=-uvz^{-(N+1)}(t/q)^{3(N+1)/4}\widetilde{x}^+(z),
\en
one obtains
\bea
&&
\left(-uvz^{-(N+1)}(t/q)^{N+1}\right)^{-1}\left(\tPhi_\la(u){x}^+((t/q)^{-1/4}z)\Bigl|_{\F^{(1,N+1)}_{-uv}}-
{x}^+((t/q)^{-1/4}z)\Bigl|_{\F^{(1,N)}_v}\tPhi_\la(u)\right)\nn\\
&&
=\left(
(1-t^{{\ell(\la)}}qz/tu)\prod_{j=1}^{\ell(\la)}\frac{1-qz/tu_j}{1-qz/u_j}\Bigl|_{|z/u_j|<1}+(t/q)^{-1}z/u
(1-t^{-{\ell(\la)}}tu/qz)\prod_{j=1}^{\ell(\la)}\frac{1-tu_j/qz}{1-u_j/qz}\Bigl|_{|u_j/z|<1}
\right)\nn\\
&&
\hspace{7cm}\times:\tPhi_\la(u)\widetilde{x}^+((t/q)^{1/4}(t/q)^{-1/2}z):\nn\\
&&
=\left((1-t^{{\ell(\la)}}qz/tu)\prod_{j=1}^{\ell(\la)}\frac{1-qz/tu_j}{1-qz/u_j}\Bigl|_{|z/u_j|<1}-(1-t^{{\ell(\la)}}qz/tu)\prod_{j=1}^{\ell(\la)}\frac{1-qz/tu_j}{1-qz/u_j}\Bigl|_{|u_j/z|<1}\right)\nn\\
&&\hspace{7cm}\times:\tPhi_\la(u)\widetilde{x}^+((t/q)^{1/4}(t/q)^{-1/2}z):\nn\\
&&
=\sum_{k=1}^{\ell(\la)}\delta(q^{-1}u_k/z)
\prod_{j=1}^{k-1}\frac{1-u_k/tu_j}{1-u_k/u_j}\prod_{j=k+1}^{{\ell(\la)}}\frac{1-u_k/tu_j}{1-u_k/u_j}
(1-u_k/tu_{{\ell(\la)}+1})
:\tPhi_\la(u)\widetilde{x}^+((t/q)^{1/4}(t/q)^{-1/2}z):\nn\\
&&
=\sum_{k=1}^{\ell(\la)}\delta(q^{-1}u_k/z)
(1-t^{-1})(q/t)t^{-k+1}\frac{c_\la}{c_{\la-{\bf1}_k}}A^-_{\la,k}
:\tPhi_{\la-{\bf1}_k}(u)\widetilde{\psi}^-((t/q)^{-1/4}q^{-1}u_k):.\lb{PhixpmxpPhi}
\ena
In the last equality we used \eqref{Philaxp}. 
The third equality follows from the formula
\bea
&&(1-s/b_{n+1})\prod_{j=1}^n\frac{1-s/b_j}{1-s/a_j}\Bigl|_{|s/a_j|<1}-
(1-s/b_{n+1})\prod_{j=1}^n\frac{1-s/b_j}{1-s/a_j}\Bigl|_{|a_j/s|<1}\nn\\
&&=\sum_{k=1}^n\delta(s/a_k)(1-a_k/b_{n+1})\frac{\prod_{j=1}^n(1-a_k/b_j)}{
\prod_{j=1\atop j\not=k}^n(1-a_k/a_j)}.\lb{dprodDelta2}
\ena
This is obtained from the partial fraction formula for $\prod_{j=1}^n\frac{1-s/b_j}{1-s/a_j}$ and 
\be
&&\frac{1}{1-z}\Bigr|_{|z|<1}-\frac{1}{1-z}\Bigr|_{|z|>1}=\delta(z).
\en

On the other hand we have
\be
&&\widetilde{\psi}^-((t/q)^{-1/4}z)\Phi_\emptyset(u)=\frac{(pu/z;p)_\infty}{(ptu/qz;p)_\infty}
:\widetilde{\psi}^-((t/q)^{-1/4}z)\Phi_\emptyset(u):,\\
&&\widetilde{\psi}^-((t/q)^{-1/4}z)\widetilde{x}^-((t/q)^{1/4}w)=
h(w/z)
:\widetilde{\psi}^-((t/q)^{-1/4}z)\widetilde{x}^-((t/q)^{1/4}w):,
\en
where
\be
&&h(w)=\frac{(pqw;p)_\infty(ptw/q;p)_\infty(pw/t;p)_\infty}{(pw/q;p)_\infty(pqw/t;p)_\infty(ptw;p)_\infty}. 
\en
We have
\be
&&\prod_{(i,j)\in \la-{\bf1}_k}h(u_{i,j}/z)\Bigl|_{z=q^{-1}u_k}=-t\frac{a^-(p)}{a^+(p)}
\frac{(ptu/u_k;p)_\infty}{(pqu/u_k;p)_\infty}
\prod_{j=1}^{\ell(\la)}\frac{(pu_j/u_k;p)_\infty}{(pqu_j/tu_k;p)_\infty}
\prod_{j=1}^{{\ell(\la)}+1}\frac{(pqu_j/u_k;p)_\infty}{(ptu_j/u_k;p)_\infty}.
\en
Hence 
\be
&&\widetilde{\psi}^-((t/q)^{-1/4}z)\tPhi_{\la-{\bf1}_k}(u)\Bigl|_{z=q^{-1}u_k}\\
&&=-t\frac{a^-(p)}{a^+(p)}
\prod_{j=1}^{\ell(\la)}\frac{(pu_j/u_k;p)_\infty}{(pqu_j/tu_k;p)_\infty}
\prod_{j=1}^{{\ell(\la)}+1}\frac{(pqu_j/u_k;p)_\infty}{(ptu_j/u_k;p)_\infty}
:\widetilde{\psi}^-((t/q)^{-1/4}z)\tPhi_{\la-{\bf1}_k}(u):\Bigl|_{z=q^{-1}u_k}.
\en
Substituting this into \eqref{PhixpmxpPhi}, one obtains
\be
&&
\tPhi_\la(u){x}^+((t/q)^{-1/4}z)\Bigl|_{\F^{(1,N+1)}_{-uv}}-
{x}^+((t/q)^{-1/4}\Bigl|_{\F^{(1,N)}_v}\tPhi_\la(u)\nn\\
&&=uvz^{-(N+1)}(t/q)^{N+1}\sum_{k=1}^{\ell(\la)}\delta(q^{-1}u_k/z)
t^{-k}a^+(p)\frac{c_\la}{c_{\la-{\bf1}_k}}A^-_{\la,k}(p)\frac{N_{\la-{\bf1}_k}(p)}{N_\la(p)}
(t/q)^{-N/2}\\
&&\hspace{7cm}\times{\psi}^-((t/q)^{-1/4}z)\tPhi_{\la-{\bf1}_k}(u)\nn\\
&&=\frac{c_\la}{q^{n(\la')}N_{\la}(p) t^*(\la,u,v,N)}q\sum_{k=1}^{\ell(\la)}\delta(q^{-1}u_k/z)
t^{-1}a^+(p)A^-_{\la,k}(p){\psi}^-((t/q)^{-1/4}z)\nn\\
&&\hspace{7cm}\times \frac{q^{n((\la-{\bf1}_k)')}N_{\la-{\bf1}_k}(p)t^*(\la-{\bf1}_k,u,v,N)}{c_{\la-1_k}
}\tPhi_{\la-{\bf1}_k}(u).
\en
Here we used
\be
&&\frac{t^*(\la-{\bf1}_k,u,v,N)}{t^*(\la,u,v,N)}=q^{-1}v(q^{-1}u_k)^{-N}(t/q)^{N/2}.
\en
\qed
\noindent

\noindent
Proof of  \eqref{Phixm}:\ 
From Theorem \ref{level1N} we have
\be
&&\tPhi_\lambda(u)\widetilde{x}^-((t/q)^{1/4}z)=\prod_{i=1}^{\ell(\la)}\frac{(tz/u_i;p)_\infty}{(pqz/u_i;p)_\infty}\prod_{i=1}^{\ell(\la)+1}\frac{(pqz/tu_i;p)_\infty}{(z/u_i;p)_\infty}:\widetilde{x}^-((t/q)^{1/4}z)\tPhi_\lambda(u):\quad \mbox{for}\ |z/u|<1,\\
&&\widetilde{x}^-((t/q)^{1/4}z)\tPhi_\lambda(u)=\prod_{i=1}^{\ell(\la)}\frac{(q^{-1}u_i/z;p)_\infty}{(pt^{-1}u_i/z;p)_\infty}\prod_{i=1}^{\ell(\la)+1}\frac{(pu_i/z;p)_\infty}{(tu_i/qz;p)_\infty}:\widetilde{x}^-((t/q)^{1/4}z)\tPhi_\lambda(u):\quad \mbox{for}\ |u/z|<1
\en
and
\be
&&B^+_\la(u/z,p)\widetilde{x}^-((t/q)^{1/4}z)\tPhi_\lambda(u)=-z/u
\prod_{i=1}^{\ell(\la)}\frac{(tz/u_i;p)_\infty}{(pqz/u_i;p)_\infty}\prod_{i=1}^{\ell(\la)+1}\frac{(pqz/tu_i;p)_\infty}{(z/u_i;p)_\infty}:\widetilde{x}^-((t/q)^{1/4}z)\tPhi_\lambda(u):.
\en
Hence noting 
\be
&&x^-(z)\Bigl|_{\F^{(1,N+1)}_{-uv}}=(-uv)^{-1}z^{N+1}(t/q)^{-3(N+1)/4}\widetilde{x}^-(z),\\
&&x^-(z)\Bigl|_{\F^{(1,N)}_{v}}=v^{-1}z^{N}(t/q)^{-3N/4}\widetilde{x}^-(z),
\en
one gets
\bea
&&\tPhi_\lambda(u)x^-((t/q)^{1/4}z)-(t/q)^{-1/2}B^+_\la(z,p)x^-((t/q)^{1/4}z)\tPhi_\lambda(u)\nn\\
&&=-(uv)^{-1}z^{N+1}(t/q)^{-(N+1)/2}\nn\\
&&\times\left(
\prod_{i=1}^{\ell(\la)}\frac{(tz/u_i;p)_\infty}{(pqz/u_i;p)_\infty}\prod_{i=1}^{\ell(\la)+1}\frac{(pqz/tu_i;p)_\infty}{(z/u_i;p)_\infty}\Bigl|_{|z/u|<1}-
\prod_{i=1}^{\ell(\la)}\frac{(tz/u_i;p)_\infty}{(pqz/u_i;p)_\infty}\prod_{i=1}^{\ell(\la)+1}\frac{(pqz/tu_i;p)_\infty}{(z/u_i;p)_\infty}\Bigl|_{|u/z|<1}
\right)\nn\\
&&\hspace{9cm} \times:\tPhi_\lambda(u)\widetilde{x}^-((t/q)^{1/4}z):\nn\\
&&=-(uv)^{-1}z^{N+1}(t/q)^{-(N+1)/2}\nn\\
&&\times\left(
\frac{\prod_{i=1}^{{\ell(\la)}}(1-tz/u_i)}{\prod_{i=1}^{{\ell(\la)}+1}(1-z/u_i)}\Bigl|_{|z/u_i|<1}-
\frac{\prod_{i=1}^{{\ell(\la)}}(1-tz/u_i)}{\prod_{i=1}^{{\ell(\la)}+1}(1-z/u_i)}\Bigl|_{|u_i/z|<1}
\right)
\prod_{i=1}^{{\ell(\la)}}\frac{(ptz/u_i;p)_\infty}{(pqz/u_i;p)_\infty}\prod_{i=1}^{{\ell(\la)}+1}\frac{(pqz/tu_i;p)_\infty}{(pz/u_i;p)_\infty}
\nn\\
&&\hspace{9cm} \times:\tPhi_\lambda(u)\widetilde{x}^-((t/q)^{1/4}z):
. \lb{PhixmmxmPhi}
\ena
Applying  the formula
\bea
&&\frac{\prod_{j=1}^{n-1}(1-s/b_j)}{\prod_{j=1}^{n}(1-s/a_j)}\Bigl|_{|s/a|<1}-
\frac{\prod_{j=1}^{n-1}(1-s/b_j)}{\prod_{j=1}^{n}(1-s/a_j)}\Bigl|_{|a/s|<1}
=\sum_{k=1}^n\delta(a_k/s)
\frac{\prod_{j=1}^{n-1}(1-a_k/b_j)}{\prod_{j=1\atop \not=k}^{n}(1-a_k/a_j)}
\lb{dprodDelta1}
\ena
and using
\be
&&:\tPhi_\lambda(u)\widetilde{x}^-((t/q)^{1/4}z):\Bigl|_{z=u_k}=\tPhi_{\la+{\bf1}_k}(u), 
\en
 the RHS of \eqref{PhixmmxmPhi} is
\be
&&(uv)^{-1}z^{N+1}(t/q)^{-(N+1)/2}\sum_{k=1}^{{\ell(\la)}+1}\delta(u_k/z)
\prod_{i=1}^{{\ell(\la)}}\frac{(tu_k/u_i;p)_\infty}{(pqu_k/u_i;p)_\infty}
\frac{\prod_{i=1}^{{\ell(\la)}+1}(pqu_k/tu_i;p)_\infty}{\prod_{i=1\atop \not=k}^{{\ell(\la)}+1}(u_k/u_i;p)_\infty}\frac{1}{(p;p)_\infty}
\tPhi_{\la+{\bf1}_k}(u)\nn\\
&&=(uv)^{-1}z^{N+1}(t/q)^{-(N+1)/2}\sum_{k=1}^{{\ell(\la)}+1}\delta(u_k/z)t^ka^-(p)\frac{c_\la}{c_{\la+{\bf1}_k}}{A}^+_{\la,k}(p)\frac{N_{\la+{\bf1}_k}(p)}{N_\la(p) }\tPhi_{\la+{\bf1}_k}(u)\nn\\
&&=\frac{c_\la}{q^{n(\la')}N_\la(p) t^*(\la,u,v,N)}q^{-1}f(1;p)a^+(p)(t/q)^{-1/2}\\
&&\qquad\quad\times\sum_{k=1}^{{\ell(\la)}+1}\delta(u_k/z){A}^+_{\la,k}(p)\frac{q^{n((\la+{\bf1}_k)')}N_{\la+{\bf1}_k}(p)t^*(\la+{\bf1}_k,u,v,N)}{c_{\la+{\bf1}_k}}
\tPhi_{\la+{\bf1}_k}(u).
\en

\qed

 \section{Proof of Proposition \ref{typeIdual}}\lb{prooftypeIdualVO}

\noindent
Proof of \eqref{Phidualpsip}:\ From \eqref{Phipsip}, we have
\be
&&\psi^+((t/q)^{1/4}z)\tPhi_\la(p^{-1}u)^{-1}=(q/t)^{1/2}B^+_\la(p^{-1}u/z;p)\tPhi_\la(p^{-1}u)^{-1}\psi^+((t/q)^{1/4}z).
\en
\qed

\noindent
Proof of \eqref{Phidualxp}:\ Similar to the calculations in Appendix \ref{prooftypeIVO}, we have 
\be
&&\tPhi_\la(p^{-1}u)^{-1}\widetilde{x}^+((t/q)^{-1/4}z)=\frac{1}{1-pt^lqz/tu}\prod_{j=1}^l\frac{1-pqz/u_j}{1-pqz/tu_j}:\tPhi_\la(p^{-1}u)^{-1}\widetilde{x}^+((t/q)^{-1/4}z):\nn\\
&&\hspace{13cm}\mbox{for}\ |z/u_j|<1,\\
&&\widetilde{x}^+((t/q)^{-1/4}z)\tPhi_\la(p^{-1}u)^{-1}=\frac{1}{1-p^{-1}t^{-l}tu/qz}\prod_{j=1}^l\frac{1-p^{-1}u_j/qz}{1-p^{-1}tu_j/qz}:\tPhi_\la(p^{-1}u)^{-1}\widetilde{x}^+((t/q)^{-1/4}z):\nn\\
&&\hspace{13cm}\mbox{for}\ |u_j/z|<1
\en
and
\bea
&&:\tPhi_\la(p^{-1}u)^{-1}\widetilde{x}^+((t/q)^{-1/4}p^{-1}tu_k/q):\nn\\
&&=:\Phi_\emptyset(p^{-1}u)^{-1}\prod_{(i,j)\in \la}\widetilde{x}^-((t/q)^{1/4}p^{-1}u_{i,j})^{-1}\widetilde{x}^+((t/q)^{-1/4}tu_k/q):\nn\\
&&=:\Phi_\emptyset(p^{-1}u)^{-1}\prod_{(i,j)\in \la+{\bf1}_k}\widetilde{x}^-((t/q)^{1/4}p^{-1}u_{i,j})^{-1}
\times\widetilde{x}^-((t/q)^{1/4}u_{k,\la_k+1})\widetilde{x}^+((t/q)^{-1/4}p^{-1}tu_k/q):\nn\\
&&=:\tPhi_{\la+{\bf1}_k}(p^{-1}u)\widetilde{\psi}^+((t/q)^{1/4}p^{-1}u_k):.\lb{Phiduallaxp}
\ena
Therefore noting
\be
&&x^+(z)\Bigl|_{\F^{(1,N)}_v}=vz^{-N}(t/q)^{3N/4}\widetilde{x}^+(z),\\
&&x^+(z)\Bigl|_{\F^{(1,N+1)}_{-p^{-1}uv}}=-p^{-1}uvz^{-(N+1)}(t/q)^{3(N+1)/4}\widetilde{x}^+(z),
\en
one obtains
\bea
&&
{x}^+((t/q)^{-1/4}z)\Bigl|_{\F^{(1,N+1)}_{-p^{-1}uv}}\tPhi_\la(p^{-1}u)^{-1}-
\tPhi_\la(p^{-1}u)^{-1}{x}^+((t/q)^{-1/4}z)\Bigl|_{\F^{(1,N)}_v}\nn\\
&&\hspace{0cm}=-p^{-1}uvz^{-(N+1)}(t/q)^{N+1}\\
&&\times\left(\frac{1}{1-p^{-1}t^{-l}tu/qz}\prod_{j=1}^{\ell}\frac{1-p^{-1}u_j/qz}{1-p^{-1}tu_j/qz}\Bigl|_{|u_j/z|<1}+p(t/q)^{-1}z/u
\frac{1}{1-pt^{\ell}qz/tu}\prod_{j=1}^{\ell}\frac{1-pqz/u_j}{1-pqz/tu_j}\Bigl|_{|z/u_j|<1}\right)\nn\\
&&\hspace{7cm}\times:\tPhi_\la(p^{-1}u)^{-1}\widetilde{x}^+((t/q)^{-1/4}z):\nn\\
&&\hspace{0cm}=-p^{-1}uvz^{-(N+1)}(t/q)^{N+1}pqz/tu\\
&&\times\left(-\frac{1}{1-t^{\ell}qz/tu}\prod_{j=1}^\ell\frac{1-pqz/u_j}{1-pqz/tu_j}\Bigl|_{|u_j/z|<1}
+\frac{1}{1-pt^{\ell}qz/tu}\prod_{j=1}^\ell\frac{1-pqz/u_j}{1-pqz/tu_j}\Bigl|_{|z/u_j|<1}
\right)\nn\\
&&\hspace{7cm}\times:\tPhi_\la(p^{-1}u)^{-1}\widetilde{x}^+((t/q)^{-1/4}z):\nn\\
&&\hspace{0cm}=-p^{-1}uvz^{-(N+1)}(t/q)^{N+1}pqz/tu\\
&&\times\sum_{k=1}^{\ell+1}\delta(p^{-1}tu_k/qz)
\frac{\prod_{j=1}^{\ell}(1-tu_k/u_j)}{\prod_{j=1\atop \not=k}^{\ell+1}(1-u_k/u_j)}
:\tPhi_{\la+{\bf 1}_k}(p^{-1}u)^{-1}\widetilde{\psi}^-((t/q)^{-1/4}u_k):.
\lb{PhidualxpmxpPhidual}
\ena
The last equality follows from \eqref{dprodDelta1}.

On the other hand we have
\be
&&\tPhi_\emptyset(u)^{-1}\widetilde{\psi}^-((t/q)^{-1/4}z)=\frac{(z/u;p)_\infty}{(qz/tu;p)_\infty}
:\widetilde{\psi}^-((t/q)^{-1/4}z)\tPhi_\emptyset(u)^{-1}:,\\
&&\widetilde{x}^-((t/q)^{1/4}u)^{-1}\widetilde{\psi}^-((t/q)^{-1/4}z)=
h(p^{-1}z/u)^{-1}
:\widetilde{\psi}^-((t/q)^{-1/4}z)\widetilde{x}^-((t/q)^{1/4}u)^{-1}:,
\en
where
\be
&&h(w)=\frac{(pqw;p)_\infty(ptw/q;p)_\infty(pw/t;p)_\infty}{(pw/q;p)_\infty(pqw/t;p)_\infty(ptw;p)_\infty}. 
\en
We have
\be
&&\tPhi_{\la+\bo_k}(p^{-1}u)^{-1}\widetilde{\psi}^-((t/q)^{-1/4}u_k)\\
&&\qquad=\frac{(pu_k/u;p)_\infty}{(pqu_k/tu;p)_\infty}\prod_{(i,j)\in\la+\bo_k}h(u_k/u_{i,j})^{-1}
:\tPhi_{\la+\bo_i}(p^{-1}u)^{-1}\widetilde{\psi}^-((t/q)^{-1/4}u_k):
\en
We have $u_{k,(\la+{\bf1}_k)_k}=u_k$, 
\be
\prod_{(i,j)\in \la+{\bf1}_k}h(u_k/u_{i,j})^{-1}
&=&\prod_{(i,j)\in \la}h(u_k/u_{i,j})^{-1}\times h(u_k/u_{k,(\la+{\bf1}_k)_k})^{-1}\\
&=&-t\frac{a^-(p)}{a^+(p)}
\frac{(pqu_k/tu;p)_\infty}{(pu_k/u;p)_\infty}
\prod_{j=1}^l\frac{(pqu_k/u_j;p)_\infty}{(ptu_k/u_j;p)_\infty}
\prod_{j=1}^{l+1}\frac{(pu_k/u_j;p)_\infty}{(pqu_k/tu_j;p)_\infty}.
\en
Hence 
\be
&&\tPhi_{\la+\bo_k}(p^{-1}u)^{-1}\widetilde{\psi}^-((t/q)^{-1/4}u_k)\\
&&=-t\frac{a^-(p)}{a^+(p)}
\prod_{j=1}^l\frac{(pqu_k/u_j;p)_\infty}{(ptu_k/u_j;p)_\infty}
\prod_{j=1}^{l+1}\frac{(pu_k/u_j;p)_\infty}{(pqu_k/tu_j;p)_\infty}.
:\tPhi_{\la+\bo_k}(p^{-1}u)^{-1}\widetilde{\psi}^-((t/q)^{-1/4}u_k):.
\en
Therefore
\be
&&\frac{\prod_{j=1}^{\ell}(1-tu_k/u_j)}{\prod_{j=1\atop \not=k}^{\ell+1}(1-u_k/u_j)}
:\tPhi_{\la+\bo_k}(p^{-1}u)^{-1}\widetilde{\psi}^-((t/q)^{-1/4}u_k):\\
&&=t^{k-1}a^+(p)\frac{c'_\la N'_{\la+\bo_k}(p)}{c'_{\la+\bo_k}N'_\la(p)}\tA^{+'}_{\la,k}(p)
\tPhi_{\la+\bo_k}(p^{-1}u)^{-1}\widetilde{\psi}^-((t/q)^{-1/4}u_k).
\en
Substituting this into \eqref{PhixpmxpPhi}, one obtains
\be
&&
{x}^+((t/q)^{-1/4}z)\Bigl|_{\F^{(1,N+1)}_{-p^{-1}uv}}\tPhi_\la(p^{-1}u)^{-1}-
\tPhi_\la(p^{-1}u)^{-1}{x}^+((t/q)^{-1/4}z)\Bigl|_{\F^{(1,N)}_v}\nn\\
&&\hspace{0cm}=\frac{c'_\la}{q^{n(\la')}t(\la,v,p^{-1}u,N)N'_\la(p)}(t/q)^{-1/2}a^+(p)\\
&&\times \sum_{k=1}^{\ell+1}\delta(p^{-1}tu_k/qz)\tA^{+'}_{\la,k}(p)
\frac{q^{n((\la+\bo_k)')}t(\la+\bo_k,v,p^{-1}u,N)N'_{\la+\bo_k}(p)}{c'_{\la+\bo_k}}
\tPhi_{\la+{\bf 1}_k}(p^{-1}u)^{-1}{\psi}^+((t/q)^{1/4}qz/t).
\en
Here we used
\be
&&\frac{t(\la+{\bf1}_k,v,p^{-1}u,N)}{t(\la,v,p^{-1}u,N)}=-p^{-1}uvz^{-(N+1)}(t/q)^{3(N+1)/2},\\
&&pqz/tu\cdot t^{k-1}=q^{\la_k}=\frac{q^{n((\la+\bo_k)')}}{q^{n(\la')}}.
\en
\qed
\noindent

\noindent
Proof of \eqref{Phisxm}:\ 
We have 
\be
&&\tPhi_\lambda(u)^{-1}\widetilde{x}^-((t/q)^{1/4}z)=\prod_{i=1}^{l(\la)}\frac{(pqz/u_i;p)_\infty}{(tz/u_i;p)_\infty}\prod_{i=1}^{l(\la)+1}\frac{(z/u_i;p)_\infty}{(pqz/tu_i;p)_\infty}:\widetilde{x}^-((t/q)^{1/4}z)\tPhi_\lambda(u)^{-1}:\quad \mbox{for}\ |z/u|<1,\\
&&\widetilde{x}^-((t/q)^{1/4}z)\tPhi_\lambda(u)^{-1}=\prod_{i=1}^{l(\la)}\frac{(pu_i/tz;p)_\infty}{(u_i/qz;p)_\infty}\prod_{i=1}^{l(\la)+1}\frac{(tu_i/qz;p)_\infty}{(pu_i/z;p)_\infty}:\widetilde{x}^-((t/q)^{1/4}z)\tPhi_\lambda(u)^{-1}:\quad \mbox{for}\ |u/z|<1
\en
and
\be
&&B^+_\la(p^{-1}u/z,p)\tPhi_\lambda(p^{-1}u)^{-1}\widetilde{x}^-((t/q)^{1/4}z)\\
&&=-pz/u
\prod_{i=1}^{l(\la)}\frac{(u_i/tz;p)_\infty}{(p^{-1}u_i/qz;p)_\infty}
\prod_{i=1}^{l(\la)+1}\frac{(p^{-1}tu_i/qz;p)_\infty}{(u_i/z;p)_\infty}:\tPhi_\lambda(p^{-1}u)^{-1}
\widetilde{x}^-((t/q)^{1/4}z):.
\en
Noting 
\be
&&x^-((t/q)^{1/4}z)\Bigl|_{\F^{(1,N+1)}_{-p^{-1}uv}}=(-p^{-1}uv)^{-1}z^{N+1}(t/q)^{-(N+1)/2}\widetilde{x}^-((t/q)^{1/4}z),\\
&&x^-((t/q)^{1/4}z)\Bigl|_{\F^{(1,N)}_{v}}=v^{-1}z^{N}(t/q)^{-N/2}\widetilde{x}^-((t/q)^{1/4}z),
\en
one gets
\bea
&&x^-((t/q)^{1/4}z)\Bigl|_{\F^{(1,N+1)}_{-p^{-1}uv}}\tPhi_\lambda(p^{-1}u)^{-1}-(t/q)^{-1/2}B^+_\la(p^{-1}u/z,p)\tPhi_\lambda(p^{-1}u)^{-1}x^-((t/q)^{1/4}z)\Bigl|_{\F^{(1,N)}_{v}}\nn\\
&&=-(uv)^{-1}pz^{N+1}(t/q)^{-(N+1)/2}\nn\\
&&\times\left(
\prod_{i=1}^{l(\la)}\frac{(u_i/tz;p)_\infty}{(p^{-1}u_i/qz;p)_\infty}
\prod_{i=1}^{l(\la)+1}\frac{(p^{-1}tu_i/qz;p)_\infty}{(u_i/z;p)_\infty}\Bigl|_{|u/z|<1}-
\prod_{i=1}^{l(\la)}\frac{(u_i/tz;p)_\infty}{(p^{-1}u_i/qz;p)_\infty}
\prod_{i=1}^{l(\la)+1}\frac{(p^{-1}tu_i/qz;p)_\infty}{(u_i/z;p)_\infty}\Bigl|_{|z/u|<1}
\right)\nn\\
&&\hspace{9cm} \times:\tPhi_\lambda(p^{-1}u)^{-1}\widetilde{x}^-((t/q)^{1/4}z):\nn\\
&&=-(uv)^{-1}pz^{N+1}(t/q)^{-(N+1)/2}\nn\\
&&\times\left(
(1-p^{-1}tu_{\ell+1}/qz)\prod_{i=1}^{l(\la)}\frac{(1-p^{-1}tu_i/qz)}{(1-p^{-1}u_i/qz)}\Bigl|_{|u_i/z|<1}-
(1-p^{-1}tu_{\ell+1}/qz)\prod_{i=1}^{l(\la)}\frac{(1-p^{-1}tu_i/qz)}{(1-p^{-1}u_i/qz)}\Bigl|_{|z/u_i|<1}
\right)\nn\\
&&\qquad\times 
\prod_{i=1}^{l(\la)}\frac{(u_i/tz;p)_\infty}{(u_i/qz;p)_\infty}\prod_{i=1}^{l(\la)+1}\frac{(tu_i/qz;p)_\infty}{(u_i/z;p)_\infty}
:\tPhi_\lambda(p^{-1}u)^{-1}\widetilde{x}^-((t/q)^{1/4}z):
. \lb{PhixsmmxmPhi}
\ena
Applying  the formula \eqref{dprodDelta2} 
and using
\be
&&:\tPhi_\lambda(p^{-1}u)^{-1}\widetilde{x}^-((t/q)^{1/4}z):\Bigl|_{z=p^{-1}u_{k,\la_k}}
=\tPhi_{\la-{\bf1}_k}(p^{-1}u)^{-1}, 
\en
 the RHS of \eqref{PhixsmmxmPhi} is
\be
&&-(uv)^{-1}pz^{N+1}(t/q)^{-(N+1)/2}ta^-(p)\\
&&\qquad\times\sum_{k=1}^{\ell(\la)}\delta(pqz/u_k)
\prod_{i=1\atop \not=k}^{\ell(\la)}\frac{(pqu_i/tu_k;p)_\infty}{(u_i/u_k;p)_\infty}
\prod_{i=1\atop \not=k}^{l(\la)+1}\frac{(tu_i/u_k;p)_\infty}{(pqu_i/u_k;p)_\infty}
\times\tPhi_{\la-{\bf1}_k}(p^{-1}u)^{-1}\nn\\
&&=-(uv)^{-1}pz^{N+1}(t/q)^{-(N+1)/2}
\sum_{k=1}^{\ell(\la)}\delta(pqz/u_k)q^{-\la_k}a^-(p)\frac{c'_\la N'_{\la-{\bf1}_k}(p)}{c'_{\la-{\bf1}_k}N'_\la(p)}A^{-'}_{\la,k}(p)\tPhi_{\la-{\bf1}_k}(p^{-1}u)^{-1}\nn\\
&&=\frac{c'_\la}{q^{n(\la')} t(\la,v,p^{-1}u,N)N'_\la(p)}(t/q)a^-(p)\\
&&\qquad\quad\times\sum_{k=1}^{\ell(\la)}\delta(pqz/u_k)A^{-'}_{\la,k}(p)
\frac{q^{n((\la-{\bf1}_k)')}t(\la+{\bf1}_k,v,p^{-1}u,N)N'_{\la-{\bf1}_k}(p)}{c'_{\la-{\bf1}_k}}
\tPhi_{\la-{\bf1}_k}(p^{-1}u)^{-1}.
\en

\qed

\section{Useful formulas for the Nekrasov function}\lb{NekF}
For readers' convenience we list some useful formulas for the 
Nekrasov function $N_{\la\mu}(x)$.  
Let $\la, \mu\in \cP^+$. Then for any integer $\ell\geq \ell(\la), \ell(\mu)$, one has 
the following formulas\cite{AK,AFS}.  
\be
N_{\la\mu}(x)&=&\prod_{\Box\in \la}(1-xq^{-a_\mu(\Box)-1}t^{-\ell_\la(\Box)})
\prod_{\blacksquare\in \mu}(1-xq^{a_\la(\blacksquare)}t^{\ell_\mu(\blacksquare)+1})\\
&=&\prod_{1\leq i\leq j\leq \ell}\frac{(xq^{-\mu_i+\la_{j+1}}t^{i-j};q)_\infty}{(xq^{-\mu_i+\la_{j}}t^{i-j};q)_\infty}\prod_{1\leq r\leq s\leq \ell}\frac{(xq^{\la_r-\mu_{s}}t^{-r+s+1};q)_\infty}{(xq^{\la_r-\mu_{s+1}}t^{-r+s+1};q)_\infty}\\
&=&\prod_{1\leq i\leq j\leq \ell}{(xq^{-\mu_i+\la_{j+1}}t^{i-j};q)_{\la_j-\la_{j+1}}}\prod_{1\leq r\leq s\leq \ell}{(xq^{\la_r-\mu_{s}}t^{-r+s+1};q)_{\mu_s-\mu_{s+1}}}\\
&=&\prod_{i,j=1}^\ell \frac{(xq^{\la_i-\mu_{j}}t^{-i+j+1};q)_\infty}{(xq^{\la_i-\mu_{j}}t^{-i+j};q)_\infty}\prod_{r=1}^\ell  \frac{(xq^{-\mu_{r}}t^{-\ell+r};q)_\infty}{(xq^{\la_r}t^{-r+\ell+1};q)_\infty},
\en
where
\be
&&(x;q)_n=(1-x)(1-xq)\cdots(1-xq^{n-1}). 
\en
The second equality follows from the formula 
\bea
&&(1-q)\sum_{(i,j)\in \la}q^{\mu_i-j}t^{\la'_j-i+1}=t\sum_{1\leq i\leq \ell(\la)}q^{\mu_i-\la_j}t^{j-i}
-\sum_{1\leq i< \ell(\la)+1}q^{\mu_i-\la_j}t^{j-i},
\ena
which is derived by using\cite{MacBook}
\bea
&&(1-q)\sum_{j=1}^{\la_i}q^{-j}t^{\la_j'}=\sum_{j=i}^{\ell(\la)}q^{-\la_j}t^{j}(1-q^{\la_j-\la_{j+1}}).
\ena

\end{appendix}

\renewcommand{\baselinestretch}{0.7}

\end{document}